\documentclass[]{article}
\usepackage[utf8]{inputenc}
\usepackage[T1]{fontenc}
\usepackage[english]{babel}
\usepackage{amsmath}
\usepackage{amsthm, amssymb}
\usepackage{shadethm}
\usepackage{algpseudocode, algorithm, algorithmicx}
\usepackage{graphicx}
\usepackage{xcolor}
\usepackage{mathrsfs}
\usepackage{geometry}
\usepackage{enumerate}
\usepackage{cleveref}
\usepackage{mathtools}
\usepackage{thmtools}
\usepackage[nothm]{thmbox}
\usepackage{authblk}

\usepackage{subcaption}
\usepackage{soul}

\usepackage{array}
\newcolumntype{L}{>{$}l<{$}}

\newtheorem{theorem}{Theorem}[section]
\newtheorem{corollary}[theorem]{Corollary}
\newtheorem{proposition}[theorem]{Proposition}
\newtheorem{lemma}[theorem]{Lemma}
\newtheorem{definition}[theorem]{Definition}
\newtheorem{remark}[theorem]{Remark}
\newenvironment{Proof}[1][]{\proof[#1]\normalsize}{\endproof}

\newcommand*\Let[2]{\State #1 $\gets$ #2}
\algrenewcommand\algorithmicrequire{\textbf{Given:}}
\algrenewcommand\algorithmicensure{\textbf{Return:}}

\newcommand{\dist}{\mathrm{dist}}
\newcommand{\rank}{\mathrm{rank}}

\newcommand{\supp}{\mathrm{supp}}
\newcommand{\graph}{\mathrm{graph}}
\newcommand{\id}{\mathrm{I}}
\newcommand{\crit}{\mathrm{crit}}
\newcommand{\dom}{\mathrm{dom}}
\newcommand{\trace}{\mathrm{trace}}
\newcommand{\A}{\mathcal{A}}

\renewcommand{\H}{\mathcal{H}}

\newcommand{\St}{\mathbb{S}}

\newcommand{\R}{\mathbb{R}}

\newcommand{\abg}{{\alpha, \beta}}

\newcommand{\eps}{\varepsilon}
\newcommand{\de}{\; d\eps}
\newcommand{\tab}{\;\;\;\;}
\newcommand{\red}[1]{#1} 
\newcommand{\blue}[1]{#1}
\newcommand{\BLUE}[1]{#1}
\newcommand{\purple}[1]{#1}

\newcommand{\ol}[1]{\bar{#1}}
\renewcommand{\overline}[1]{\ol{#1}}
\newcommand{\vect}{\mathrm{vec}}
\newcommand{\E}[2][]{\textnormal{\;\textsf{E}}_{#1}\!\left[#2\right]}

\renewcommand{\Pr}[2][]{\textnormal{\;\textsf{Pr}}_{#1}\!\left[#2\right]}

\DeclareMathOperator*{\argmin}{arg\,min}

\author[1]{Massimo Fornasier}
\author[2]{Johannes Maly}
\author[3]{Valeriya Naumova}
\affil[1]{\small Applied and Numerical Analysis and Optimization and Data Analysis, Technical University Munich, Germany}
\affil[2]{\small Chair for Mathematics of Information Processing, RWTH Aachen University, Germany}
\affil[3]{\small Machine Intelligence Department, Simula Metropolitan Center for Digital Engineering, Norway}
\date{}
\begin{document}




\title{Robust Recovery of Low-Rank Matrices with Non-Orthogonal Sparse Decomposition from Incomplete Measurements}
\maketitle

\begin{abstract}
We consider the problem of  recovering an unknown effectively \red{$(s_1,s_2)$-}sparse low-rank-$R$ matrix $X$ \red{with possibly non-orthogonal  rank-$1$ decomposition} from incomplete and inaccurate linear measurements of the form $y = \mathcal A (X) + \eta$, \red{where $\eta$ is an {\it ineliminable} noise\footnote{With ``ineliminable'' we mean that the focus of the paper is not on recovery from exact $y = \mathcal A (X)$ measurements,  rather on the stable recovery under severe measurement noise.}.}
We \red{first derive} an optimization formulation for matrix recovery under the considered model and propose a novel algorithm, called {\it \textbf{A}lternating \textbf{T}ikhonov regularization and \textbf{Las}so} (A-T-LA$\rm{S}_{2,1}$), to solve it. The algorithm is based on a multi-penalty regularization, which is able to leverage both structures (low-rankness and sparsity) simultaneously. 
The algorithm is a fast first order method, and straightforward to implement. We prove global convergence for {\it any} linear measurement model to stationary points and local convergence to global minimizers. 
By adapting the concept of restricted isometry property from compressed sensing to our novel model class, we prove error bounds between global minimizers and ground truth, up to noise level, from a number of subgaussian measurements scaling as $R(s_1+s_2)$, up to log-factors in the dimension, and relative-to-diameter distortion. Simulation results demonstrate both the accuracy and efficacy of the algorithm, as well as its superiority to the state-of-the-art algorithms in strong noise regimes and for matrices whose singular vectors do not possess exact (joint-) sparse support.
\end{abstract}


\section{Introduction}

Due to data deluge, the existing amounts of data increasingly exceed the processing capacities, leading much of the recent research to focus on compressive data acquisition and storage. In this case, data recovery typically requires finding a solution for an underdetermined system of linear equations, which becomes tractable only when the data possess some special structure.   
This general paradigm encompasses two important problem classes, which have received significant scientific attention recently: compressed sensing and (sparse) principal component analysis.  The work of this paper stands at the intersection of these problems, incorporating the challenges from both of them and extending their framework. Specifically, we are interested in recovering an effectively sparse low-rank matrix $X$ from an incomplete and inaccurate set of linear measurements $y =  \mathcal A (X) + \eta$. This problem is relevant in several areas such as blind deconvolution, machine learning, and data mining \red{\cite{haykin1994blind,zou2006sparse,d2005direct}}, as we illustrate with a couple of motivating examples \red{in Section \ref{sec:Examples} below}.

\paragraph{Related work} \blue{The recovery from linear measurements of low-rank matrices without sparsity constraints has been well-studied as an extension of classical compressed sensing theory, i.e., compressed sensing of sparse vectors \cite{candes2011tight, recht2010} -- a vector is called $s$-sparse if at most $s$ of its entries are nonzero -- from seemingly under-determined linear measurements.} When the unknown matrix is assumed to have both low-rankness and sparsity, Oymak et.\ al.\ in \cite{oymak2015simultaneously} showed that a mere convex combination of regularizers for different sparsity structures  does not allow,  in general, outperform the recovery guarantees of the ``best'' one of them alone. Consequently, in order to improve recovery further, one has to go beyond linear combinations of already known convex regularizers. In \cite{bahmani2016near} the authors overcome the aforementioned limitations of purely convex approaches by assuming a nested structure of the measurement operator $\A$ and applying basic solvers for low-rank resp. row-sparse recovery in two consecutive steps which is an elegant approach but clearly restricts possible choices for $\A$. Lee et.\ al.\ \cite{lee2013near}  propose {and} analyze the so-called {\it Sparse Power Factorization} (SPF), a modified version of Power Factorization (see \cite{jain2013low}) for recovery of low-rank matrices by representing them as product of two \red{orthogonal} matrices $X = UV^T$ and then applying alternating minimization over the (de)composing matrix $U,V$. SPF introduces Hard Thresholding Pursuit to each of the alternating steps to enforce additional sparsity of the columns of $U$ and/or $V$. Lee et. al.\ were able to show that using suitable initializations and assuming the noise level to be small enough, SPF approximates low-rank and row- and/or column-sparse matrices $X$ from a nearly optimal number of measurements: \blue{if $X \in \mathbb{C}^{n_1\times n_2}$ is rank-$R$, has $s_1$-sparse columns and $s_2$-sparse rows, $m \gtrsim R(s_1 + s_2) \log (\max \{en_1/s_1, en_2/s_2\})$ measurements suffice for robust recovery ($e$ is here the base of the natural logarithm),  which is up to the log-factor at the information theoretical bound.} \red{Despite the theoretical optimality, the setting of SPF is actually quite restrictive as all columns (resp.\ rows) need to share a common support, and this seems to be an empirically necessary requirement, see Section \ref{Numerics3}, and the matrices $U,V$ need to be simultaneously orthogonal}. On the one hand, empirically, it has also been shown in  \cite{lee2013near} that SPF outperforms methods based on convex relaxation. On the other side, SPF is heavily based on the assumption that the operator $\mathcal A$ possesses a suitable restricted isometry property and cannot be applied to arbitrary inverse problems of type \eqref{eq:measModel}, as it may even fail to converge otherwise. The reason is that  SPF is based on hard-thresholding  \cite{BLUMENSATH2009265}, which is not a Lipschitz continuous (non-expansive) map.

\red{Another related line of work comes from statistical literature under the name sparse principal component analysis (SPCA) \cite{zou2006sparse,d2005direct}. SPCA estimates the principal subspaces of a covariance matrix when the singular vectors are sparse in order to defeat the curse of dimensionality. However, observations in SPCA are 
provided directly from noisy samples, whereas in our case the matrix is only observed indirectly, through linear measurements that mix the components. Therefore, the problem considered in this paper is, in general, much harder than SPCA.
}

\paragraph{Contribution} Recent works provide theoretical and numerical evidence of superior performance of multi-penalty regularization, see \cite{naumova2014minimization, naumova2016, daubechies2016} and references therein, for correct modeling and separation \red{of the additive superposition of signals $u+v$.
Motivated by these results, we extend them to a multiplicative superposition model $u v^T$ and propose to recover and} decompose $X$ by a variational approach based on alternating minimization of  the following multi-penalty functional $J_\abg^R \colon \R^{n_1} \times ... \times \R^{n_1} \times \R^{n_2} \times ... \times \R^{n_2} \rightarrow \R$ defined, for $\alpha, \beta > 0$, by
\begin{equation} \label{Jab0}
J_{\alpha,\beta}^R(u^1, \dots, u^R, v^1, \dots, v^R) := \left\Vert y - \mathcal A\left(\sum_{r=1}^R u^r (v^r)^T\right) 
\right\Vert_2^2 + \alpha \sum_{r=1}^R \Vert u^r \Vert_2^2 + \beta \sum_{r=1}^R \Vert v^r  
\Vert_1, 
\end{equation}
where $\alpha, \beta$ are regularization parameters. We  denote the global minimizer of \eqref{Jab0} by $$(u_{\alpha, \beta}^1,\dots, u_{\alpha, \beta}^R,v_{\alpha, \beta}^1,\dots,v_{\alpha, \beta}^R).$$ 
\red{The main contributions of this paper are summarized in the following list of highlights:
\begin{itemize}
    \item we address the theoretical and numerical analysis of an iterative alternating minimization algorithm, based on simple iterative soft-thresholding,  for the minimisation of \eqref{Jab0}, which we dub {(A-T-LA$\rm S_{2,1}$)} for \textbf{A}lternating \textbf{T}ikhonov regularization and \textbf{Las}so or simply ATLAS throughout the rest of the paper. ATLAS performs the robust recovery of effectively sparse low-rank matrices from incomplete measurements;
    \item  we provide general convergence  guarantees,  in particular  Theorem \ref{LocalConvergence},  for {\it any} inverse problem, where neither restricted isometry property (RIP) of the measurement operator $\mathcal A$  nor conditions on the support distribution of $X$ are assumed. This is achieved by using convex relaxation ($\ell_1$-norm minimization) at each iteration and by  virtue of the Lipschitz continuity of soft-thresholding.  While we show always global convergence to stationary points, we are able to show local convergence to global minimizers within a very specifically determined radius, cf.\ Remark \ref{rem:Radius}. Let us stress that the state-of-the-art Sparse Power Factorization (SPF)  would not even converge in such a general setting (here we do not yet assume any RIP);
    \item we prove approximation guarantees for  global minimizers of \eqref{Jab0}  with no RIP assumptions on the measurements, see Proposition \ref{Bound-y2} and Lemma \ref{SparsityControl};
    \item recovery results for low-rank and sparse matrices available in the literature hold for low level of noise or for vanishing noise only, which is not always a realistic assumption in practice (see for instance, the grocery store problem in Section \ref{sec:Examples} below). We show recovery guarantees for global minimizers of \eqref{Jab0} under a new RIP adapted to the novel matrix class we introduce in this paper, see Theorem \ref{ApproxX} and the model classes \eqref{S} and \eqref{K} in Section \ref{sec:RecoveryPropertiesRIP}; 
    \item in fact, differently from SPF, ATLAS works for a significantly wider class of matrices, see again \eqref{S} and \eqref{K}. In particular, ATLAS does not require exact sparsity, common support, or orthogonality of columns (resp. rows). Moreover we prove that our RIP will be fulfilled for sub-Gaussian measurement matrices, while in \cite{lee2013near} only Gaussian measurements have been so far considered;
    \item we demonstrate in the high level noise regime superior empirical performance of ATLAS as compared to the state-of-the-art algorithm SPF for sparse low-rank matrix recovery. 
\end{itemize}
We stress at this point that, to our knowledge there is no other algorithm or nonconvex program as \eqref{Jab0} available in the literature that enjoy all the above mentioned features. Presently, this comes yet with a price: in fact, so far our analysis falls short in proving global convergence to global minimizers and we intend to address this issue in follow up work. We will approach it by more explicitly estimating the radius of convergence of ATLAS, which would result from a more careful inspection of the Kurdyka-Lojasiewicz property, see Section \ref{LocalConvergenceSection}. This approach is significantly different from recent work on initializations \cite{ma2017implicit,li2018rapid} and it will need novel research.  
}

\blue{
\paragraph{Outline} The organization of the paper is as follows. Section \ref{ProblemSection} presents the setting of this paper, clarifies notation, and introduces the algorithm for the minimisation of \eqref{Jab0}. In Sections \ref{sec:Main1}-\ref{LocalConvergenceSection} we give an overview of the main results. The corresponding proofs can be found in Section \ref{Proofs}. \purple{In} Section \ref{Numerics}, \purple{we present} the actual implementation of \red{ATLAS} and \purple{provide numerical experiments, confirming our theoretical results and showing extensive comparisons to SPF.} We conclude in Section \ref{Conclusion} with a discussion on open problems and future work.
}

\subsection{Some Motivating Examples} \label{sec:Examples}

\red{Before moving to the main part of the paper, let us \red{consider}
a couple of motivating \red{examples and applications of the considered model}. \\
The first example views low-rank matrices with sparsity constraints from a machine learning perspective and extends the classical setting of sparse PCA to incomplete linear observations of the data matrix. \\
The second example is a classical problem in signal processing, that is blind deconvolution. By now this problem has been widely explored in the literature \cite{ahmed2014blind,li2018rapid,lee2016blind}\\
In particular, in both these examples we do not expect necessarily that the measurements fulfill an RIP condition.
}

\paragraph{Example 1: Sparse Principal Component Analysis from inaccurate and incomplete linear measurements} Principal Component Analysis (PCA) \cite{jolliffe2011principal} is a classical tool for processing large amounts of data and performing data analysis such as dimensionality reduction and factor extraction. Its scope of application ranges from engineering and technology to social sciences, and biology.\\
PCA and, more generally, matrix completion \cite{Cand2009} has been widely used for recommendation systems as popularised by the so-called Netflix prize problem \cite{Bennett07thenetflix}. We illustrate PCA by considering a simple example of such recommendation system for a grocery store, which has $n_1$ regular customers and $n_2$ products. Let $X \in \R^{n_1\times n_2}$ be a matrix whose components $X_{i,j}$ encode the probability of customer $i$ buying product $j$. 
It is reasonable to assume that there are only $R \ll \min\{n_1,n_2\}$ underlying basic factors like age, income, family size, etc.\ which govern the customer's purchase behavior. For each basic factor $r \in [R] \coloneqq \{1,...,R\}$ one defines two vectors: a vector $u^r \in \R^{n_1}$ of components $u_i^r$ encoding for each user $i \in [n_1]$ how much they are affected by the factor $r$, and a vector $v^r \in \R^{n_2}$ encoding the probability of buying product $j$ if having factor $r$. Then, one can decompose
\begin{align} \label{eq:PCA}
	X \approx UV^T = \sum_{r=1}^{R} u^r (v^r)^T
\end{align}
as the product of two matrices $U \in \R^{n_1\times R}$ and $V \in \R^{n_2\times R}$ with columns $u^r$ and $v^r$. Even if the product $UV^T$ is only approximately $X$, the decomposition into orthogonal principal components $U$ and loadings $V$ is appealing for more interpretability and having less data to store ($\mathcal{O}(\max \{n_1, n_2\}R)$ instead of $\mathcal{O}(n_1n_2)$).\\
However, if we want to understand  which factors mostly affect customer's behaviour, PCA might not be the best option, since principal components are usually a linear combination of all original variables. To further improve interpretability and reduce the number of explicitly used variables,  sparse PCA \cite{zou2006sparse,d2005direct}, which promotes sparsity of the loadings $v^r$ in \eqref{eq:PCA}, has been proposed. Sparse PCA trades orthogonality of the principal components for sparse solutions. In the aforementioned example of the grocery store, it is quite reasonable to assume sparsity of the probability distributions $v^r$, as certain factors normally are more correlated with the probability of purchase of few specific items.

For some applications one may not have access to the complete matrix $X$ but only to a partial indirect information, i.e., one has only $m \ll n_1n_2$ scalars encoding information about $X$. In the example of the grocery store this may model the situation where customers do not all possess a fidelity card, which allows to identify them individually, and the grocery store still wishes to learn the matrix $X$ from aggregated revenues. Each day $d \in [D]$ the store caches in a certain amount of money $y^\ell_d$ corresponding to purchases of a {\it random} subset $T_d \subset [n_1]$ of its customers ($\ell \in \mathbb{N}$ is a fixed index, whose role will soon become clear). If $P^\ell \in \R^{n_2}$ is a vector encoding the prices $P_j^\ell$ of each product $j$ and $\mathcal{P}_{i,d} \subset [n_2]$ is the {\it random}  set of products purchased by customer $i$ on a day $d$, we can express the takings as 
\begin{align*}
	y^\ell_d = \sum_{i \in T_d} \sum_{j \in \mathcal{P}_{i,d}} P_j^\ell.
\end{align*} 
 If we assume that each customer $i$ visits the grocery store with probability $q_i$, we can compute the expected takings as 
\begin{align*}
	\mathbb E_{T_d,\mathcal{P}_{\cdot,d}} \left [ {\sum_{i \in T_d}  \sum_{j \in \mathcal{P}_{i,d}} P_j^\ell} \right ]= \sum_{i=1}^{n_1} q_i \sum_{j=1}^{n_2} X_{i,j} P_j^\ell.
\end{align*}
Choosing $D$ sufficiently large, the law of large numbers guarantees that
\begin{align*}
	\lim_{D\to \infty} \frac{1}{D} \sum_{d=1}^{D} y^\ell_d =\mathbb E_{T_d,\mathcal{P}_{\cdot,d}} \left [ {\sum_{i \in T_d}\sum_{j \in \mathcal{P}_{i,d}} P_j^\ell} \right ],
\end{align*}
in probability and almost surely. Moreover, by Central Limit Theorem, we may model the average takings of $D$ days as 
\begin{align*}
	 \frac{1}{D} \sum_{d=1}^{D} y^\ell_d= \sum_{i=1}^{n_1} q_i \sum_{j=1}^{n_2} X_{i,j} P_j^\ell + \eta_{D}^\ell,
\end{align*}
for a suitable Gaussian noise $\eta_{D}^\ell$. By defining $y^\ell = \frac{1}{D} \sum_{d=1}^{D} y^\ell_d$, we can rewrite the above equation as
\begin{align*}
	y^\ell= \sum_{i=1}^{n_1} \sum_{j=1}^{n_2} (q_i P_j^\ell) X_{i,j} + \eta_{D,l} = \langle A_\ell, X \rangle_F + \eta_{D}^\ell
\end{align*}
where the matrix $A_\ell \in \R^{n_1\times n_2}$ has entries $(q_i P_j^\ell)_{i,j}$, and $\langle \cdot, \cdot \rangle_F$ is the Frobenius scalar product.\\
 Tracking the daily sales over a time period of $m\cdot D$ days and perturbing the prizes in each subperiod $\ell \in [m]$ randomly\footnote{The random fluctuation of prizes is applied by groceries also for rotating promotions on products. Periodic price reductions, or sales,  constitute a widely observed phenomenon in retailing. Sales occur on a regular basis, which  suggests that they are not entirely due to random variations such as shocks to inventory holdings or demand.} would result in $m$ inaccurate linear measurements, where each single measurement is a random average over the entries of $X$ with an ineliminable additive noise $\eta_{D}^\ell$. The whole measurement process can be written as
\begin{align} \label{eq:measModel}
y = \A (X) + \eta
\end{align}
where $\A \colon \R^{n_1\times n_2} \rightarrow \R^m$ is a linear operator defined by the matrices $A_1,...,A_m$ and $\eta =(\eta_{D}^1,\dots,\eta_{D}^m)^T \in \R^m$ models the  noise. 
As we  demonstrate in this paper, it is possible  by means of our resource efficient algorithm to recover a low-rank-$R$ matrix $X$ with effectively $(s_1,s_2)$-sparse non-orthogonal rank-$1$ decomposition, from a number $m \approx R(s_1+ s_2)$ of random noisy measurement. This would offer a plausible solution to the grocery store problem. (We should stress that in general our algorithm will converge to a data fitting solution of low rank and sparse components  for arbitrary linear measurement operators $\A$.)

\paragraph{Example 2: Blind deconvolution in signal processing} In \emph{blind deconvolution} \cite{haykin1994blind} one is interested in recovering two unknown vectors $w$ and $s$ solely from their (cyclic) convolutional product
\begin{align} \label{eq:BlindDeconvolution}
	y = w \ast s + \eta= \left( \sum_{i=1}^{m} w_i s_{(k-i)\, \mathrm{mod}\, m} \right)_{k = 1}^m + \eta,
\end{align}
where $\eta$ is again measurement noise.
In imaging applications, $s$ represents the picture and $w$ -- an unknown blurring kernel \cite{stockham1975blind}. In signal transmission, $s$ is a coded message and $w$ models the properties of the transmission channel \cite{godard1980self}. Independently of the concrete application, problem \eqref{eq:BlindDeconvolution} is highly under-determined and contains ambiguities.\\ 
In \cite{ahmed2014blind} the authors used that by bilinearity of the convolution,  \eqref{eq:BlindDeconvolution} can be represented as a linear map acting on the tensor product $ws^T$, a technique commonly known as \emph{lifting}. They assumed in addition that the channel properties $w$ and the message $s$ are drawn from lower dimensional subspaces and are of the form $w = Bh$ and $s = Cx$ with $h \in \mathbb{C}^{n_1}$ and $x \in \mathbb{C}^{n_2}$ being coefficient vectors encoding channel and message ($B$ and $C$ are suitable transformation matrices). Accordingly, they  re-write \eqref{eq:BlindDeconvolution} as
\begin{align*}
	y = \A (X)+ \eta,
\end{align*}
where the rank-1 matrix $X=hx^T \in \mathbb{C}^{n_1\times n_2}$ has to be recovered from $m$ linear measurements, a quite popular model in compressed sensing literature \cite{Recht:2010fk}. Under suitable assumptions on $\A$ the recovery of $X$ is solved by (convex) nuclear norm minimization.\\
In \emph{blind demixing} \cite{ling2017blind,ling2017regularized,JungKrahmerStoeger2017} or MIMO channel identification \cite{dili01} a receiver gets the overlay of $R$ different convolutions which translates the above mentioned formulation into the recovery of rank-$R$ matrices from linear measurements of the type
\begin{align*}
	y = \A \left( \sum_{r=1}^{R} h_r x_r^T \right)+ \eta.
\end{align*} 
As already mentioned in \cite{JungKrahmerStoeger2017}, one can typically impose extra structure like sparsity on the channel impulse responses $h$ to further reduce the number of measurements $m$. In this case one wants to benefit from exploiting two different structures at the same time, low-rankness and sparsity. \\


\section{Problem Formulation and Notation} \label{ProblemSection}

We recall that the \textit{Singular Value Decomposition (SVD)} of a matrix $Z \in \R^{n_1 \times n_2}$ is given by
\begin{align} \label{SVD}
Z = U\Sigma V^T = \sum_{r = 1}^{\rank(Z)} \sigma_r u^r (v^r)^T,
\end{align}
where $\Sigma$ is a diagonal matrix containing the singular values $\sigma_1 \ge ... \ge \sigma_{\rank(Z)} > 0,$ while $U \in \R^{n_1 \times \rank(Z)}$ and $V \in \R^{n_2 \times \rank(Z)}$ have orthonormal columns which are called left and right singular vectors. 
In the following, we assume $\hat{X}$ is of rank $R>0$ and possesses a decomposition of the form
\begin{align} \label{xSVD}
\hat{X} = \sum_{r=1}^R \hat{u}^r (\hat{v}^r)^T, 
\end{align}
where $\hat{v}^r$ are effectively $s$-sparse, \red{a useful concept introduced by Plan and Vershynin in \cite{Plan2013LP}. 
\begin{definition}[Effectively Sparse Vectors] \label{def:EffectivelySparse}
	Let 
	\begin{align*}
	K_{n,s} = \{ z \in \R^n \colon \Vert z \Vert_2 \le 1 \text{ and } \Vert z \Vert_1 \le \sqrt{s} \}.
	\end{align*}
	\blue{The set of effectively $s$-sparse vectors of dimension $n$ is defined by $\{ z \in \R^n \colon z / \| z \|_2 \in K_{n,s} \}$.}
\end{definition}
\begin{remark}
	\blue{Recall that $z \in \R^n$ is $s$-sparse if it has at most $s$ nonzero entries and note that any $s$-sparse vector is also effectively $s$-sparse.} Effectively sparse vectors are well approximated by sparse vectors as made precise in \cite[Lemma 3.2]{Plan2013LP}. 
\end{remark}
}
We call the vectors $\hat u^r$ (resp. $\hat v^r$) the left (resp. right)  component vectors of $\hat X$. From the context it will be clear to which decomposition they are referred. In fact, \eqref{xSVD} does not need to be the SVD of $\hat X$, although this case is also \purple{covered} by our analysis, \red{as we do not require $\hat{v}^r$ to be mutually orthogonal}.
We focus on decompositions \eqref{xSVD} with effectively sparse right component vectors. Conceptually straight-forward, but perhaps tedious modifications of the arguments lead to similar results in the left-sided and both-sided sparse case (see Section \ref{MainProofSection}). \blue{We mention here that in the case of right-sided $s$-sparsity a natural dimensional setting is $R \ll s \approx n_1 \ll n_2$ (see Remark \ref{rem:dim} for a more detailed discussion).}\\
Furthermore, we are given some linear measurement operator $\A \colon \R^{n_1 \times n_2} \rightarrow \R^m$ and the vector of measurements $y$, which is obtained from $\hat{X}$ by
\begin{align} \label{A}
y = \A (\hat{X}) + \eta = {\frac{1}{\sqrt m}}\begin{pmatrix}
\langle A_1,\hat{X} \rangle_F \\
\vdots \\
\langle A_m,\hat{X} \rangle_F	
\end{pmatrix} + \eta.
\end{align}
The operator $\A$ is completely characterized by the $m$ matrices $A_i \in \R^{n_1 \times n_2}$ and individual measurements correspond to Frobenius products $\langle A_i,\hat{X} \rangle_F = \trace(A_i \hat{X}^T)$.  Noise comes into play by the additive vector $\eta \in \R^m$ of which only the $\ell_2$-norm is known.

	\begin{table}[h]
		\centering
		\blue{
		\begin{tabular}{ || m{4.5cm} | m{10cm} || } 
			\hline
			\multicolumn{2}{ || c || }{Notation} \\ 
			\hline \hline
			$\hat{X} \in \R^{n_1\times n_2}$ & Low-rank and sparse ground-truth \\ 
			$\hat{u}^r \in \R^{n_1}, \hat{v}^r \in \R^{n_2}$ & Left and right components of $\hat{X}$, cf.\ \eqref{xSVD} \\
			$\A \colon \R^{n_1\times n_2} \rightarrow \R^m$ & Measurement operator, cf.\ \eqref{A} \\
			$y \in \R^m$ & Observed measurements \\
			$\eta \in \R^m$ & Additive measurement noise \\
			$J_\abg^R \colon \R^{n_1} \times ... \times \R^{n_2} \rightarrow \R$ & Objective functional depending on $\alpha, \beta > 0$ and $R \in \mathbb{N}$, cf.\ \eqref{Jab0} \\
			$X_\abg \in \R^{n_1\times n_2}$ & Matrix of left and right components $u_\abg^r$ and $v_\abg^r$ \\
			$u_\abg^r \in \R^{n_1}, v_\abg^r \in \R^{n_2}$ & Components of global (non-unique) minimizer of $J_\abg^R$ \\
			\hline
		\end{tabular}
	}
	\end{table}

\paragraph{Notation} {\purple{For a matrix $Z \in \R^{n_1 \times n_2},$ we denote its  transpose by $Z^T.$ A variety of norms are used throughout this paper: $\Vert Z \Vert_p$ is the Schatten-$p$ quasi-norm ($\ell_p$-quasi-norm of the vector of singular values $(\sigma_1,...,\sigma_R)^T$); $\Vert Z \Vert_F$ is the Frobenius norm (the $\ell_2$-norm of the vector of singular values); $\Vert \cdot \Vert_{2 \rightarrow 2}$ is the operator norm of $Z$ (the top singular value). 
Note that for $0<p < 1$, the Schatten-$p$ norm is only a quasi-norm. For $p = 2$, the Schatten norm is equal to the Frobenius norm, whereas the $\infty$-Schatten norm corresponds to the operator norm.}
We use the shorthand notation $[R] = \{ 1,...,R \}$ to write index sets. \red{We denote the index set of the non-zero entries of $v$ as $\supp(v).$}
The relation $a \gtrsim b$ is used to express $a \ge Cb$ for some positive constant $C$, and $a \simeq b$ stands for {$a \gtrsim b$} and {$b \gtrsim a$}.} \\
\blue{If for $\hat X$ of rank $R>0$ the sparse decomposition in \eqref{xSVD} agrees with the SVD for $\|\hat v^r \|_2= \sigma_r$, $r \in [R]$, then for any $0<p< \infty$}
\begin{equation}\label{qneq0}
\|X\|_p^p = \sum_{r=1}^R(\| \hat u^r \|_2 \| \hat v^r \|_2)^p.
\end{equation}
If the decomposition \eqref{xSVD} does not coincide with the SVD of $\hat X$, then $\hat u^1,\dots, \hat u^R$ are anyhow linearly independent and
\begin{eqnarray*}
\|X\|_F^2 &=& \sum_{j=1}^{n_2} \sum_{i=1}^{n_1} \left | \sum_{r=1}^R \hat u^r_i \hat v^r_j \right |^2
= \sum_{j=1}^{n_2}  \left \| \sum_{r=1}^R \frac{\hat u^r}{\| \hat u^r \|_2} \| \hat u^r \|_2 \hat v^r_j \right\|^2_2 \\ 
&\simeq& \sum_{j=1}^{n_2}  \sum_{r=1}^R \| \hat u^r \|_2^2 |\hat v^r_j |^2=\sum_{r=1}^R (\| \hat u^r \|_2 \|\hat v^r \|_2)^2.
\end{eqnarray*}
From this  and the equivalence of $\ell_p$-quasi-norms and Schatten-$p$-quasi-norms for $0<p \leq 2$, one further obtains as a relaxation of \eqref{qneq0}
\begin{equation}\label{qneq}
c_{\hat U}^{-1} R^{p/2-1}  \sum_{r=1}^R(\| \hat u^r \|_2 \| \hat v^r \|_2)^p \leq \|X\|_p^p \leq C_{\hat U} R^{1-p/2} \sum_{r=1}^R(\| \hat u^r \|_2 \| \hat v^r \|_2)^p,
\end{equation}
for positive constants $c_{\hat U}, C_{\hat U}> 0$, which depend on the largest and smallest eigenvalues of the Gramian of the vectors $\hat u^1/\| \hat u^1 \|_2,\dots, \hat u^R/\| \hat u^R \|_2$.  Below we shall use \eqref{qneq} mostly for $p=2/3$.

\paragraph{Recovery algorithm} 
Following promising results on multi-penalty functionals \red{for unmixing problems \cite{naumova2016, naumova2014minimization}, we propose to approximate} $\hat{X}$ by global minimizers of the functional $J_\abg^R$ defined in \eqref{Jab0}, which \purple{combines} one quadratic least-squared error term on the measurements with several convex regularizers applied to vectors (not matrices). Note that $J_\abg^R$ does apply to matrices implicitly by viewing each $2R$-tuple $(u^1, \dots, u^R, v^1, \dots, v^R)$ as the matrix $X=\sum_{r=1}^R u^r(v^r)^T$, and we denote $X_{\alpha, \beta}=\sum_{r=1}^R u^r_{\alpha, \beta}(v^r_{\alpha, \beta})^T$ the one corresponding to a global minimizer $(u_{\alpha, \beta}^1,\dots, u_{\alpha, \beta}^R,v_{\alpha, \beta}^1,\dots,v_{\alpha, \beta}^R)$.
In such a way, instead of combining convex regularizers for sparse and low-rank matrices, we enforce low-rankness by restricting the domain of $J_\abg^R$ properly (the decomposition can only consist of $R$ vector pairs) and promote sparsity by $\ell_1$-norm regularization directly on vectors of the decomposition.

\blue{Despite the  convex multi-penalty regularization term $\alpha \sum_{r=1}^R \Vert u^r \Vert_2^2 + \beta \sum_{r=1}^R \Vert v^r  
\Vert_1$, the functional \eqref{Jab0} is highly non-convex, hence, in view of the negative results on convex approaches to multi-structural recovery in \cite{oymak2015simultaneously} it provides hope for better performances than lifting and convex relaxation.} At the same time, one notices that $J_\abg^R$ becomes convex when all but one $u^r$ and/or $v^r$ are fixed. Hence, we can minimize the functional efficiently by the alternating scheme $(\text{A-T-LAS}_{2,1})$
\begin{equation}
	 \label{ATLAS0}
\begin{cases}
u_{k+1}^1 &=  \argmin_u \left\Vert \left( y - \A \left( \sum_{r=2}^R u_k^r {v_k^r}^T 
\right) \right) - \A(u {v_k^1}^T) \right\Vert^2_2 + \alpha \Vert u \Vert^2_2
+ \frac{1}{2\lambda_k^1} \Vert u - u_k^1 \Vert^2_2, \\
v_{k+1}^1 &=  \argmin_v \left\Vert \left( y - \A \left( \sum_{r=2}^R u_k^r {v_k^r}^T 
\right)	\right) - \A(u_{k+1}^1 v^T) \right\Vert^2 + \beta \Vert v \Vert_1
+ \frac{1}{2\mu_k^1} \Vert v - v_k^1 \Vert^2_2, \\
&\vdots \\
u_{k+1}^R &=  \argmin_u \left\Vert \left( y - \A \left( \sum_{r=1}^{R-1} u_{k+1}^r 
{v_{k+1}^r}^T \right) \right) - \A(u {v_k^R}^T) \right\Vert^2_2 + \alpha \Vert u \Vert^2_2
+ \frac{1}{2\lambda_k^R} \Vert u - u_k^R \Vert^2_2, \\
v_{k+1}^R &=  \argmin_v \left\Vert \left( y - \A \left( \sum_{r=1}^{R-1} u_{k+1}^r 
{v_{k+1}^r}^T \right) \right) - \A(u_{k+1}^R v^T) \right\Vert^2 + \beta \Vert v \Vert_1
+ \frac{1}{2\mu_k^R} \Vert v - v_k^R \Vert^2_2, \\
\end{cases}
\end{equation}

In each iteration above, the  terms $\Vert u - u_k^r \Vert^2_2$ and $\Vert v - v_k^r \Vert^2_2$ are added to provide theoretical convergence guarantees for the sequence $(u_k^1,...,v_k^R)$ with suitable choice of the $2R$ positive sequences of parameters $(\lambda_k^1)_{k \in \mathbb{N}},\dots ,(\lambda_k^R)_{k \in \mathbb{N}}$, $(\mu_k^R)_{k \in \mathbb{N}},\dots ,(\mu_k^R)_{k \in \mathbb{N}} > 0$. 
In practice, ATLAS converges \purple{without those terms.}
As most of the non-convex minimization algorithms, empirical performances of ATLAS likely depends on \blue{a proper initialization $(u_0^1,\dots,v_0^R)$. Setting $(u_0^1,\dots,u_0^R)$ to the leading left and $(v_0^1,\dots,v_0^R)$ to the leading right singular vectors of $\mathcal A^* (y) = \sum_{i=1}^m y_i A_i$, where $\A^*$ denotes the adjoint of $\A$, ensures empirically stable recovery in the experiments (Section \ref{Numerics}).} However, we do not provide any theoretical guarantees for this observation. 



\blue{We are now ready to state the main results of the paper. First, we show in Section \ref{sec:Main1} how minimizers of $J_\abg^R$ yield, under mild assumptions, solutions to the inverse problem \eqref{A}. Second, to explain the reconstruction performance observed in numerical simulations, we introduce in Section \ref{sec:RecoveryPropertiesRIP} a versatile matrix model class and come up with a suitable and novel restricted isometry property (RIP), which captures both low-rankness and sparsity, \red{and} provide in Section \ref{SectionRIP} bounds on a sufficient number of measurements for subgaussian operators to fulfill the RIP with high probability. Finally, local convergence of ATLAS to global minimizers of $J^R_{\alpha,\beta}$ is discussed in Section \ref{LocalConvergenceSection}.}

\section{Properties of Minimizers of $J_\abg^R$} \label{sec:Main1}

Let us begin with some basic properties that minimizers of $J_\abg^R$ have under very general assumptions. For a given minimizer 
$(u_\abg^1,...,u_\abg^R,v_\abg^1,...,v_\abg^R)$ of $J_\abg^R$ we denote
\begin{align} \label{Xab}
X_\abg = U_\abg \Sigma_\abg V_\abg^T = \sum_{r=1}^{R} (\sigma_\abg)_r \frac{u_\abg^r}{\Vert u_\abg^r \Vert_2} \left( \frac{v_\abg^r}{\Vert v_\abg^r \Vert_2} \right)^T
\end{align} 
where $(\sigma_\abg)_r = \Vert u_\abg^r \Vert_2 \Vert v_\abg^r \Vert_2$ for all $r \in [R]$, and $\Sigma_\abg$ is the diagonal matrix defined by the vector $\sigma_\abg$.  The first result bounds  measurement misfit by $X_\abg$.
\begin{proposition}[Measurement misfit] \label{Bound-y2}
	Assume $(u_\abg^1,...,v_\abg^R)$ is a global minimizer of $J_\abg^R$ and $\hat{X}$ is fulfilling the noisy measurements $y~=~\A(\hat{X})~+~\eta$. Then,
	\begin{align}\label{measfit}
	\Vert y - \A(X_\abg) \Vert_2^2 \le \Vert \eta \Vert_2^2 + C_{2,1} \sqrt[3]{\alpha \beta^2} \sum_{r=1}^{R} \left( \Vert \hat{u}^r \Vert_2 \Vert \hat{v}^r \Vert_1 \right)^{\frac{2}{3}},
	\end{align}
	where $C_{2,1}$ is the constant from Lemma \ref{fpq} below.	
\end{proposition}
\blue{The estimate in Proposition \ref{Bound-y2} has perhaps a counterintuitive form. In particular, the exponent $\frac{2}{3}$ comes from the choice of regularizers ($\ell_2$ on the left and $\ell_1$ on the right components), cf.\ proof of Proposition \ref{Bound-y2} and Lemma \ref{fpq}. Note that, for $\hat{X}$ replaced with $\lambda \hat{X}$ and $\lambda > 0$, the same bound can be obtained by replacing $\alpha,\beta$ with $\lambda^{-\frac{2}{3}} \alpha, \lambda^{-\frac{2}{3}} \beta$. Using similar scalings one may control the minimizer's norm as well.}
%
%
\begin{lemma}[Boundedness] \label{Bound-uv2}
	Assume $(u_\abg^1,...,v_\abg^R)$ is a global minimizer of $J_\abg^R$ and $\hat{X}$ is fulfilling the noisy measurements $y~=~\A(\hat{X})~+~\eta$. If $\| y - \A(X_\abg) \|_2 \ge \| \eta \|_2$, we have
	\begin{align} \label{eq:Bounds}
	\begin{split}
	\sum_{r=1}^{R} \Vert u_\abg^r \Vert_2^2 &\le C_{2,1} \sqrt[3]{ \frac{\beta^2}{\alpha^2}} \sum_{r=1}^{R} \left( \Vert \hat{u}^r \Vert_2 \Vert \hat{v}^r \Vert_1 \right)^\frac{2}{3}, \\
	\sum_{r=1}^{R} \Vert v_\abg^r \Vert_1 &\le C_{2,1} \sqrt[3]{ \frac{\alpha}{\beta}} \sum_{r=1}^{R} \left( \Vert \hat{u}^r \Vert_2 \Vert \hat{v}^r \Vert_1 \right)^\frac{2}{3},
	\end{split}
	\end{align}
	and
	\begin{align*}
	\sum_{r=1}^{R} \left( \Vert u_\abg^r \Vert_2 \Vert v_\abg^r \Vert_1 \right)^\frac{2}{3} \le \sum_{r=1}^{R} \left( \Vert \hat{u}^r \Vert_2 \Vert \hat{v}^r \Vert_1 \right)^\frac{2}{3}
	\end{align*}
	where $C_{2,1}$ is the constant from Lemma \ref{fpq}.
\end{lemma}
The two estimates in \eqref{eq:Bounds} point out an interesting property of $J_\abg^R$. If one chooses the parameters $\alpha$ and $\beta$ of different magnitude, either the left or the right components of a minimizer $(u_\abg^1,...,v_\abg^R)$ can be forced to become smaller in norm, while the grip on the others is lost. If $\alpha$ and $\beta$ are chosen to be equal the norm bounds are balanced and one obtains
\begin{align*}
\sum_{r=1}^{R} \left( \Vert u_\abg^r \Vert_2^2 + \Vert v_\abg^r \Vert_1 \right) \le C_{2,1} \sum_{r=1}^{R} \left( \Vert \hat{u}^r \Vert_2 \Vert \hat{v}^r \Vert_1 \right)^\frac{2}{3}.
\end{align*}
The assumption $\| y - \A(X_\abg) \|_2 \ge \| \eta \|_2$ is not restrictive. As soon as $\| y - \A(X_\abg) \|_2 = \| \eta \|_2$ one does not have to decrease $\alpha$ and $\beta$ any further \red{as this will lead to overfitting}. 
\red{We can further control effective sparsity of the minimizer's right components.}
%
\begin{lemma}[Sparsity control] \label{SparsityControl}
	Assume $\A \colon \R^{n_1\times n_2} \rightarrow \R^m$ is a linear operator and $y \in \R^m$. Let $(u_\abg^1,...,v_\abg^R)$ be a minimizer of $J_\abg^R$. For all $r \in [R]$ we have that if $\Vert v_\abg^r \Vert_2 \ge \Vert y \Vert_2^2/\gamma$ for some $\gamma > 0$, then 
	\begin{align*}
	\frac{\Vert v_\abg^r \Vert_1}{\Vert v_\abg^r \Vert_2} < \frac{\gamma}{\beta}.
	\end{align*}
\end{lemma}
The above Lemma states that those vector components $v_\abg^r$, which lie not too close to zero, are effectively sparse. Numerical experiments suggest that if $\hat{X}$ has $s$-sparse right components $\hat{v}^r$, ATLAS yields solutions with exactly sparse right components $v_\abg^r$. The theoretical necessity of considering 
{\it effective} sparsity also when $\hat{X}$ has $s$-sparse right components  is caused by the difficulty of obtaining better bounds on the support size of the vectors $v_\abg^r$.\\
To conclude, 
we can claim that $X_\abg$ is, even without any more specific requirements on $\A$, a reasonable approximation of $\hat{X}$, in the sense that it is of rank $R$, fulfills the measurements up to noise level, and has effectively sparse right components. However, the parameters $\alpha$ and $\beta$ have to be chosen with care, neither too small nor too large. Moreover, Lemma \ref{Bound-uv2} shows that $\alpha$ and $\beta$ have to be chosen of similar magnitude. Otherwise either left or right components of $X_\abg$ cannot be controlled.

\section{Recovery Properties of Minimizers of $J_\abg^R$ with RIP}\label{sec:RecoveryPropertiesRIP}

To explain the performance of ATLAS illustrated in Figure \ref{fig:Greyscale} we introduce two sets of matrices, which are sums of few rank-one matrices with sparse singular vectors.
\red{We stress here that we are not requiring the orthogonality of the components.}
 We also define corresponding additive RIPs, which are useful for proving the approximation result and can be seen as a generalization of the rank-$R$ and $(s_1,s_2)$-sparse RIP of Lee et.\ al.\ in  \cite{lee2013near}. 

\paragraph{Matrix models} The first matrix set is, \blue{for $\Gamma \ge 0$},
\begin{align} \label{S}
\begin{split}
S_{s_1,s_2}^{R,\Gamma} = \{ Z \in \R^{n_1\times n_2} \colon \exists\; u^1,...,u^R &\in \R^{n_1},\; v^1,...,v^R \in \R^{n_2}, \\
\text{ and } \sigma &= (\sigma_1,\dots,\sigma_R)^T \in \R^R, \text{ s.t. } \\
Z &= \sum_{r=1}^{R} \sigma_r u^r (v^r)^T, \\
\text{where } |\supp(u^r)| &\le s_1,\; |\supp(v^r)| \le s_2,\; \Vert u^r \Vert_2 = \Vert v^r \Vert_2 = 1,\\ 
\text{ for all } r \in [R],& \text{ and } \Vert \sigma \Vert_2 \le \Gamma \}. 
\end{split}
\end{align}
It contains all matrices $Z$ which can be decomposed into three matrices $U\Sigma V^T$ such that $U \in \R^{n_1\times R}$ and $V \in \R^{n_2 \times R}$ have $s_1$-sparse (resp.\ $s_2$-sparse) unit norm columns and $\Sigma \in \R^{R\times R}$ is the diagonal matrix defined by $\sigma$. The set is restricted to decompositions with $\Vert \Sigma \Vert_F \le \Gamma$.

The important difference w.r.t.\ \cite{lee2013near} is that the columns do not need to share a common support. Moreover, we do not require $U$ and $V$ to be orthogonal matrices. Nevertheless, all matrices $X$ with rank less or equal $R$, $s_1$-sparse (resp. $s_2$-sparse) left and right singular vectors, and $\Vert X \Vert_F \le \Gamma$ are in $S_{s_1,s_2}^{R,\Gamma}$. In this case $\Vert \Sigma \Vert_F = \Vert X \Vert_F$. We call such an admissible decomposition $U\Sigma V^T$ in \eqref{S} a Sparse Decomposition (SD) of $Z$. Note that the SD is not unique and that the SVD of $Z$ is not necessarily a SD of $Z$.\\

\red{We further generalize  $S_{s_1,s_2}^{R,\Gamma}$ to effectively sparse vectors. Recall the definition of $K_{n,s}$ in Definition \ref{def:EffectivelySparse}.} 
%
\red{\blue{For $\Gamma \ge 0$}, we define}
\begin{align}\label{K}
\begin{split}
K_{s_1,s_2}^{R,\Gamma} = \{ Z \in \R^{n_1\times n_2} \colon \exists\; u^1,...,u^R &\in K_{n_1,s_1},\; v^1,...,v^R \in 
K_{n_2,s_2}, \\
\text{ and } \sigma &= (\sigma_1,\dots,\sigma_R)^T \in \R^R, \text{ s.t. } \\
Z &= \sum_{r=1}^{R} \sigma_r u^r (v^r)^T, \\
\text{where } \Vert u^r \Vert_2 = \Vert v^r \Vert_2 &= 1, \text{ for all } r \in [R], \text{ and } \Vert \sigma \Vert_2 \le \Gamma \}
\end{split}
\end{align}
which is a relaxed version of $S_{s_1,s_2}^{R,\Gamma}$ as $S_{s_1,s_2}^{R,\Gamma} \subset K_{s_1,s_2}^{R,\Gamma}$. 
One of the most important features of the class $K_{s_1,s_2}^{R,\Gamma}$  is that it is to a certain extent closed under summation: in fact if $Z \in  K_{s_1,s_2}^{R,\Gamma}$ and $\hat Z \in  K_{\hat s_1,\hat s_2}^{R,\hat\Gamma}$ then 
\begin{equation}\label{sumrule}
Z - \hat Z \in  K_{\max\{s_1,\hat s_1\},\max\{s_2,\hat s_2\}}^{2R,\sqrt{\Gamma^2+\hat{\Gamma}^2}}.
\end{equation}
We call such an admissible decomposition $Z = U\Sigma V^T$ in \eqref{K} an effectively Sparse Decomposition of $Z$ and use the same shorthand notation, i.e., SD. The context makes clear which decomposition is meant. Any $\hat{X}$ decomposed as in \eqref{xSVD} belongs to $K_{n_1,s}^{R,\Gamma}$ if $\sum_{r=1}^{R} \Vert \hat{u}^r \Vert_2^2 \Vert \hat{v}^r \Vert_2^2 \le \Gamma^2$. Having the sets $S_{s_1,s_2}^{R,\Gamma}$ and $K_{s_1,s_2}^{R,\Gamma}$ at hand we now define corresponding RIPs.
\begin{definition}[Additive Rank-$R$ and (effectively) $(s_1,s_2)$-sparse RIP$_\Gamma$] \label{apprRIPDef}
	A linear operator $\A : \R^{n_1 \times n_2} \rightarrow \R^m$ satisfies the additive rank-$R$ and $(s_1,s_2)$-sparse RIP$_\Gamma$ with isometry constant $\delta > 0$ if
	\begin{align} \label{apprRIP}
	\left| \Vert \A(Z) \Vert_2^2 - \Vert Z \Vert_F^2 \right| \le \delta,
	\end{align}
	for all $Z \in S_{s_1,s_2}^{R,\Gamma}$. \\
	If \eqref{apprRIP} holds for all $Z \in K_{s_1,s_2}^{R,\Gamma}$, we say $\A$ has the additive rank-$R$ and effectively $(s_1,s_2)$-sparse RIP$_\Gamma$. Note that the rank-$R$ and effectively $(s_1,s_2)$-sparse RIP$_\Gamma$ implies the rank-$R$ and $(s_1,s_2)$-sparse RIP$_\Gamma$ as $S_{s_1,s_2}^{R,\Gamma} \subset K_{s_1,s_2}^{R,\Gamma}$.
\end{definition}
\blue{
\begin{remark}\label{Gammacond2}
By not enforcing orthogonality, the SDs allow certain ambiguities. In particular, any $\hat X \in K_{s_1,s_2}^{R,\Gamma}$ could as well be decomposed as follows
$$
\hat X=  \sum_{r=1}^{R} \sigma_r u^r (v^r)^T =\sum_{r=1}^{R} \sum_{j=1}^K  \frac{\sigma_r}{K} u^r (v^r)^T,  \quad \left ( \sum_{r=1}^{R} \sum_{j=1}^K  \frac{\sigma_r^2}{K^2} \right)^{1/2}  \leq \frac{\Gamma}{\sqrt K},
$$
for any $K \in \mathbb N$, implying $\hat X \in K_{s_1,s_2}^{R K,\Gamma/\sqrt{K}}$. Since by this argument $K_{s_1,s_2}^{R,\Gamma} \subset K_{s_1,s_2}^{R K,\Gamma/\sqrt{K}}$ and consequently an RIP on $K_{s_1,s_2}^{R K,\Gamma/\sqrt{K}}$ is harder to satisfy, one is in general interested in choosing an SD of minimal complexity ($K = 1$).
\end{remark}
}
%


\paragraph{Recovery results} We are ready now to state the main recovery result: If one assumes RIP, any appropriate global minimizer of $J_\abg^R$ provides an approximation to $\hat{X}$, with an error bound depending on the magnitude of $\alpha$ and $\beta$, the sparsity $s$, the RIP constant $\delta$, and the magnitude of $\hat{X}$ measured in an appropriate Schatten quasi-norm. The approximation is worsened in an additive way by noise level.
\begin{theorem}[Approximation of $\hat{X}$] \label{ApproxX}
	Fix the positive constants $\alpha, \beta > 0$, \blue{$\Gamma \ge 0$}, and the effective sparsity indicator level $1\leq s \leq n_2$. Let $\A$ have the additive rank-$2R$ effectively $(n_1,{\max\{s,(\gamma/\beta)^2\}})$-sparse RIP$_{(c+1)\Gamma}$ with RIP-constant $0 < \delta < 1$, for a fixed choice of $\gamma > 0$ and $c \ge 1$.\\
	If $\hat{X} \in K_{n_1,s}^{R,\Gamma}$ is of rank $R$ and $y = \A(\hat{X}) + \eta \in \R^m$, then
	{\begin{align} \label{Approximation}
	\Vert \hat{X} - X_\abg \Vert_F \le \sqrt{s^{\frac{1}{3}} R^{\frac{2}{3}} C_{2,1}c_{\hat U}} \sqrt[6]{\alpha \beta^2} \Vert \hat{X} \Vert_{\frac{2}{3}}^\frac{1}{3} + 2\Vert \eta \Vert_2 + \sqrt{\delta},
	\end{align}}
	for any global minimizer $(u_\abg^1,...,v_\abg^R)$ of $J_\abg^R$ that fulfills $\Vert v_\abg^r \Vert_2 \ge (\| \hat{X} \|_{F} + \| \eta \|_2 + \sqrt{\delta})^2/\gamma$ for all $r \in [R]$ and $\Vert \sigma_\abg \Vert_F \le c\Gamma$ in \eqref{Xab}. In this case, in particular, $X_\abg \in K_{n_1,(\gamma/\beta)^2}^{R,c\Gamma}$ with the SD in \eqref{Xab}. 
\end{theorem}

There are some aspects of this result we would like to discuss before we proceed:
\begin{enumerate}
    \item[(a)] If we could take the limits $\alpha \rightarrow 0$ and $\beta \rightarrow 0$, the error in \eqref{Approximation} would vanish up to noise-level and RIP-constant. However, {this limit cannot be performed as} there are important restrictions dictated by the need of fulfilling simultaneously the RIP and the assumptions on $X_\abg$. If $\beta$ is getting small the conditions for having RIP degenerate, i.e., reconstruction for a fixed number of measurements only works up to a minimal $\beta$. Letting $\alpha$ go to zero while keeping $\beta$ fixed leads to minimizers which violate the lower bound $\Vert v_\abg^r \Vert_2 \ge (\| \hat{X} \|_{F} + \| \eta \|_2 + \sqrt{\delta})^2/\gamma$ or the upper bound $\Vert \sigma_\abg \Vert_F \le c\Gamma$. To see this, note that by Lemma \ref{Bound-uv2} small $\alpha$ leads to strict bounds on $\| v_\abg^r \|_2$ and weak bounds on $\| u_\abg^r \|_2$.\\
    \item[(b)]  Let us mention that in case $\hat{X} \in K_{n_1,s}^{R,\Gamma}$ and the SD of $\hat{X}$ coincides with its SVD, then in view of the identity \eqref{qneq0} the factor $c_{\hat U} R^{2/3}$ in the error estimates \eqref{Approximation} and \eqref{err1eq} can be substituted by $1$, hence there would be no dependence on the rank $R$.\\
    \item[(c)] In order to clarify how $(\gamma/\beta)^2$ and $s$ are related in the RIP in Theorem \ref{ApproxX} (and Corollary \ref{ApproxXcor} below), let us assume for simplicity that the SD of $\hat{X}$ coincides with its SVD and $\alpha = \beta$. Consequently, to get an error bound independent of $s$ in \eqref{Approximation}, $\alpha$ and $\beta$ have to be chosen of order $\mathcal{O}(s^{-\frac{1}{3}})$, i.e., $(\gamma/\beta)^2$ is of order $\mathcal{O}(s^{\frac{2}{3}})$  which means that an $(n_1,\gamma^2s)$-sparse RIP$_{(c+1)\Gamma}$ is sufficient for recovery.\\
    \item[(d)] The result only applies to minimizers whose scaling matrix $\Sigma_\abg$ is bounded in Frobenius norm and whose right components $v_\abg^r$ are not too close to zero. The first requirement is necessary as the RIP is restricted to SDs with scaling matrices within a ball around zero. The second one is needed to show some level of effective sparsity of the minimizers $X_\abg$ (see also the discussion in Section \ref{sec:Main1}). 
    While effective sparsity of (right) component vectors of $X_\abg$ is naturally wished and expected if $\hat X \in K_{n_1,s}^{R,\Gamma}$, we were not able in all cases to show {\it exact} sparsity of (right) component vectors of $X_\abg$ if $\hat X \in S_{n_1,s}^{R,\Gamma}$, but again only their effective sparsity. Hence, we are bound to using as an artifact of the proof the stronger effectively $(s_1,s_2)$-sparse RIP$_\Gamma$ for theoretical analysis also in this case. In numerical experiments, however, for $\hat X \in S_{n_1,s}^{R,\Gamma}$ the obtained minimizers $X_\abg$ are empirically {\it exactly} sparse (not just effectively sparse) and, hence, the weaker rank-$2R$ $(s_1,s_2)$-sparse RIP$_\Gamma$ might suffice in practice. The latter can already be guaranteed for a smaller number of measurements.\\
    \item[(e)] As argued in Section \ref{MainProofSection} the above theorem can be straightforwardly extended to sparsity on left component vectors. In this case $J_\abg^R$ has to be adapted by considering $\ell_1$-norm penalties on the $u$-components.\\
    \item[(f)] It is important to require $\rank(\hat{X}) = R$ as otherwise the equivalence of Schatten-norm and normed SD cannot be guaranteed as \eqref{qneq}. If the SD of $\hat{X}$ coincides with its SVD though, the rank condition may be dropped.
\end{enumerate}

\blue{By choosing $\alpha$ and $\beta$ in relation to the noise-to-signal ratio $\Vert \eta \Vert_2^2/\Vert \hat{X} \Vert_{\frac{2}{3}}^{\frac{2}{3}}$ we obtain the following version of Theorem \ref{ApproxX}, which has the form of a typical compressed sensing recovery bound. Assuming the RIP, the approximation error is linear in noise level while the slope of the linear function depends on sparsity level and possibly the rank. However, peculiarly, for a fixed number of measurements the RIP fails for exceedingly small noise and correspondingly small $\alpha$ and $\beta$, cf.\ (a) in the discussion of Theorem \ref{ApproxX}.} \BLUE{To be more precise, if $\alpha \approx \beta \approx \Vert \eta \Vert_2^2/\Vert \hat{X} \Vert_{\frac{2}{3}}^{\frac{2}{3}}$ and $\eta$ is small, it might happen that there exists no viable choice of $\gamma > 0$ simultaneously fulfilling the requirements on $\A$ and $\| v_\abg^r \|_2$ in Theorem \ref{ApproxX}. Hence, the result is valid only for sufficiently small signal-to-noise ratios. On the one hand, we will show in Section \ref{Numerics} with numerical experiments, this apparently counterintuitive result is factual and not an \purple{artifact} of the proof technique. A possible intuitive explanation is that $J_\abg^R$ becomes {a mere least-squares} without sparsifying effect for $\alpha$ and $\beta$ close to zero, which is caused by vanishing noise. On the other hand, Section \ref{Numerics} and, in particular, the discussion in Section \ref{Numerics1} demonstrate that by slightly overestimating the noise-level, ATLAS is still practical in low-noise settings. }
\begin{corollary} \label{ApproxXcor}
	Let {$\hat{X} \in K_{n_1,s}^{R,\Gamma}$ with $\rank(\hat{X}) = R$ fulfill the noisy measurements $y = \A(\hat{X}) + \eta$} and let $\alpha = \beta = \Vert \eta \Vert_2^2/\Vert \hat{X} \Vert_{\frac{2}{3}}^{\frac{2}{3}} < 1$. Assume $\A$ has for some $\gamma>0$ and $c \ge 1$ the additive rank-$2R$ effectively $\left(n_1,{\max}\{s,\gamma^2 (\Vert \hat{X} \Vert_{\frac{2}{3}}^{\frac{2}{3}}/\Vert \eta \Vert_2^2)^2 \}\right)$-sparse RIP$_{(c+1)\Gamma}$ with  RIP-constant $0 < \delta < 1$. Then, for $X_{\alpha, \beta}$ with $\Vert \Sigma_\abg \Vert_F \le c\Gamma$ and $\Vert v_{\alpha, \beta}^r \Vert_2 \ge (\| \hat{X} \|_{F} + \| \eta \|_2 + \sqrt{\delta})^2/\gamma$, $r \in [R]$, we have
	\begin{align}\label{err1eq}
	\Vert \hat{X} - X_\abg \Vert_F \le \left( 2 \sqrt{c_{\hat U} R^{2/3} s^{1/3} } + 2 \right) \Vert \eta \Vert_2 + \sqrt{\delta}.
	\end{align}
\end{corollary} 
\blue{
\begin{remark} \label{ApproxXcorRemark}
One could object that the simple zero solution $\bar X=0$ is already a competitor in case of large noise $\|\eta \|_2 \geq c \| \hat X\|_F$, for $c > 0$, i.e.,
\begin{equation}\label{naive}
\Vert \hat{X} -  \bar X \Vert_F \le c^{-1} \|\eta \|_2.
\end{equation}
However, for a larger number $m$ of measurements we can consider lower level of noise, i.e., $ c \to 0$ and the bound \eqref{naive} would explode, while \eqref{err1eq} would remain effective.  Moreover, our numerical experiments shows empirically that also in case of larger noise level, computing $X_\abg$ gives a solution, which outperforms not only trivial competitors  as $\bar X$, but also state-of-the-art methods such as SPF.
\end{remark}
}

\section{RIP Results for Subgaussian Operators} \label{SectionRIP}

\blue{As already mentioned in the end of Section \ref{ProblemSection}, a linear operator $\A$ of the form \eqref{A} which is drawn from a subgaussian distribution fulfills the above introduced RIPs with high probability.} This is stated in the following Lemma. We first recall the definition of subgaussian random variables (for further details see \cite{vershynin2010introduction}).

\begin{definition}[Subgaussian Random Variable]
	A random variable $\xi \in \R$ is called $\mathcal{K}$-subgaussian if the tail bound $\Pr[]{|\xi| > t} \le C\exp(-ct^2/\mathcal{K}^2)$ holds where $c,C > 0$ are absolute constants. The smallest possible number for $\mathcal{K} > 0$ is called subgaussian norm of $\xi$ and denoted by $\Vert \xi \Vert_{\psi_2}$. 
\end{definition}

\begin{remark}
	The class of subgaussian random variables covers important special cases as Gaussian, Bernoulli, and more generally all bounded random variables (see \cite{vershynin2010introduction}).
\end{remark}
\begin{lemma}[RIP for Subgaussian Operators] \label{GaussianRIP}
	\blue{Let $\Gamma \ge 0$} and let $\A \colon \R^{n_1\times n_2} \rightarrow \R^m$ be the linear measurement operator of form \eqref{A}. Assume, all $A_i$, for $1 \le i \le m$, have i.i.d.\ $\mathcal{K}$-subgaussian entries $a_{i,j,k}$ with mean $0$ and variance $1$. If
	\begin{align} \label{meas1}
	     m \gtrsim \left( \frac{\delta}{\Gamma^2 R} \right)^{-2} R (s_1+s_2) \log \left( \max \{n_1,n_2\} \right)
	\end{align}
	then ${\A}$ has the additive rank-$R$ and $(s_1,s_2)$-sparse RIP$_\Gamma$ with isometry constant $\delta \in (0,\Gamma^2 R)$ with probability at least $1-2\exp(-C (\delta/\Gamma^2 R) m)$ where $C > 0$ is a constant depending on $\mathcal{K}$. If
	\begin{align} \label{meas2}
	     m \gtrsim \left( \frac{\delta}{\Gamma^2 R} \right)^{-2} R (s_1+s_2) \log^3 \left( \max \{n_1,n_2\} \right)
	\end{align}
	then ${\A}$ has the additive rank-$R$ and effectively $(s_1,s_2)$-sparse RIP$_\Gamma$ with isometry constant $\delta \in (0,\Gamma^2 R)$ with probability at least $1-2\exp(-C'(\delta/\Gamma^2 R) m)$ where $C' > 0$ is a constant depending on $\mathcal{K}$.
\end{lemma} 
\begin{remark}\label{rem:dim}
	Lemma \ref{GaussianRIP} states, for $\delta = \Delta (\Gamma^2 R)$, $\Delta \in (0,1)$, that, up to log-factors, $m \approx \mathcal{O} \left( \Delta^{-2} R(s_1+s_2)) \right)$ subgaussian measurements are sufficient to have $\delta$-stable embeddings of $S_{s_1,s_2}^{R,\Gamma}$ and $K_{s_1,s_2}^{R,\Gamma}$ (cf.\ \cite[Def. 1.1 \& Thm. 1.5]{Plan2014}). Note that $\Gamma^2 R$ is the squared Frobenius diameter of $S_{s_1,s_2}^{R,\Gamma}$ and $K_{s_1,s_2}^{R,\Gamma}$. 
	
	As we restrict ourselves below to {$s$-effective \purple{sparse right component vectors}} of $\hat{X}$, we only use the rank-$R$ and (effectively) $(n_1,s)$-sparse RIP$_\Gamma$. 
    \purple{For the presented results to have some meaning, a typical dimensional setting is} $R \ll s \approx n_1 \ll n_2$. In fact, if $n_1$ were close to $n_2$ in magnitude, the sparsity $s$ of the right component vectors would not be useful to reduce the order of the measurements $m\approx \mathcal O( R(n_1+s)) \approx \mathcal O( Rn_1) \approx  \mathcal O( Rn_2)$. Moreover, if $R$ were close to $n_1$, the matrix would not be low-rank as $n_1$ would be the maximal possible rank.\\ 
  In \cite{lee2013near} the authors give information theoretical lower bounds on the necessary number of measurements for reconstructing low-rank matrices with sparse singular vectors sharing a common support, namely $m \gtrsim R(s_1+s_2)$. As we do not require orthogonality of SDs in $S_{s_1,s_2}^{R,\Gamma}$ resp. $K_{s_1,s_2}^{R,\Gamma}$ (excluding a scaling invariant RIP which is independent of the set diameter, see Remark \ref{DominiksRemark}), the bounds in \eqref{meas1} and \eqref{meas2} are up to $\log$-factors at the information theoretic limit for the class of matrices in \cite{lee2013near}. We are not aware of any information theoretical lower bounds for the more general class of matrices considered in the present paper.
\end{remark}
\begin{remark} \label{DominiksRemark}
	The additive RIP in \eqref{apprRIP} differs from the commonly used multiplicative RIPs of the form
	\begin{align} \label{multRIP}
	(1-\delta) \Vert Z \Vert_F^2 \le \Vert \A(Z) \Vert_2^2 \le (1+\delta) \Vert Z \Vert_F^2
	\end{align}
	as it is not scaling invariant and $\A(Z) = \A(Z')$ does not imply $Z = Z'$ but only $\Vert Z - Z' \Vert_2^2 \le \delta$. In fact it is not possible to derive a classical scaling invariant RIP like \eqref{multRIP} on $K_{s_1,s_2}^{R,\Gamma}$ under similar conditions as \eqref{meas2}. The main problem is non-orthogonality of the SD. A simple example illustrates this point: Assume $R=2$ and $m \simeq 2 (n_1+s) \log^3 \left( \max \{n_1,n_2\} \right)$ and the linear operator $\A$ fulfills \eqref{multRIP} for all $Z \in K_{n_1,s}^{2,1}$. Choose some $u \in \R^{n_1}, v_1 \in \R^{n_2}$ of unit norm and $\Vert v_1 \Vert_1 \le \sqrt{s}/2$. Define $v_2 := -v_1 + \eps w$ for any $w \in \R^{n_2}$ and choose $\eps > 0$ sufficiently small to ensure $\Vert v_2 \Vert_1 \le \sqrt{s}$ and $\Vert v_2 \Vert_2 \approx 1$. Then $Z := (1/2)uv_1^T + (1/2)uv_2^T \in K_{n_1,s}^{2,1}$ and \eqref{multRIP} holds. But this implies by definition of $Z$ and scaling invariance of \eqref{multRIP} that
	\begin{align*}
		(1-\delta) \Vert uw^T \Vert_F^2 \le \Vert \A(uw^T) \Vert_2^2 \le (1+\delta) \Vert uw^T \Vert_F^2
	\end{align*}
	which means the RIP directly extends to all rank-$1$ matrices (not only those with sparse right component). If $n_1,s \ll n_2$, this is a clear contradiction to information theoretical lower bounds, as corresponding RIPs would require at least $m \simeq \max\{n_1,n_2\}$ (see \cite[Section 2.1]{candes2011tight}).
\end{remark}

\section{Convergence of ATLAS} \label{LocalConvergenceSection}
\red{In the following by adapting results of Attouch et.\ al.\ in \cite{attouch2010proximal} we show convergence of ATLAS.}
\red{Specifically,} there is a neighborhood $\mathcal{U}_{(u_\abg^1,...,v_\abg^R)}$ of a global minimizer $(u_\abg^1,...,v_\abg^R)$ such that the sequence $(u_k^1,...,v_k^R)$ defined by \eqref{ATLAS0} converges to $(u_\abg^1,...,v_\abg^R)$  if the initialization lies within $\mathcal{U}_{(u_\abg^1,...,v_\abg^R)}$. However, we do not provide proof for any initialization to fulfill the requirement \red{and we leave this open issue} for future research, \red{cf.\ Remark \ref{rem:Radius} below.}  The techniques in \cite{attouch2010proximal} also \purple{might be adjusted} for an analysis of rate of convergence of ATLAS, but this would go beyond the scope of this work and {is} \purple{a topic for future investigation.} 
We begin by a generalization of the basic conditions of \cite{attouch2010proximal}. Let $L$ be a functional of the following form:
\begin{align*}
(H) \tab 
&\begin{cases}
L(u^1,\dots,u^R,v^1,\dots,v^R) = \sum_{r=1}^R f_r(u^r) + Q(u^1,\dots ,v^R) + \sum_{r=1}^R
g_r(v^r), \\
f_r: \mathbb{R}^{n_1} \rightarrow {\mathbb{R} \cup \{ \infty \}}, \; g_r:\mathbb{R}^{n_2} \rightarrow {\mathbb{R} \cup \{ \infty \}} 
\text{ are proper lower semicontinuous, for } 1 \leq r 
\leq R, \\
Q: \mathbb{R}^{n_1} \times \cdots \times \mathbb{R}^{n_1} \times \mathbb{R}^{n_2} \times
\cdots \times \mathbb{R}^{n_2} \rightarrow \mathbb{R} \text{ is a }
C^1 \text{ function}, \\
\nabla Q \text{ is Lipschitz continuous on bounded subsets of } \mathbb{R}^{n_1} \times
\cdots \times \mathbb{R}^{n_1} \times \mathbb{R}^{n_2} \times \cdots \times 
\mathbb{R}^{n_2}.
\end{cases}
\intertext{For given $(u_0^1,\dots ,v_0^R) \in (\mathbb{R}^{n_1})^R \times (\mathbb{R}^{n_2})^R$ and fixed sequences $(\lambda_k^1)_{k \in \mathbb{N}},\dots ,(\lambda_k^R)_{k \in \mathbb{N}}$ and
$(\mu_k^R)_{k \in \mathbb{N}},\dots ,(\mu_k^R)_{k \in \mathbb{N}}$ assume that}
(H1) \tab 
&\begin{cases}
\inf L > -\infty, \\
L(\cdot ,u_0^2,\dots ,v_0^R) \text{ is proper}, \\
\text{for some positive } r_- < r_+ \text{ the sequences } \lambda_k^1,\dots,\mu_k^R
\text{ belong to } (r_- , r_+).
\end{cases}
\end{align*}
The adapted main result of \cite{attouch2010proximal} now guarantees convergence of the so-called Proximal Alternating Minimization
\begin{align} \label{ProxAlgo0}
	(\text{PAM}) \tab
	&\begin{cases}
	u_{k+1}^1 = \argmin_{u \in \R^{n_1}} L(u,u_k^2,\dots,u_k^R,v_k^1,\dots,v_k^R) + \frac{1}{2\lambda_k^1} \Vert u - u_k^1 \Vert_2^2,\\
	v_{k+1}^1 = \argmin_{v \in \R^{n_2}} L(u_{k+1}^1,u_k^1,\dots,u_k^R,v,v_k^2\dots,v_k^R) + \frac{1}{2\mu_k} \Vert v - v_k^1 \Vert_2^2,\\
	\vdots \\
	u_{k+1}^R = \argmin_{u \in \R^{n_1}} L(u_{k+1}^1,\dots,u_{k+1}^{R-1},u,v_{k+1}^1,\dots,v_{k+1}^{R-1},v_k^R) + \frac{1}{2\lambda_k} \Vert u - u_k^R \Vert_2^2,\\
	v_{k+1}^R = \argmin_{v \in \R^{n_2}} L(u_{k+1}^1,\dots,u_{k+1}^R,v_{k+1}^1,\dots,v_{k+1}^{R-1},v) + \frac{1}{2\mu_k} \Vert v - v_k^R \Vert_2^2,
	\end{cases}
\end{align} 
to a stationary point of $L$ (resp. convergence to a global minimizer $(u_\ast^1,\dots,v_\ast^R)$ of $L$ if the initialization $(u_0^1,\dots,v_0^R)$ of (PAM) lies sufficiently close to $(u_\ast^1,\dots,v_\ast^R)$) if $L$ fulfills $(H)$, $(H1)$ and the so called Kurdyka-Lojasiewicz Property, which requires $L$ to behave well around stationary points.
\begin{definition}[Kurdyka-Lojasiewicz Property]
	A proper lower semicontinuous function $f:\mathbb{R}^n \rightarrow {\mathbb{R} \cup \{ \infty \}}$
	is said to have the KL-property at $\ol{x} \in \dom \, \partial f$\footnote{Here $\partial f$ denotes the subdifferential of $f$ and $\dom \, \partial f$ the domain on which $\partial f$ takes finite values.} if there exist $\eta \in 
	\left( 0,\infty \right]$, a neighborhood $U$ of $\overline{x}$ and a continuous concave
	function $\varphi : \left[ 0,\infty \right) \rightarrow \mathbb{R}_+$ such that
	\begin{enumerate}
		\item[-] $\varphi(0) = 0$,
		\item[-] $\varphi$ is $C^1$ on $(0,\eta)$, 
		\item[-] $\varphi'(t) > 0$, for all $t \in (0,\eta)$,
		\item[-] and, for all $x \in U \cap \{ x \in \mathbb{R}^n : f(\overline{x}) < f(x) <
		f(\overline{x}) + \eta \}$, the KL-inequality holds:
		\begin{align*}
		\varphi'(f(x)-f(\overline{x}))\; \dist(0,\partial f(x)) \geq 1.
		\end{align*}
	\end{enumerate}
\end{definition}
\begin{theorem}[Local Convergence to Global {Minimizers}] \label{LocalConvergence}
	Assume that $L$ satisfies $(H)$, $(H1)$. If $L$ has the Kurdyka-Lojasiewicz property at its global minimizer $(u_\ast^1,\dots,v_\ast^R)$, {then} there exist $\eps, \eta > 0$, such that {the initial conditions
	\begin{align*}
	\Vert (u_0^1,\dots,v_0^R) - (u_\ast^1,\dots,v_\ast^R) \Vert_2 < \eps, \tab \min L < L(u_0,v_0) < \min L + \eta,
	\end{align*}
	imply that} the {iterations $(u_k^1,\dots,v_k^R)$ generated by} (PAM) converge to $(u_*^1,\dots,v_*^R)$\\ $\in \argmin L$. If $L$ has the Kurdyka-Lojasiewicz at each point of its domain, {then} either $\Vert (u_k^1,\dots,v_k^R) \Vert_2 \rightarrow \infty$ or $(u_k^1,\dots,v_k^R)$ converges to a stationary point of $L$.
\end{theorem}

\begin{remark}{} \label{rem:Radius}
	Let us briefly note two observations with regard to Theorem \ref{LocalConvergence}:
    \begin{enumerate}
    	\item[(a)] The main difficulty in characterizing the convergence radius is to characterize the KL-parameters $U$ and $\eta$ of $L$. Doing so for a non-convex functional like $J_\abg^R$ is a challenging task on its own and thus the main reason for us to defer the treatment of initialization to future work.\\ 
    	\item[(b)] We will see below that $J_\abg^R$ has the KL-property with $\varphi(t) = ct^{1-\theta}$, for $c > 0$ and $\theta \in [0,1)$. As \cite{attouch2010proximal} shows, a characterization of $\theta$ would determine the convergence speed of the alternating minimization of $L$. While \cite{li2013global} can be used to compute $\theta$ for piecewise convex polynomials, it is unclear how to do the same for non-convex polynomials. Addressing this more general issue would in particular provide a convergence speed analysis of ATLAS.
    \end{enumerate} 
\end{remark}

Theorem \ref{LocalConvergence} is a straight-forward adaption of the results in \cite{attouch2010proximal}. We defer the details to the Appendix. By applying Theorem \ref{LocalConvergence} to {$L=J_\abg^R$} and ATLAS we obtain convergence to stationary points and local convergence to global minimizers as the sequence $(u_k^1,\dots,v_k^R)$ is bounded by coercivity of $J_\abg^R$. One can check that conditions $(H)$, $(H1)$ are fulfilled by $J_\abg^R$ and ATLAS for a suitable choice of the sequences $(\lambda_k^1)_{k \in \mathbb{N}},\dots ,(\lambda_k^R)_{k \in \mathbb{N}}$, $(\mu_k^R)_{k \in \mathbb{N}},\dots ,(\mu_k^R)_{k \in \mathbb{N}}$. It remains to validate the KL-property. As mentioned in \cite[Section 4.3]{attouch2010proximal}, all semialgebraic functions satisfy the KL-property at each point with $\varphi(t) = ct^{1-\theta}$ for some $\theta \in [0,1) \cap \mathbb{Q}$ and $c > 0$. Hence, by showing that $J_\abg^R$ is semialgebraic, we get the KL-property for free. But we pay the price of having no better knowledge on the parameters $\eps$ and $\eta$ in Theorem \ref{LocalConvergence}, which characterize the convergence radius. {Therefore, let us conclude} by showing that $J_\abg^R$ is semialgebraic, i.e., $\graph(J_\abg^R) \subset \mathbb{R}^{Rn_1+Rn_2} \times \mathbb{R}$ is a semialgebraic set.\\
A set in $\mathbb{R}^d$ is called semialgebraic if it can be written as a finite union of sets of the form
\begin{align*}
\{ x \in \mathbb{R}^d \; : \; p_i(x) = 0, \; q_i(x) \purple{>} 0, \; i = 1,\dots,p \},
\end{align*}
where $p_i,q_i$ are real polynomials. 
First, the absolute value of one component of a vector $h(x) := \vert x_l \vert$ is a semialgebraic function as
\begin{align*}
\graph (h) = \{ (x,r) \in \mathbb{R}^d \times \mathbb{R} : x_i + r = 0, \; x_i < 0 \}
\cup \{ (x,r) \in \mathbb{R}^d \times \mathbb{R} : x_i = 0, \; r = 0 \} \\
\cup \{ (x,r) \in \mathbb{R}^d \times \mathbb{R} : x_i - r = 0, \; -x_i < 0 \}.
\end{align*}
Second, it is clear that polynomials $p$ are semialgebraic as $\graph (p) = \{ (x,r) \in \mathbb{R}^d \times \mathbb{R} : p(x) - r = 0 \}$ and, third, composition, finite sums and finite products of semialgebraic functions are semialgebraic. The semialgebraicity of $J_\abg^R$ follows as
\begin{align*}
J_\abg^R(u^1,\dots ,v^R) = \sum_{l = 1}^m \vert y_l - \sum_{r = 1}^R \langle A_l , u^r {v^r}^T 
\rangle_F \vert^2 + \alpha \sum_{r = 1}^R \sum_{l = 1}^{n_1} \vert u_l^r \vert^2 + \beta 
\sum_{r=1}^R \sum_{l = 1}^{n_2} \vert v_l^r \vert
\end{align*}
is just a finite composition of semialgebraic basic units.


\section{Proofs} \label{Proofs}

This section provides proofs for the main results from Sections \ref{sec:Main1}-\ref{LocalConvergenceSection}. Some merely technical parts are moved to the Appendix to ease the reading. We begin by showing the general properties of global minimizers $(u_\abg^1,...,v_\abg^R)$ of $J_\abg^R$ and proving Theorem \ref{ApproxX}. Then, we present the proof of \Cref{GaussianRIP}. The proof of \Cref{LocalConvergence} can be found in the Appendix, as it is based on straightforward modifications of the {arguments in} \cite{attouch2010proximal}.

\subsection{Bounds on Minimizers} \label{BoundSection}

Recall the SD related representation $\hat{X} = \sum_{r=1}^R \hat{u}^r (\hat{v}^r)^T$ in \eqref{xSVD} where $\hat{\sigma}_r = \Vert \hat{u^r} \Vert_2 \Vert \hat{v}^r \Vert_2$ and the notation $X_\abg = \sum_{r=1}^{R} u_\abg^r (v_\abg^r)^T$. For proving Proposition \ref{Bound-y2} and Lemma \ref{Bound-uv2} we need following technical lemma.

\begin{lemma} \label{fpq}
	Let $\alpha,\beta,a,b,p,q > 0$. Then 
	\begin{align*}
	f: \R^+ \rightarrow \R, \tab f(\lambda) := \lambda^p \alpha a + \frac{1}{\lambda^q} \beta b,
	\end{align*}
	attains its minimum at $\tilde{\lambda} = \left( \frac{q}{p} \frac{\beta b}{\alpha a} \right)^\frac{1}{p+q}$ and has the minimal value
	\begin{align*}
	\min f = f(\tilde{\lambda}) = C_{p,q} (\alpha a)^\frac{q}{p+q} (\beta b)^\frac{p}{p+q},
	\end{align*}
	where $C_{p,q} = \left( \frac{q}{p} \right)^\frac{p}{p+q} + \left( \frac{p}{q} \right)^\frac{q}{p+q}$.
\end{lemma}

\begin{Proof}[of Lemma \ref{fpq}]
	The result is obtained by differentiation of $f$ and by searching for its derivative's zeros.
\end{Proof}

\begin{Proof}[of Proposition \ref{Bound-y2}]
	By applying Lemma \ref{fpq} $R$ times using $p=2, q=1, a=\Vert \hat{u}^r \Vert_2^2, b = \Vert \hat{v}^r \Vert_1$ we get $\tilde{\lambda}_1,...,\tilde{\lambda}_R$, such that
	\begin{align}
	\begin{split} \label{boundEstimate2}
	J_\abg^R (\tilde{\lambda}_1 \hat{u}^1,...,\tilde{\lambda}_R \hat{u}^R, \frac{1}{\tilde{\lambda}_1} \hat{v}^1,..., \frac{1}{\tilde{\lambda}_R} \hat{v}^R) &= \Vert y - \A(\hat{X}) \Vert_2^2 + \sum_{r=1}^{R} C_{2,1} \sqrt[3]{\alpha \beta^2} \sqrt[3]{\Vert \hat{u}^r \Vert_2^2 \Vert \hat{v}^r \Vert_1^2} \\ 
	&= \Vert \eta \Vert_2^2 + C_{2,1} \sqrt[3]{\alpha \beta^2} \sum_{r=1}^{R} \left( \Vert \hat{u}^r \Vert_2 \Vert \hat{v}^r \Vert_1 \right)^\frac{2}{3}.
	\end{split} 
	\end{align}
	Note that, although not explicitly labeled, each $\tilde{\lambda}_r$ depends on the choice of $\alpha$ and $\beta$ as well as on $a,b,p$, and $q$. The minimality of $(u_\abg^1,...,v_\abg^R)$ implies
	\begin{align*}
	\Vert y - \A(X_\abg) \Vert_2^2 &\le J_\abg^R(u_\abg^1,...,v_\abg^R) \le J_\abg^R(\tilde{\lambda}_1 \hat{u}^1,...,\tilde{\lambda}_R \hat{u}^R, \frac{1}{\tilde{\lambda}_1} \hat{v}^1,..., \frac{1}{\tilde{\lambda}_R} \hat{v}^R) \\
	&= \Vert \eta \Vert_2^2 + C_{2,1} \sqrt[3]{\alpha \beta^2} \sum_{r=1}^{R} \left( \Vert \hat{u}^r \Vert_2 \Vert \hat{v}^r \Vert_1 \right)^\frac{2}{3}
	\end{align*}
	which is the claim.
\end{Proof}

The proof of Lemma \ref{Bound-uv2} works in a similar way.
\begin{Proof}[of Lemma \ref{Bound-uv2}]
	From \eqref{boundEstimate2} in the proof of Proposition \ref{Bound-y2} we obtain
	\begin{align*} 
	\| y - \A(X_\abg) \|_2^2 + \sum_{r=1}^{R} \left( \alpha \Vert u_\abg^r \Vert_2^2 + \beta \Vert v_\abg^r \Vert_1 \right) &= J_\abg^R(u_\abg^1,...,v_\abg^R) \\ 
	&\le J_\abg^R(\tilde{\lambda}_1 \hat{u}^1,...,\tilde{\lambda}_R \hat{u}^R, \frac{1}{\tilde{\lambda}_1} \hat{v}^1,..., \frac{1}{\tilde{\lambda}_R} \hat{v}^R) \\
	&= \Vert \eta \Vert_2^2 + C_{2,1} \sqrt[3]{\alpha \beta^2} \sum_{r=1}^{R} \left( \Vert \hat{u}^r \Vert_2 \Vert \hat{v}^r \Vert_1 \right)^\frac{2}{3}
	\end{align*}
	The first part of the claim follows by subtracting $\| y - \A(X_\abg) \|_2^2$ on both sides, leaving out half of the terms on the left-hand side, and dividing by $\alpha$ (resp.~$\beta$). To show the second part, note that by minimality of $(u_\abg^1,...,v_\abg^R)$ and Lemma \ref{fpq}
	\begin{align*}
		\sum_{r=1}^{R} \left( \alpha \Vert u_\abg^r \Vert_2^2 + \beta \Vert v_\abg^r \Vert_1 \right) = C_{2,1} \sqrt[3]{\alpha \beta^2} \sum_{r=1}^{R} \left( \Vert u_\abg^r \Vert_2 \Vert v_\abg^r \Vert_1 \right)^\frac{2}{3}
	\end{align*}
	and hence
	\begin{align*}
		\| y - \A(X_\abg) \|_2^2 + C_{2,1} \sqrt[3]{\alpha \beta^2} \sum_{r=1}^{R} \left( \Vert u_\abg^r \Vert_2 \Vert v_\abg^r \Vert_1 \right)^\frac{2}{3} &= J_\abg^R(u_\abg^1,...,v_\abg^R) \\
		&\le J_\abg^R(\tilde{\lambda}_1 \hat{u}^1,...,\tilde{\lambda}_R \hat{u}^R, \frac{1}{\tilde{\lambda}_1} \hat{v}^1,..., \frac{1}{\tilde{\lambda}_R} \hat{v}^R) \\
		&= \Vert \eta \Vert_2^2 + C_{2,1} \sqrt[3]{\alpha \beta^2} \sum_{r=1}^{R} \left( \Vert \hat{u}^r \Vert_2 \Vert \hat{v}^r \Vert_1 \right)^\frac{2}{3}.
	\end{align*}
	Subtracting $\| y - \A(X_\abg) \|_2^2$ on both sides and dividing by $C_{2,1} \sqrt[3]{\alpha \beta^2}$ concludes the proof.
\end{Proof}

To show the effective sparsity as in Lemma \ref{SparsityControl}, we combine the fact that $X_\abg$ is a minimizer with the assumed lower bound on $v_\abg^r$.
\begin{Proof}[of Lemma \ref{SparsityControl}]
	\newcommand{\rt}{{r}}
	By comparing $J_\abg^R (u_\abg^1,...,v_\abg^r)$ to $J_\abg^R (0,...,0)$, we get
	\begin{align*}
		\sum_{r=1}^{R} \left( \alpha \Vert u_\abg^r \Vert_2^2 + \beta \Vert v_\abg^r \Vert_1 \right) &\le J_\abg^R(u_\abg^1,...,v_\abg^R) \le J_\abg^R(0,...,0) = \Vert y \Vert_2^2.
	\end{align*}
	This implies $\Vert v_\abg^r \Vert_1 < \Vert y \Vert_2^2/\beta$. As by assumption $\Vert v_\abg^r \Vert_2 \ge \Vert y \Vert_2^2/\gamma$, we conclude
	\begin{align*}
	\frac{\Vert v_\abg^\rt \Vert_1}{\Vert v_\abg^\rt \Vert_2}  < \frac{\Vert y \Vert_2^2}{\beta} \frac{\gamma}{\Vert y \Vert_2^2} = \frac{\gamma}{\beta}.
	\end{align*}
\end{Proof}

\subsection{Proof of Theorem \ref{ApproxX}} \label{MainProofSection}

We have now all necessary tools at hand to prove our main approximation result. Most of the \purple{technical} work has been \purple{already presented} in {Proposition \ref{Bound-y2} and Lemma \ref{SparsityControl}}. By combining the RIP with the above bounds on norms and sparsity of minimizers, we can estimate the worst-case distance between $\hat{X}$ and $X_\abg$ depending on the size of $\alpha$ and $\beta$, the sparsity $s$, the RIP constant $\delta$, and the size of $\hat{X}$ measured in a Schatten quasi-norm.

As the reader may notice, all technical results of \Cref{BoundSection} can be adapted to effective sparsity of the left components $(u_\abg^1,...,u_\abg^R)$ as well. \purple{This can be done by  replacing $\ell_2$-norms by corresponding $\ell_1$-norms in $J_\abg^R$}. The proof Lemma \ref{SparsityControl}, which guarantees effective sparsity of the right components, is independent of the minimization of the left components. Therefore, \Cref{SparsityControl} applies also to the left components if $\ell_2$-norms are replaced by $\ell_1$-norms in $J_\abg^R$. \Cref{ApproxX} then can be adapted to this setting in a straightforward way.

\begin{Proof}[of Theorem \ref{ApproxX}]
	As $\Vert y \Vert_2 \le \Vert \A(\hat{X}) \Vert_2 + \| \eta \|_2 \le (\| X \|_F + \sqrt{\delta}) + \| \eta \|_2$, \Cref{SparsityControl} applies and yields that $X_\abg$ is in $K_{n_1,(\gamma/\beta)^2}^{R,c\Gamma}$. Combined with {$\hat{X} \in K_{n_1,s}^{R,\Gamma}$}, we know from {\eqref{sumrule}} that the difference $\hat{X} - X_\abg \in K_{n_1,{\max} \{s,(\gamma/\beta)^2\} }^{2R,(c+1)\Gamma}$. Hence, we  apply the rank-$2R$ and effectively $(n_1,{\max} \{ s,(\gamma/\beta)^2 \})$-sparse RIP$_{(c+1)\Gamma}$ of $\A$ to obtain (note that $|a^2 - b^2| \le \delta$ implies $|a - b| \le \sqrt{\delta}$, for $a,b > 0$)
	\begin{align*}
	\Vert \hat{X} - X_\abg \Vert_F &\le \Vert \A(\hat{X}) - \A(X_\abg) \Vert_2 + \sqrt{\delta} \le \left( \Vert y - \A(X_\abg) \Vert_2 + \Vert \eta \Vert_2 \right) + \sqrt{\delta} \\
	&\le \sqrt{s^{\frac{1}{3}} R^{\frac{2}{3}}C_{2,1}c_{\hat U} \sqrt[3]{\alpha \beta^2} \Vert \hat{X} \Vert_{\frac{2}{3}}^\frac{2}{3} + \Vert \eta \Vert_2^2 } +  \Vert \eta \Vert_2 + \sqrt{\delta} \\
	&\le \sqrt{s^{\frac{1}{3}} R^{\frac{2}{3}}C_{2,1}c_{\hat U}} \sqrt[6]{\alpha \beta^2} \Vert \hat{X} \Vert_{\frac{2}{3}}^\frac{1}{3} + 2 \Vert \eta \Vert_2 + \sqrt{\delta}.
	\end{align*} 
	In the third inequality we used Proposition \ref{Bound-y2} in combination with $\Vert \hat{v}^r\Vert_1  \leq \sqrt s \Vert \hat{v}^r \Vert_2$ and
	\begin{align*}
	\sum_{r=1}^{R} \left( \Vert \hat{u}^r \Vert_2 \Vert \hat{v}^r \Vert_1 \right)^\frac{2}{3} \le s^\frac{1}{3} \sum_{r=1}^{R} \left( \Vert \hat{u}^r \Vert_2 \Vert \hat{v}^r \Vert_2 \right)^\frac{2}{3} \leq c_{\hat U} R^{\frac{2}{3}}s^\frac{1}{3} \Vert \hat{X} \Vert_\frac{2}{3}^\frac{2}{3},
	\end{align*}
    where we used again \eqref{qneq} for  $p=2/3$.
\end{Proof}

\subsection{Proof of Lemma \ref{GaussianRIP}}\label{sec:coveringnumbers}

For proving \Cref{GaussianRIP} we need bounds on the covering numbers of $S_{s_1,s_2}^{R,\Gamma}$ and $K_{s_1,s_2}^{R,\Gamma}$. The covering number $N(M,\Vert \cdot \Vert,\eps)$ of a set $M$ is the minimal number of $\Vert \cdot \Vert$-balls of radius $\eps$ that are needed to cover the set $M$ completely. The cardinality of any $\eps$-net $\tilde{M}$ of $M$, i.e., for all $z \in M$ there is $\tilde{z} \in \tilde{M}$ with $\Vert z - \tilde{z} \Vert < \eps$, yields an upper bound for $N(M,\Vert \cdot \Vert,\eps)$. The bound for $N(S_{s_1,s_2}^{R,\Gamma},\Vert \cdot \Vert_F,\eps)$ below is an adaption of Lemma 3.1 in \cite{candes2011tight} and its proof can be found in the Appendix.
%
\begin{lemma}[Covering Number for Low-Rank Matrices with Sparse Rank-$R$ Decomposition] \label{CoveringNumber}
	Let $S_{s_1,s_2}^{R,\Gamma}$ be the set defined in \eqref{S}. Then, for all $0 < \eps < 1$, one has
	\begin{align} \label{CoveringCardinality}
	\log(N(S_{s_1,s_2}^{R,\Gamma},\Vert \cdot \Vert_F,\eps)) \le R(s_1+s_2+1) \log\left( \frac{18\Gamma R}{\eps} \right) + Rs_1 \log\left( \frac{e n_1}{s_1} \right) + Rs_2 \log\left( \frac{e n_2}{s_2} \right).
	\end{align}
\end{lemma}
\paragraph{} To derive a similar bound on $N(K_{s_1,s_2}^{R,\Gamma},\Vert \cdot \Vert_F,\eps)$ we need information on the covering number of the set of effectively $s$-sparse vectors $K_{n,s} \subset \R^n$. Plan and Vershynin derived several interesting properties of $K_{n,s}$ in \cite{Plan2013LP}. Among those  \cite[Lemma 3.4]{Plan2013LP} gives the following bound for $N(K_{n,s},\Vert \cdot \Vert_2,\eps)$.

\begin{lemma}\label{effectivelySparse}
	For $0 < \eps < 1$ the covering number of $K_{n,s}$ is bounded by
	\begin{align*}
	\log N(K_{n,s},\Vert \cdot \Vert_2,\eps) \le
	\begin{cases}
	n \log \left(\frac{6}{\eps}\right) & 0 < \eps < 2\sqrt{\frac{s}{n}},\\
	\frac{4s}{\eps^2} \log \left( \frac{9\eps n}{s} \right) & 2\sqrt{\frac{s}{n}} {\leq} \eps < 1.
	\end{cases}
	\end{align*}
\end{lemma}

\begin{lemma}[Covering Number for Matrices with effectively Sparse Decomposition] \label{CoveringNumber2}
	Let $K_{s_1,s_2}^{R,\Gamma}$ be the set defined in \eqref{K}. Assume w.l.o.g.\ that $s_1/n_1 \le s_2/n_2$. \red{Then, for all $0 < \eps < 6\Gamma \sqrt{R}$, one has}
	\begin{align} \label{CoveringCardinality2}
	\log(N(K_{s_1,s_2}^{R,\Gamma}, \Vert \cdot \Vert_F, \eps)) \le \begin{cases}
	R(n_1+n_2+1) \log\left( \frac{36\Gamma R}{\eps} \right) & 0 < \eps < 12\Gamma\sqrt{\frac{Rs_1}{n_1}},\\
	\frac{144\Gamma^2 R^2s_1}{\eps^2} \log\left( \frac{9\eps n_1}{6\Gamma \sqrt{R} s_1} \right) + R(n_2+1) \log\left( \frac{36\Gamma R}{\eps} \right) & 12\Gamma\sqrt{\frac{Rs_1}{n_1}} {\le} \eps < 12\Gamma\sqrt{\frac{Rs_2}{n_2}},\\
	\frac{144\Gamma^2 R^2 (s_1 + s_2)}{\eps^2} \log\left( \frac{9\eps n_1}{6\Gamma \sqrt{R} s_1} \right) + R \log\left( \frac{18\Gamma R}{\eps} \right) & \red{12\Gamma\sqrt{\frac{Rs_2}{n_2}} {\le} \eps < 6\Gamma \sqrt{R}.}
	\end{cases}
	\end{align}
\end{lemma}
\begin{Proof}
	Let $\tilde{K}_{n,s}$ be a minimal $\eps/(6\Gamma \sqrt{R})$-net for $K_{n,s}$ in Euclidean norm. Let $D_\Gamma$ be the set of $R\times R$ diagonal matrices with Frobenius-norm less or equal $\Gamma$. It is well known that $N(D_\Gamma, \Vert \cdot \Vert_F, \eps) \le (3\Gamma/\eps)^R$. Denote by $\tilde{D}_\Gamma$ a minimal $(\eps/(6R))$-net of $D_\Gamma$ and define the sets
	\begin{align*}
	K &= \{ Z \in \R^{n_1\times n_2} \colon Z = U\Sigma V^T \text{ with } u^r \in K_{n_1,s_1}, \; v^r \in K_{n_2,s_2} \text{ for all } r \in [R], \text{ and } \Vert \Sigma \Vert_F \le \Gamma \} \\
	\tilde{K} &= \{ \tilde{Z} \in \R^{n_1\times n_2} \colon \tilde{Z} = \tilde{U}\tilde{\Sigma} \tilde{V}^T \text{ with } \tilde{u}^r \in \tilde{K}_{n_1,s_1}, \; \tilde{v}^r \in \tilde{K}_{n_2,s_2} \text{ for all } r \in [R], \text{ and } \tilde{\Sigma} \in \tilde{D}_\Gamma \}.
	\end{align*}
	We first show that $\tilde{K}$ is an $(\eps/2)$-net of $K$. Let $Z = U\Sigma V^T \in K$ be given. There exists $\tilde{Z} = \tilde{U} \tilde{\Sigma} \tilde{V}^T \in \tilde{K}$ with $\Vert u^r - \tilde{u}^r \Vert_2 \le \eps/(6\Gamma \sqrt{R})$, $\Vert v^r - \tilde{v}^r \Vert_2 \le \eps/(6\Gamma \sqrt{R})$, for all $r \in [R]$, and $\Vert \Sigma - \tilde{\Sigma} \Vert_F \le \eps/(6R)$. Therefore, $\Vert U - \tilde{U} \Vert_F^2 = \sum_{r=1}^R \Vert u^r - \tilde{u}^r \Vert_2^2 \le (\eps/(6\Gamma))^2$ and $\Vert V - \tilde{V} \Vert_F^2 \le (\eps/(6\Gamma))^2$. Moreover, $\Vert U \Vert_F^2 = \sum_{r=1}^{R} \Vert u^r \Vert_2^2 \le R$ (the same holds for $V,\tilde{U},\tilde{V}$) and $\Vert U\Sigma \Vert_F \le \Vert \Sigma \Vert_F$ (the same holds for $\Sigma V^T,\tilde{U} \Sigma,\Sigma \tilde{V}^T$). We now obtain by the triangle inequality and the fact that $\Vert AB \Vert_F \le \Vert A \Vert_F \Vert B \Vert_F$
	\begin{align*}
		\Vert Z - \tilde{Z} \Vert_F &\le \Vert (U-\tilde{U})\Sigma V^T \Vert_F + \Vert \tilde{U} (\Sigma - \tilde{\Sigma}) V^T \Vert_F + \Vert \tilde{U} \tilde{\Sigma} (V - \tilde{V})^T \Vert_F \\
		&\le \frac{\eps}{6\Gamma} \Gamma + \sqrt{R} \frac{\eps}{6R} \sqrt{R} + \Gamma \frac{\eps}{6\Gamma} \le \frac{\eps}{2}.
	\end{align*}
	Since $K_{s_1,s_2}^{R,\Gamma} \subset K$ one has $N(K_{s_1,s_2}^{R,\Gamma},\Vert \cdot \Vert_F, \eps) \le N(K,\Vert \cdot \Vert_F, \eps/2)$. Hence,
	\begin{align*}
		N(K_{s_1,s_2}^{R,\Gamma},\Vert \cdot \Vert_F, \eps) \le |\tilde{K}| \le |\tilde{K}_{n_1,s_1}|^R |\tilde{D}_\Gamma| |\tilde{K}_{n_2,s_2}|^R
	\end{align*}
	which yields the claim by applying \Cref{effectivelySparse}.
\end{Proof}

Lemma \ref{GaussianRIP} can be proven by applying the following bound on suprema of chaos processes \cite[Theorems 1.4 \& 3.1]{krahmer2014suprema} in combination with the {bounds on the covering numbers} $N(S,\Vert \cdot \Vert_F,\eps)$ and $N(K,\Vert \cdot \Vert_F,\eps)$ of $S$ and $K$ of \Cref{CoveringNumber} and \Cref{CoveringNumber2}. {We recall below the relevant result in the form presented in \cite{JungKrahmerStoeger2017}.} The appearing $\gamma_2$-functional is defined in \cite{krahmer2014suprema} and can be bounded by
\begin{equation}\label{dudley}
\gamma_2 \left( \H,  \Vert \cdot \Vert_{2 \rightarrow 2} \right) \lesssim  \int_{0}^{d_{2 \rightarrow 2} \left( \H \right) } \sqrt{ \log N \left( \H , \Vert \cdot \Vert_{2\rightarrow 2} ,\varepsilon \right) } d\varepsilon,
\end{equation}
in the case of a set of matrices $\H$ equipped with the operator norm. {Here and below $d_\boxdot(\H) = \sup_{H \in \H} \Vert H \Vert_\boxdot$, 
where $\boxdot$ is a generic norm.}

\begin{theorem} \label{KMR} 
	Let $ \mathcal{H} $ be a  symmetric set of matrices, i.e., $ \mathcal{H} = - \mathcal{H} $, and let $ \xi $ be a random vector whose entries $\xi_i$ are independent $\mathcal{K}$-subgaussian random variables with mean $0$ and variance $1$. Set
	\begin{align*}
	E&= \gamma_2 \left( \mathcal{H}, \Vert \cdot \Vert_{2\rightarrow 2}  \right) \left(  \gamma_2 \left( \mathcal{H}, \Vert \cdot \Vert_{2\rightarrow 2}  \right) + d_F (\mathcal{H}) \right) \\
	V&= d_{2 \rightarrow 2} \left( \mathcal{H} \right) \left( \gamma_2 \left( \mathcal{H}, \Vert \cdot \Vert_{2\rightarrow 2}  \right) +   d_F (\mathcal{H}) \right)\\
	U&= d^2_{2 \rightarrow 2} \left( \mathcal{H} \right) 
	\end{align*}
	Then, for $t>0$, 
	\begin{equation*}
	\Pr[]{\underset{H \in \mathcal{H}}{\sup} \big\vert \Vert H \xi \Vert_{\ell_2}^2 -  \E[]{\Vert H\xi \Vert_{2}^2}  \big\vert \ge c_1 E + t } \le 2 \exp \left( -c_2 \min \left( \frac{t^2}{V^2}, \frac{t}{U} \right) \right).
	\end{equation*}
	The constants $c_1$ and $c_2 $ are universal and only depend on $\mathcal{K}$.
\end{theorem}
{We refer the reader to \cite{krahmer2014suprema}  and \cite{JungKrahmerStoeger2017} for further details.}

\begin{Proof}[of Lemma \ref{GaussianRIP}]
	The proof consists of three main parts. We start in \textbf{(I)} by fitting our setting into the one of Theorem \ref{KMR}. In \textbf{(IIa)} resp. \textbf{(IIb)} the $\gamma_2$-functional gets bounded for $S_{s_1,s_2}^{R,\Gamma}$ and $K_{s_1,s_2}^{R,\Gamma}$, and in \textbf{(III)} we conclude by applying Theorem \ref{KMR}.
	\paragraph{(I)} We first switch the roles of our random measurement operator ${\A}$ applied to the fixed matrices $Z$ to have fixed operators $H_Z$ applied to a random vector $\xi$. Denote by $\vect(Z) \in \R^{n_1n_2}$ the vectorization of $Z$. Observe, for all $Z\in \R^{n_1\times n_2}$,
	\begin{align*}
	{\A}(Z) = \frac{1}{\sqrt{m}} \begin{pmatrix}
	\langle \vect(A_1),\vect(Z) \rangle \\
	\vdots \\
	\langle \vect(A_m),\vect(Z) \rangle
	\end{pmatrix} = \frac{1}{\sqrt{m}} \begin{pmatrix}
	\vect(Z)^T & 0 & \cdots \\
	&\ddots & \\
	\cdots & 0 & \vect(Z)^T
	\end{pmatrix} \cdot \begin{pmatrix}
	\vect(A_1) \\
	\vdots \\
	\vect(A_m)
	\end{pmatrix} = H_Z \cdot \xi
	\end{align*}
	where $H_Z \in \R^{m\times mn_1n_2}$ is a matrix depending on $Z$ and $\xi \in \R^{mn_1n_2}$ has i.i.d.\ $\mathcal{K}$-subgaussian entries $\xi_l$ of mean $0$ and variance $1$. We define $\H_S = \{ H_Z \colon Z \in S_{s_1,s_2}^{R,\Gamma} \}$. Note that the mapping $Z \mapsto H_Z$ is an isometric linear bijection. In particular, we have $\Vert H_Z \Vert_F = \Vert Z \Vert_F$ and $\Vert H_Z \Vert_{2\rightarrow 2} = \Vert Z \Vert_F / \sqrt{m}$. For $Z \in S_{s_1,s_2}^{R,\Gamma}$ it holds that $\Vert Z \Vert_F \le \Vert U \Vert_F \Vert \Sigma V^T \Vert_F \le \Gamma \sqrt{R}$. Hence, $d_F(\H_S) \le \Gamma \sqrt{R}$ and $d_{2\rightarrow 2}(\H_S) \le \Gamma \sqrt{R}/\sqrt{m}$. 
	\paragraph{(IIa)} Since $\Vert H_Z \Vert_{2\rightarrow 2} = \Vert Z \Vert_F / \sqrt{m}$ and $Z \mapsto H_Z$ is a linear bijection, it follows that $N(\H_S,\Vert \cdot \Vert_{2\rightarrow 2},\eps) = N(S,\Vert \cdot \Vert_F,\sqrt{m} \eps)$. We can estimate by \eqref{dudley} and \Cref{GammaBounds}
	\begin{align*}
	\gamma_2 \left( \H_S,  \Vert \cdot \Vert_{2 \rightarrow 2} \right) &\lesssim  \int_{0}^{\frac{\Gamma\sqrt{R}}{\sqrt{m}}} \sqrt{ \log N \left( \H_S, \Vert \cdot \Vert_{2\rightarrow 2} ,\varepsilon \right) } d\varepsilon = \int_{0}^{\frac{\Gamma \sqrt{R}}{\sqrt{m}}} \sqrt{ \log N \left( S_{s_1,s_2}^{R,\Gamma} , \Vert \cdot \Vert_F ,\sqrt{m} \varepsilon \right) } d\varepsilon \\
	&\le \sqrt{\frac{C_S \Gamma^2 R^2 (s_1+s_2) \log \left( \max \left\{ n_1,n_2 \right\} \right) }{m}} =: \mathcal{L}_S.
	\end{align*}
	for some constant $C_S > 0$.
	\paragraph{(IIb)} In the same manner we obtain a bound on $\gamma_2 (\H_K,\Vert \cdot \Vert_{2\rightarrow 2})$ where $\H_K = \{ H_Z \colon Z \in K_{s_1,s_2}^{R,\Gamma} \}$. Recall that $\Vert H_Z \Vert_F = \Vert Z \Vert_F$, $\Vert H_Z \Vert_{2\rightarrow 2} = \Vert Z \Vert_F/\sqrt{m}$ and $Z \mapsto H_Z$ is an linear bijection. This implies $N(\H_K,\Vert \cdot \Vert_{2\rightarrow 2},\eps) = N(K_{s_1,s_2}^{R,\Gamma},\Vert \cdot \Vert_F,\sqrt{m}\eps)$. Note that $d_F(\H_K) \le \Gamma \sqrt{R}$ and $d_{2\rightarrow 2}(\H_K) \le \Gamma \sqrt{R}/\sqrt{m}$. We obtain by \eqref{dudley} and \Cref{GammaBounds}
	\begin{align*}
	\gamma_2(\H_K,\Vert \cdot \Vert_{2\rightarrow 2}) &\lesssim \int_{0}^{\frac{\Gamma \sqrt{R}}{\sqrt{m}}} \sqrt{\log N(\H_K,\Vert \cdot \Vert_{2\rightarrow 2},\eps)} \;d\eps = \int_{0}^{\frac{\Gamma \sqrt{R}}{\sqrt{m}}} \sqrt{\log N(K_{s_1,s_2}^{R,\Gamma},\Vert \cdot \Vert_F,\sqrt{m}\eps)} \;d\eps \\
	&\le \sqrt{\frac{C_K \Gamma^2 R^2 (s_1 + s_2) \log^3(\max\{n_1,n_2\})}{m}} =: \mathcal{L}_K
	\end{align*}
	for some constant $C_K > 0$.
	\paragraph{(III)} The final part of the proof is now equal for both sets $S_{s_1,s_2}^{R,\Gamma}$ and $K_{s_1,s_2}^{R,\Gamma}$. We write $\mathcal{L}$ for $\mathcal{L}_S$ resp. $\mathcal{L}_K$ and assume $m \gtrsim C_S \Delta^{-2} R (s_1+s_2) \log \left( \max \left\{ n_1,n_2 \right\} \right)$ resp. $m \gtrsim C_K \Delta^{-2} R (s_1 + s_2) \log^3(\max\{n_1,n_2\})$, for some $0 < \Delta < 1$. Then, $\mathcal{L} \le \Gamma \sqrt{R}$ and 
		\begin{align} \label{eq:2}
			\mathcal{L}^2 + \Gamma \sqrt{R} \mathcal{L} \le \Gamma^2 R (\Delta^2 + \Delta) \le 2\Gamma^2 R \Delta.
	\end{align}
	We obtain the following bounds on the quantities (cf. \Cref{KMR}):
	\begin{align} \label{ParameterBound}
		E \le \mathcal{L}^2 + \Gamma \sqrt{R} \mathcal{L}, \tab\tab V \le \frac{\Gamma \sqrt{R} \mathcal{L} + \Gamma^2 R}{\sqrt{m}}, \tab\tab U \le \frac{\Gamma^2 R}{m}.
	\end{align}
	Observing now that $\E[]{\Vert H_Z \xi \Vert_2^2} = \Vert H_Z \Vert_F^2 = \Vert Z \Vert_F^2$ and recalling $\Gamma \ge 1$ we finally get, for $\delta \ge 3 c_1 \Gamma^2 R \Delta$ (which implies by \eqref{eq:2} that $\delta \ge c_1 E + c_1 \Gamma^2 R \Delta$),
	\begin{align*}
		\Pr[]{\sup_{Z \in S} \left| \Vert \A(Z) \Vert_2^2 - \Vert Z \Vert_F^2 \right| \ge \delta} &\le \Pr[]{\sup_{H_Z \in \H} \left| \Vert H_Z \xi \Vert_2^2 - \E[]{\Vert H_Z \xi \Vert_2^2} \right| \ge c_1 E + c_1 \Gamma^2 R \Delta } \\
		&\le 2\exp \left( -c_2 \min \left\{ m \frac{c_1^2 \Gamma^4 R^2 \Delta^2}{\Gamma^2R(\mathcal{L} + \Gamma \sqrt{R})^2}, m  \frac{c_1 \Gamma^2 R \Delta}{\Gamma^2 R} \right\}  \right) \\ 
		&\le 2\exp \left( -C \Delta^2 m \right)
	\end{align*}
	where $C > 0$ is a positive constant which depends on $\mathcal{K}$. In the last step we used that $\mathcal{L} + \Gamma \sqrt{R} \in [ \Gamma\sqrt{R},2\Gamma\sqrt{R} ]$ (because $0 < \mathcal{L} < \Gamma \sqrt{R}$).
\end{Proof}
\begin{remark}\label{Gammacond}
The condition $1 \leq \Gamma \leq \Gamma^2$ is used in a crucial way in part (III) of the proof above. Additionally, without condition $\Gamma\geq 1$, any $\hat X \in S_{s_1,s_2}^{R,\Gamma}$ could be decomposed also as follows
$$
\hat X=  \sum_{r=1}^{R} \sigma_r u^r (v^r)^T =\sum_{r=1}^{R} \sum_{j=1}^K  \frac{\sigma_r}{K} u^r (v^r)^T,  \quad \left ( \sum_{r=1}^{R} \sum_{j=1}^K  \frac{\sigma_r^2}{K^2} \right)^{1/2}  \leq \frac{\Gamma}{\sqrt K},
$$
for any $K \in \mathbb N$, implying $\hat X \in S_{s_1,s_2}^{R K,\Gamma/\sqrt{K}}$ as well, which would result in a larger number of necessary measurements $m \gtrsim \left( \frac{\delta}{\Gamma^2 R} \right)^{-2} R K (s_1+s_2) \log^3 \left( \max \{n_1,n_2\} \right)$. Hence, $\Gamma \geq 1$ emerges as a natural condition, in order to have a correct proof of part (III) and to avoid ambiguities on the necessary measurements.

\end{remark}


\section{Implementation and Numerical Experiments} \label{Numerics}

After having obtained some theoretical insight on the proposed optimization problem, we \purple{provide an implementation of \eqref{ATLAS0} and discuss its predicted behavior in numerical experiments}. Therefore, we begin by presenting the implementation that has been used in all experiments\footnote{The corresponding Matlab code is provided at https://www-m15.ma.tum.de/Allgemeines/SoftwareSite}. As in practice ATLAS converges even without the auxiliary terms $\frac{1}{2\lambda_k^\ell} \Vert u - u_k^\ell \Vert^2_2$ and $\frac{1}{2\lambda_k^\ell} \Vert v - v_k^\ell \Vert^2_2$ introduced in \eqref{ATLAS0}, for sake of simplicity we drop those terms. By the alternating form of \eqref{ATLAS0} one has to solve a certain number of Tikhonov regularization resp.\ $\ell_1$-LASSO problems. \purple{Note that  for the Tikhonov regularization}
\begin{align*}
	u = \argmin_{z \in \R^n} \Vert y - Az \Vert_2^2 + \alpha \Vert z \Vert_2^2,
\end{align*}
\purple{with $A \in \R^{m\times n}, y \in \R^m$, and $\alpha > 0$,} the solution is explicitly given by $u = (\alpha\id + A^TA)^{-1} A^Ty$. Solutions to $\ell_1$-LASSO
\begin{align*}
    v = \argmin_{z \in \R^n} \Vert y - Az \Vert_2^2 + \beta \Vert z \Vert_1,
\end{align*}
for some $A \in \R^{m\times n}, y \in \R^m$ and $\beta > 0$ can be well approximated by the so-called Iterative Soft-Thresholding Algorithm (ISTA) which is based on the soft-thresholding operator $\St_\beta$
\begin{align*} 
	\St_\beta (z) = \begin{pmatrix}
		S_\beta (z_1) \\ \vdots \\ S_\beta (z_{n_2}) 
	\end{pmatrix}, \tab
	\text{ where } S_\beta (z_i) = \begin{cases}
		z_i - \frac{\beta}{2} & z_i > \frac{\beta}{2} \\
		0 & |z_i| \le \frac{\beta}{2} \\
		z_i + \frac{\beta}{2} & z_i < -\frac{\beta}{2}
	\end{cases}.
\end{align*}
Hence, a suitable implementation of \eqref{ATLAS0} is given by Algorithm \ref{A1}, \purple{whereas Algorithm \ref{A2} describes ISTA for the reader's convenience.} Necessary modifications in case of sparse left {component vectors of $\hat{X}$ are rather straightforward}. \blue{We generate different ground-truths $\hat{X}$ at random by first fixing $n_1,n_2$, and $s$ and then uniformly at random drawing $R$ pairs of unit norm vectors $\hat{u}^r \in \R^{n_1}, \hat{v}^r \in \R^{n_2}$ and a Gaussian vector $\hat{\sigma} \in \R^R$ (for the right components, we first choose a support of size $s$ uniformly at random and only then fill the non-zero positions with a randomly drawn $s$-dimensional unit norm vector). The ground-truth $\hat{X}$ is obtained as $\hat{U} \hat{\Sigma} \hat{V}^T$ and re-normalized to a given value. To have a fair comparison to SPF which has been designed for orthogonal decompositions, we orthogonalize the component vectors $\hat{u}^r,\hat{v}^r$ without changing their support before composing $\hat{X}$. Experiments, however, showed that ATLAS performs in a similar way without this additional step.} 

Let us turn toward numerical simulations. First, we check if the main theoretical results stated in Theorem \ref{ApproxX} and Corollary \ref{ApproxXcor} describe the qualitative and quantitative behavior of the approximation error well. Then, we compare ATLAS to the already mentioned Sparse Power Factorization (SPF), \cite{lee2013near}. {We used the leading singular vectors of $\A^\ast(y)$ to initialize both algorithms, which is likely not an optimal choice and certainly may cause loss of performance for both algorithms, but it is nevertheless sufficient to illustrate certain comparisons numerically.}

\begin{algorithm} 
	\caption{\textbf{:}  \textbf{ATLAS}$(y,A,R,v_0^1,...,v_0^R,\alpha,\beta)$} \label{A1}
	\begin{algorithmic}[1]
		\Require{$y \in \R^m$, $A \in \R^{m\times n_1n_2}$, rank $R$, $v_0^1,...,v_0^R \in \R^{n_2}$ and $\alpha, \beta > 0$}
		\Statex
		\While{stop condition is not satisfied}
		\Let{$u_k$}{$\bigl(\alpha\id + \A_v(v_{k-1})^T \A_v(v_{k-1})\bigr)^{-1} \A_v(v_{k-1})^T y$}
		\Comment{$\A(uv^T) = \A_v(v)\cdot u$}
		\Let{$v_k$}{\textbf{ISTA}$(y,\A_u(u_k),v_{k-1},\beta)$}
		\Comment{$\A(uv^T) = \A_u(u)\cdot v$}
		\EndWhile
		\State 
		\Return{$u_\text{final}^1,...,v_\text{final}^R$}
	\end{algorithmic}
\end{algorithm}

\begin{algorithm}
	\caption{\textbf{:} \textbf{ISTA}$(y,A,v_0,\beta)$} \label{A2}
	\begin{algorithmic}[1]
		\Require{$y \in \R^m$, $A \in \R^{m\times n}$, $v_0 \in \R^n$ and $\beta > 0$ }
		\Statex
		\While{stop condition is not satisfied}
		\Let{$v_k$}{$\St_\beta \left[ v_{k-1} + A^T (y - Av_{k-1}) \right]$}
		\EndWhile
		\State 
		\Return{$v_\text{final}$}
	\end{algorithmic}
\end{algorithm}

\subsection{Validation of Corollary \ref{ApproxXcor}} \label{Numerics1}

\begin{figure}[]
	\centering
	\captionsetup{width=.8\linewidth}
	\begin{subfigure}[b]{0.45\textwidth}
		\includegraphics[width=\textwidth]{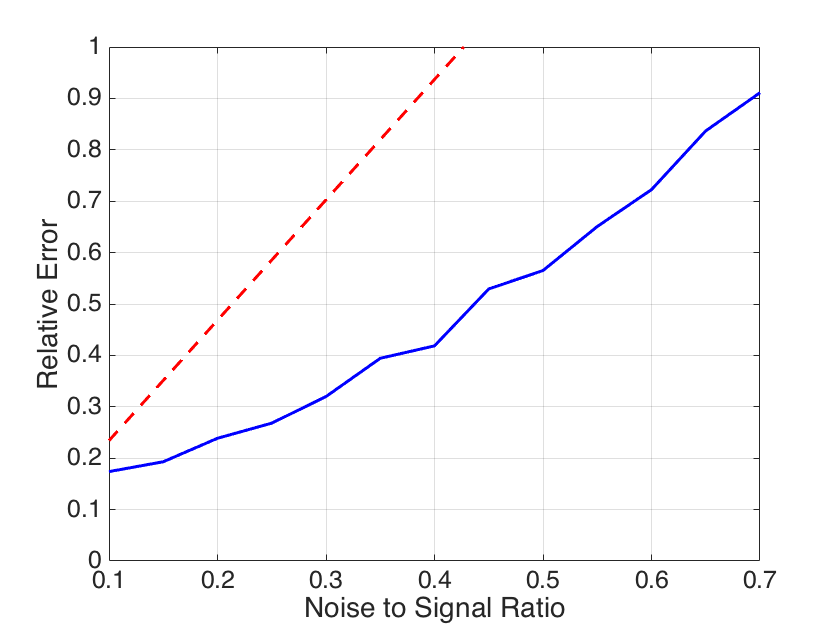}
		\caption{$R=1$}
		\label{FIG:a}
	\end{subfigure}
	\quad 
	\begin{subfigure}[b]{0.45\textwidth}
		\includegraphics[width=\textwidth]{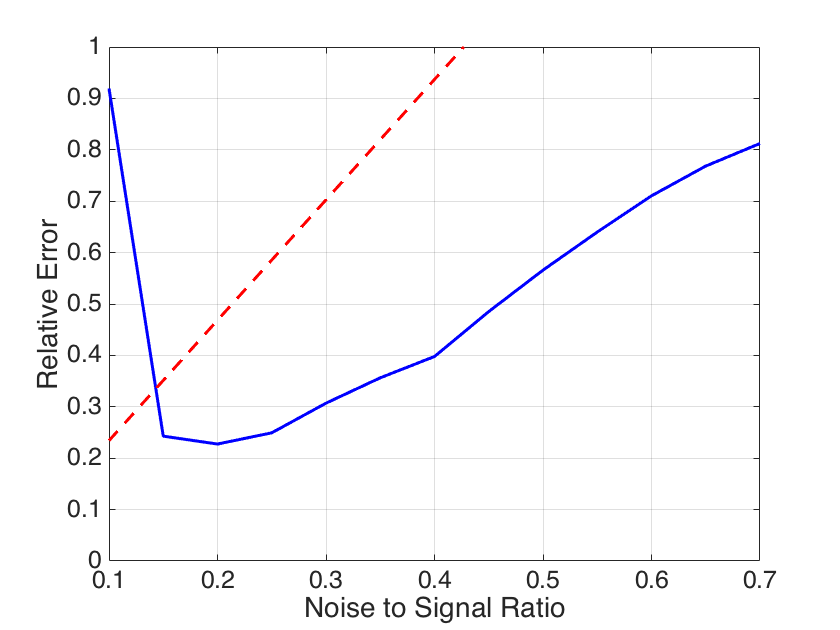}
		\caption{$R=5$}
		\label{FIG:b}
	\end{subfigure}
	\caption{Approximation quality depending on noise level (see Section \ref{Numerics1}). The $x$-axis shows noise to signal ratio $\Vert \eta \Vert_2/\Vert \hat{X} \Vert_F$ while the $y$-axis presents approximation error relative to $\Vert \hat{X} \Vert_F$. One can see the comparison of approximation results (solid blue) and theoretical bound (dashed red)}\label{Fig:NoiseTest}
\end{figure}

Figure \ref{Fig:NoiseTest} shows the average approximation error of $100$ randomly drawn $\hat{X} \in \R^{16 \times 100}$, $\Vert \hat{X} \Vert_F = 10$, with $\rank(\hat{X}) = 1$ (resp. $\rank(\hat{X}) = 5$) and $10$-sparse right singular vector(s) from $m = 90$ (resp. $m = 400$) noisy measurements $y = \A(\hat{X})+\eta$. The parameters have been chosen exemplarily for purpose of illustration. {The operator $\A$ is drawn once at random.} The error bound from Corollary \ref{ApproxXcor} is plotted as dashed red line, \purple{whereas the average approximation errors are in blue}. Though not tight the theoretical bound seems to describe the linear dependence of the approximation error on noise level {appropriately}. In addition, Figure \ref{Fig:NoiseTest} (b) shows a breakdown of approximation for noise to signal ratios below $\approx 0.25$. This occurrence is not surprising as the assumptions of Corollary \ref{ApproxXcor} include a lower-bound on the {noise-to-signal} ratio for a fixed number of measurements \blue{(if the noise-to-signal ratio becomes small, the parameters $\alpha = \beta = \Vert \eta \Vert_2^2/\Vert \hat{X} \Vert_{\frac{2}{3}}^{\frac{2}{3}}$ have to be chosen so small that the regularization weakens and the RIP requirements are harder to fulfill)}. Below a certain value the RIP requirements will be too strong for $\A$ {to fulfill it}, the RIP {breaks down,} and the recovery guarantees fail.\\
\blue{Since it is essential for practical purposes to know whether the approximation computed by ATLAS is reliable, a central question is how to judge in which of the two regimes (noise sufficiently large vs.\ noise too small) one currently is. The (effective) sparsity of the computed approximation is a good indicator for this. As long as the noise level resp. the parameters $\alpha$ and $\beta$ are sufficiently large, the approximation's (effective) sparsity will remain bounded while it explodes as soon as a critical value is crossed, cf.\ Figure \ref{Fig:Parameter} in Section \ref{Numerics2}. If the real noise level is below this threshold, one can just overestimate it by a value slightly above the threshold and choose corresponding $\alpha$ and $\beta$. For instance, the experiment in Section \ref{Numerics2} shows that even in the case of vanishing noise approximation works well up to an accuracy of $0.01 \|\hat{X}\|_F$.}

\subsection{Validation of Theorem \ref{ApproxX}} \label{Numerics2}
\begin{figure}[]
	\centering
	\captionsetup{width=.8\linewidth}
	\begin{subfigure}[b]{0.45\textwidth}
		\includegraphics[width=\textwidth]{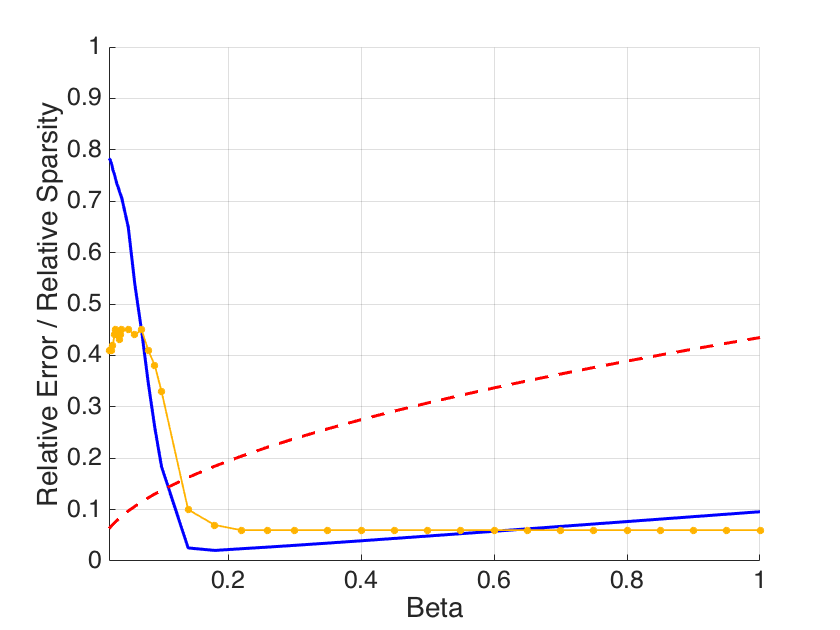}
		\caption{$\alpha = \beta$}
	\end{subfigure}
	\\ 
	\begin{subfigure}[b]{0.45\textwidth}
		\includegraphics[width=\textwidth]{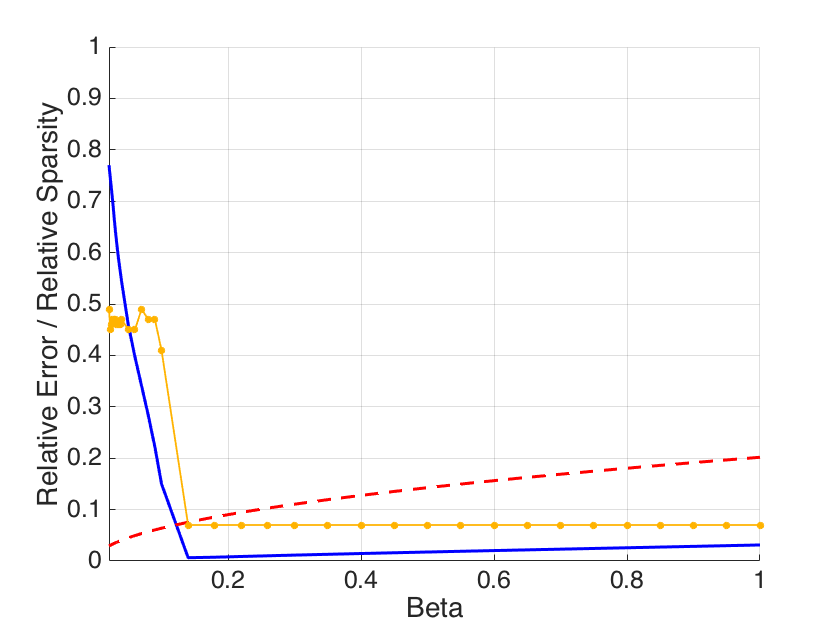}
		\caption{$\alpha = 0.01\beta$}
	\end{subfigure}
	\quad 
	\begin{subfigure}[b]{0.45\textwidth}
		\includegraphics[width=\textwidth]{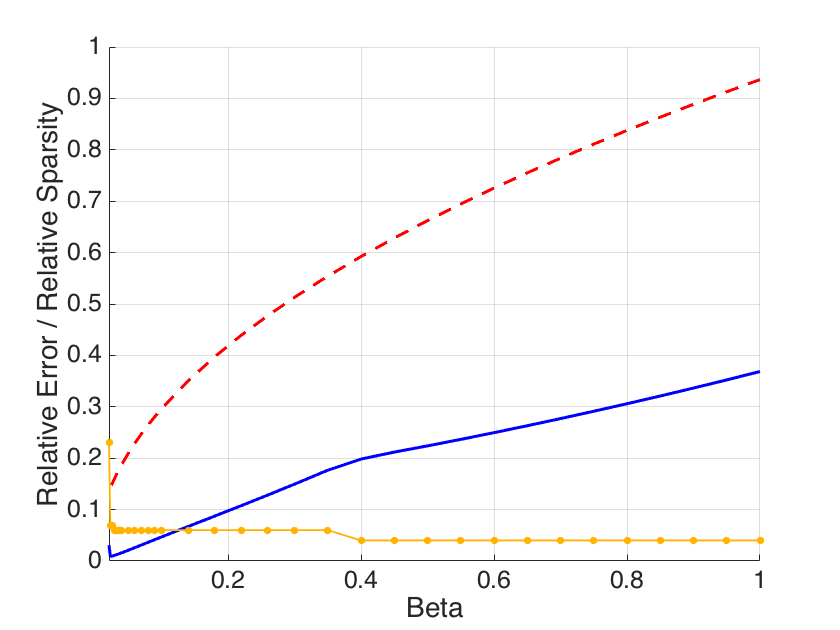}
		\caption{$\alpha = 100\beta$}
	\end{subfigure}
	\caption{Approximation quality and sparsity depending on parameter size (see Section \ref{Numerics2}). The approximation error (solid blue) and the theoretical bound (dashed red) are measured relative to $\Vert \hat{X} \Vert_F$ while sparsity of the right singular vector (dotted yellow) is  relative to $n_2$.}\label{Fig:Parameter}
\end{figure}

\purple{In the second experiment, we study the influence of parameters $\alpha$ and $\beta$ on the reconstruction accuracy. In particular, we vary} the parameters $\alpha$ and $\beta$ when reconstructing one randomly drawn $\hat{X} \in \R^{16 \times 100}$, $\Vert \hat{X} \Vert_F = 10$, with $\rank(\hat{X}) = 1$ and $10$-sparse right singular vector {from $90$ measurements without noise}. Again parameter choice is exemplary. We compare the three settings: (a) $\alpha = \beta$, (b) $\alpha = 0.01\beta$ and (c) $\alpha = 100\beta$ in Figure \ref{Fig:Parameter}. One can observe a decrease of approximation error for $\alpha,\beta \rightarrow 0$ up to a certain threshold, {under which  the approximation} seemingly fails. While this threshold lies at $\beta \approx 0.15$ in (a) and (b) it is hardly recognizable in (c). At the same time (a) and (b) show a much {smaller approximation error}. These observations suggest that the choice of $\alpha$ strongly influences the approximation quality of ATLAS. This  is consistent with Theorem \ref{ApproxX}, as a smaller $\alpha$ leads to a smaller theoretical approximation error bound.\\
Even though (a) and (b) show a linear decrease in approximation error which is in contrast to the square-root behavior of the theoretical bound, (c) suggests that the error, indeed, behaves similar to the theoretical bound.\\
\purple{Figure \ref{Fig:Parameter} shows that} the sparsity level remains stable for sufficiently large $\beta$ and breaks down precisely at the same threshold as the approximation error,  coinciding with the \purple{violation} of the RIP \purple{conditions}.

\begin{figure}[]
	\centering
	\captionsetup{width=.8\linewidth}
	\includegraphics[width=0.45\textwidth]{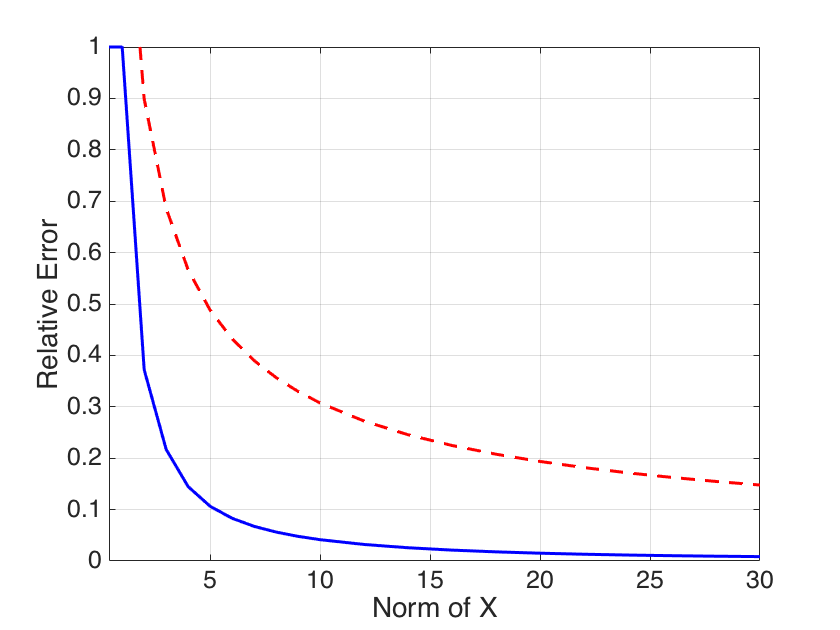}
	\caption{{Approximation error depending on the magnitude of $\hat{X}$ {in Frobenius norm} (see Section \ref{Numerics2}). Approximation error (solid blue) and theoretical bound (dashed red) are relative to $\Vert \hat{X} \Vert_F$.}}\label{Fig:NormTest}
\end{figure}


For a better understanding of ATLAS we made a third experiment reconstructing one randomly drawn $\hat{X} \in \R^{16 \times 100}$ with $\rank(\hat{X}) = 1$ and $10$-sparse right singular vector for different values of $\Vert \hat{X} \Vert_F$ {from $90$ measurements. The noise level was set to $0$ and the parameters to $\alpha = \beta = 0.5$.} The outcome is depicted in Figure \ref{Fig:NormTest}. One can see that the relative approximation error {decreasing with the magnitude of $\hat{X}$ as expected from the bound of Theorem \ref{ApproxX}.} This seemingly confirms the theoretical dependence of reconstruction error on $\Vert \hat{X} \Vert_\frac{2}{3}^\frac{1}{3}$.

\subsection{ATLAS vs SPF} \label{Numerics3}

After confirming the theoretical results numerically, we now turn to the comparison of  ATLAS  with its state-of-the-art counterpart SPF  \cite{lee2013near}. To our knowledge, SPF is the only algorithm available so far in matrix sensing, which exploits low-rankness and sparsity constraints together and comes with near-optimal recovery guarantees (not relying on a special structure of $\A$ as in \cite{bahmani2016near}). As \cite{lee2013near} contains exhaustive numerical comparisons of SPF and {low-rank (resp. sparse) recovery strategies based on convex relaxation}, SPF suffices for numerical benchmark tests. From the structure of the algorithms and their respective theoretical analysis one would expect SPF to yield more accurate reconstruction in the noiseless-to-low-noise setting, while ATLAS should prove to be more reliable if noise becomes large. This theoretical expectation is confirmed by the following experiments.


\begin{figure}[!htb]
	\centering
	\captionsetup{width=.9\linewidth}
	\begin{subfigure}[b]{0.45\textwidth}
		\includegraphics[width=\textwidth]{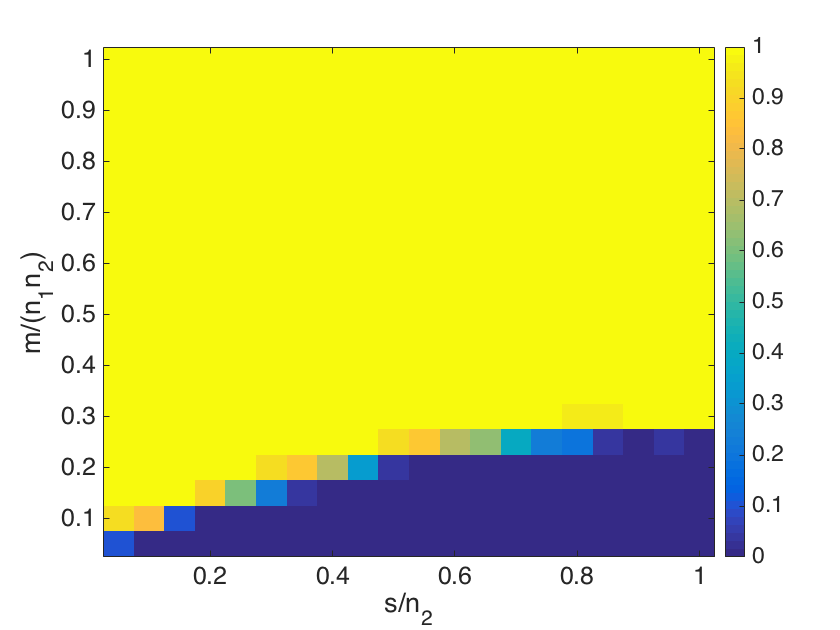}
		\caption{SPF, no noise}
		\label{fig:a}
	\end{subfigure}
	\quad 
	\begin{subfigure}[b]{0.45\textwidth}
		\includegraphics[width=\textwidth]{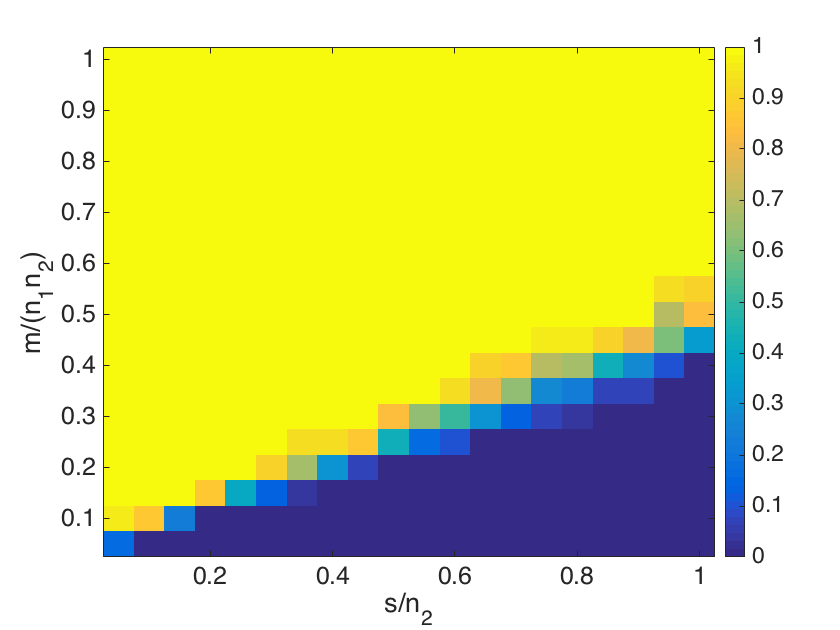}
		\caption{ATLAS, no noise}
		\label{fig:b}
	\end{subfigure}
	\\ 
	\begin{subfigure}[b]{0.45\textwidth}
		\includegraphics[width=\textwidth]{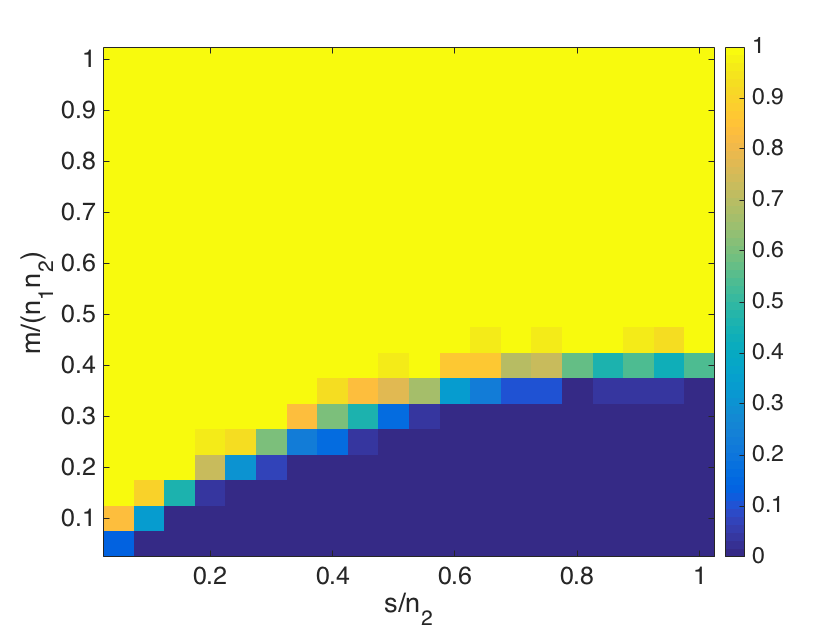}
		\caption{SPF, with relatively strong noise}
		\label{fig:c}
	\end{subfigure}
	\quad
	\begin{subfigure}[b]{0.45\textwidth}
		\includegraphics[width=\textwidth]{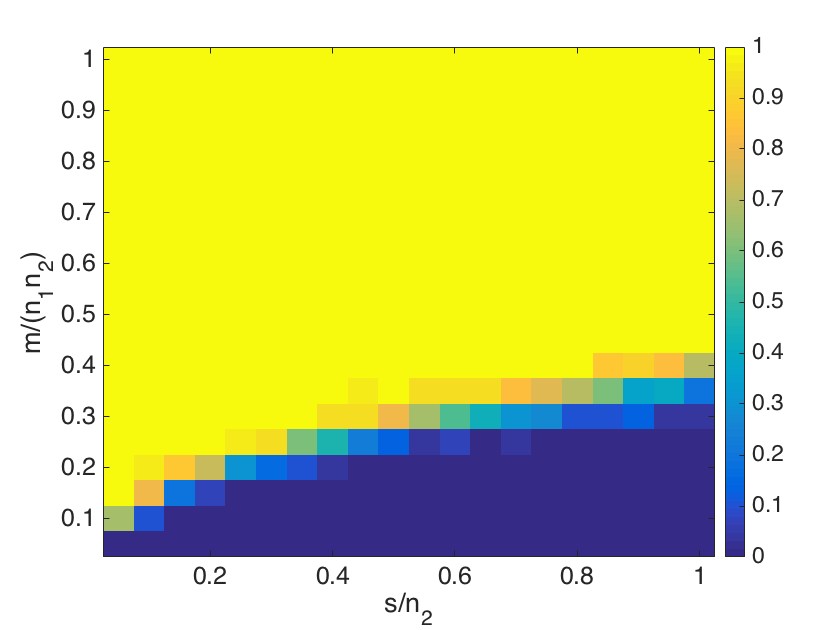}
		\caption{ATLAS, with relatively strong noise}
		\label{fig:d}
	\end{subfigure}
	\caption{Phase transition diagrams comparing SPF and ATLAS with and without noise on the measurements (see Section \ref{Numerics3}). \blue{Empirical recovery probability, i.e., percentage of successful reconstructions, is depicted by color from zero (blue) to one (yellow).}}\label{fig:Greyscale}
\end{figure}

\begin{figure}[!htb]
	\centering
	\captionsetup{width=.8\linewidth}
	\begin{subfigure}[b]{0.45\textwidth}
		\includegraphics[width=\textwidth]{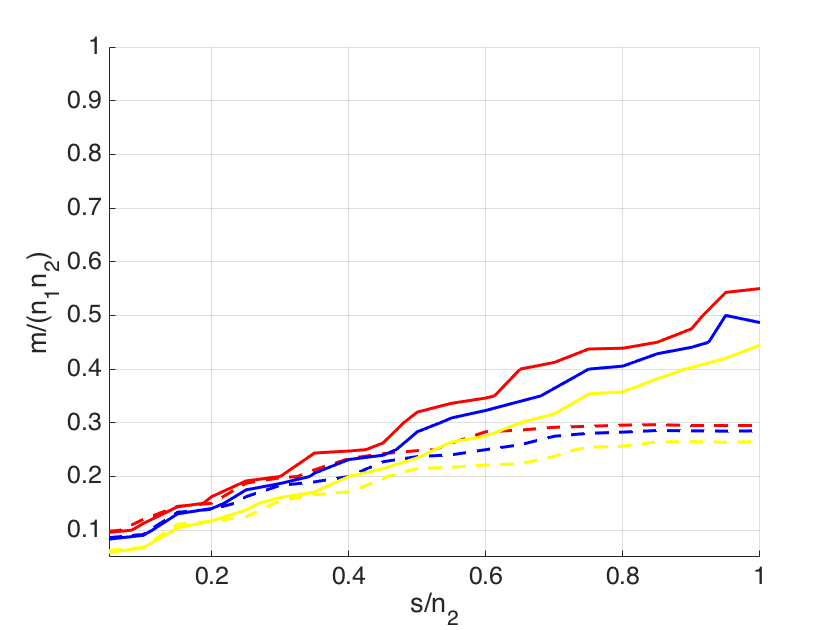}
		\caption{No noise}
	\end{subfigure}
	\quad 
	\begin{subfigure}[b]{0.45\textwidth}
		\includegraphics[width=\textwidth]{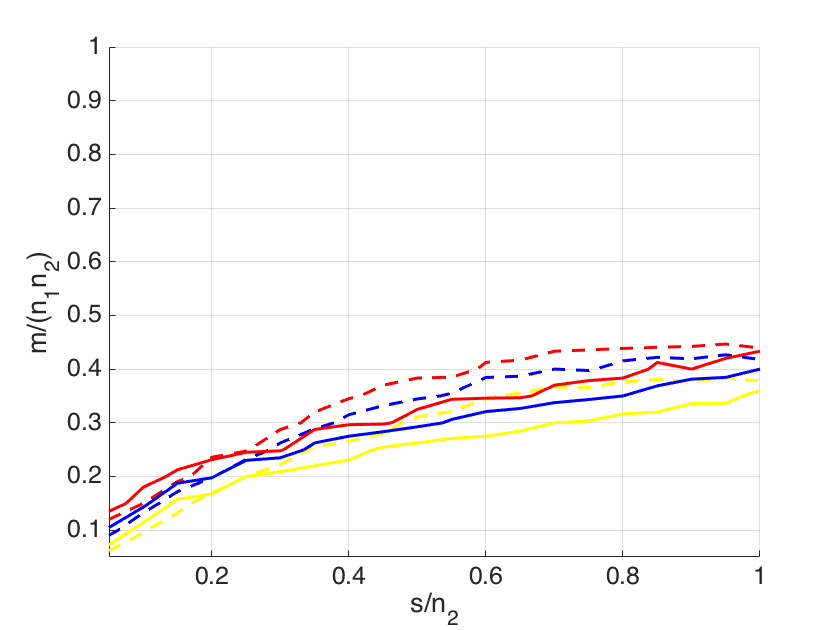}
		\caption{Noise}
	\end{subfigure}
	\caption{Recovery probability comparison of SPF (dashed) and ATLAS (solid). Plotted are the thresholds for $90\%$ (red), $70\%$ (blue) and $30\%$ (yellow) successful recoveries. A recovery was counted successful if $\Vert \hat{X} - X_\text{appr} \Vert_F/\Vert \hat{X} \Vert_F \le 0.2$ (resp. $0.4$)}\label{Fig:Greyscale}
\end{figure}

In Figure \ref{fig:Greyscale} we compare for $s/n_2 \in [0,1]$ and $m/(n_1n_2)$ the number of successful recoveries of $30$ randomly drawn $\hat{X} \in \R^{4 \times 128}$, $\Vert \hat{X} \Vert_F = 10$, with $\rank(\hat{X}) = 1$ and $s$-sparse right singular vectors from $m$ measurements. {The dimensions of $\hat{X}$ were chosen accordingly to similar experiments in \cite{lee2013near}.} We set the noise level to $0$ (resp. $0.3\Vert \hat{X} \Vert_F$) and counted the recovery successful if $\Vert \hat{X} - X_\text{appr} \Vert_F/\Vert \hat{X} \Vert_F \le 0.2$ (resp. $0.4$).  {In order to compare the noisy and noiseless cases, we fix $\alpha = \beta = 0.5$ for both, which is a reasonable choice for high noise level, but perhaps  sub-optimal if the noise level is low.} Selected quantiles are directly compared in Figure \ref{Fig:Greyscale} for convenience.\\
{As expected,} SPF outperforms ATLAS if there is no noise. In case of strong noise on the measurements, \purple{the situation changes}. In \purple{particular, we observe the improved performance of ATLAS, whereas the SPF performance remarkably deteriorates.}


\begin{figure}[!htb]
 	\centering
 	\captionsetup{width=0.9\linewidth}
 	\begin{subfigure}[b]{0.45\textwidth}
 		\includegraphics[width=\textwidth]{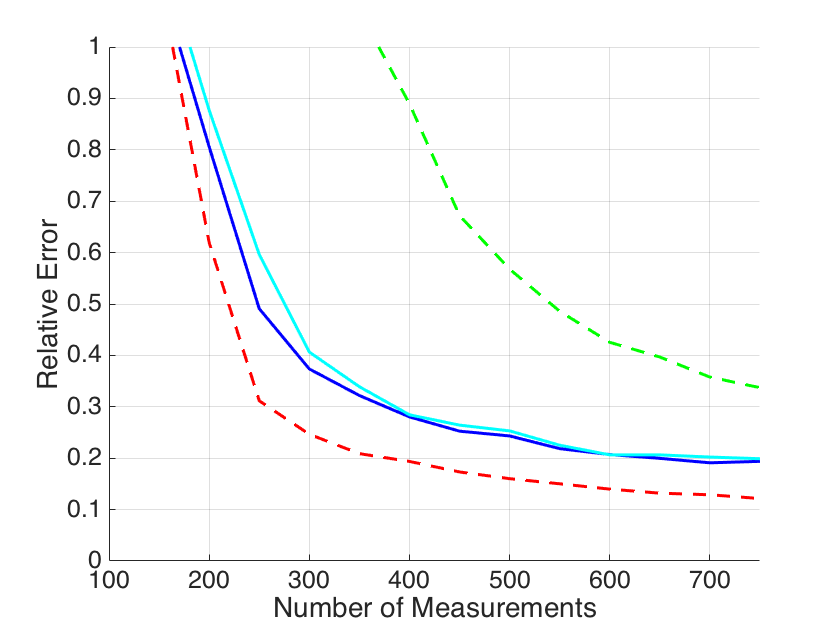}
 	\end{subfigure}
 	\quad
 	\begin{subfigure}[b]{0.45\textwidth}
 		\includegraphics[width=\textwidth]{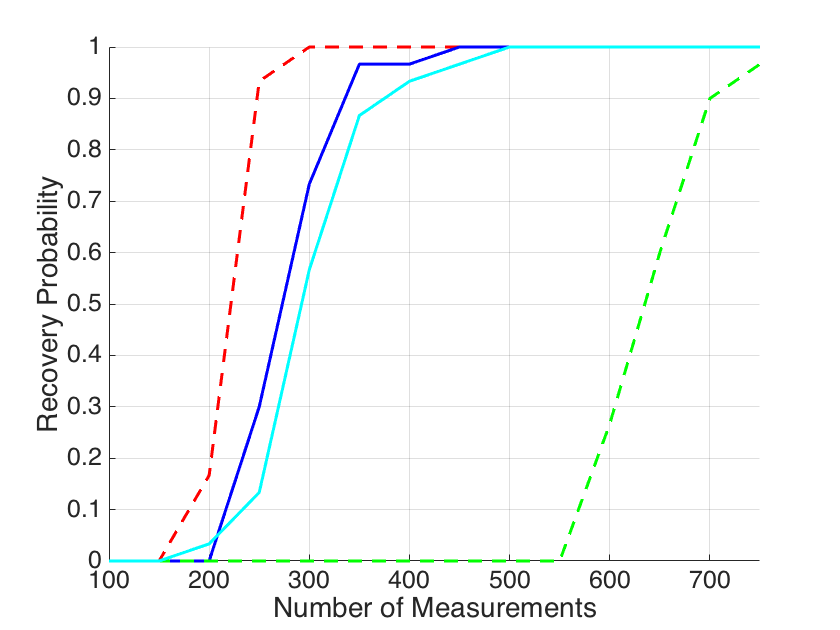}
 	\end{subfigure}
  \caption{ {\footnotesize Comparison of SPF and ATLAS with and without common support for $R=5$ {(see Section \ref{Numerics3})}. Depicted are average approximation error relative to $\Vert X \Vert_F$ and empirical recovery probabilities of SPF (dashed) and ATLAS (solid). Common Support: SPF (red) vs ATLAS (blue). Arbitrary Support: SPF (green) vs ATLAS (cyan).}}
\label{result}
\end{figure}

To further quantify this effect, we \purple{perform} the experiments \purple{reflected} in Figure \ref{result}. For varying number of measurements we compared average approximation error and recovery probability of SPF and ATLAS for $30$ randomly chosen $\hat{X} \in \R^{16 \times 100}$, $\Vert \hat{X} \Vert_F = 10$, with $\rank(\hat{X}) = 5$ and $10$-sparse right singular vectors which either share a common support or may have various support sets. The parameters are chosen as $\alpha = \beta = 0.5$. One can clearly see that SPF outperforms ATLAS even in the noisy case for common support sets of the singular vectors. This is not surprising as ATLAS makes no use of the additional information provided by shared support sets. If the singular vectors, however, do not share a common support set, ATLAS shows its strength in the noisy setting. SPF which needs pre-information on the row-/column-sparsity $\tilde{s}$ of $\hat{X}$ has to be initialized with $\tilde{s} = Rs$ as in the general case all support sets may differ. \\
\begin{figure}[!htb]
 	\centering
 	\captionsetup{width=0.9\linewidth}
 	\begin{subfigure}[b]{0.15\textwidth}
 		\includegraphics[width=\textwidth]{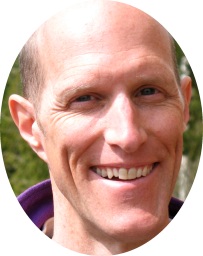}
 	\end{subfigure}
 	\quad
 	\begin{subfigure}[b]{0.15\textwidth}
 		\includegraphics[width=\textwidth]{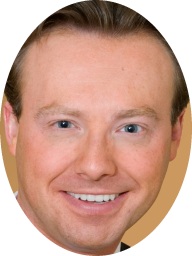}
 	\end{subfigure}
 	\quad
 	\begin{subfigure}[b]{0.15\textwidth}
 		\includegraphics[width=\textwidth]{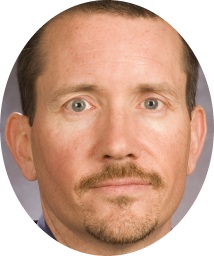}
 	\end{subfigure}
 	\quad
 	\begin{subfigure}[b]{0.15\textwidth}
 		\includegraphics[width=\textwidth]{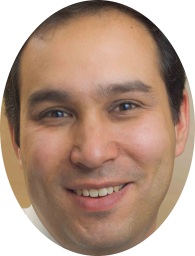}
 	\end{subfigure}
 	\quad
 	\begin{subfigure}[b]{0.15\textwidth}
 		\includegraphics[width=\textwidth]{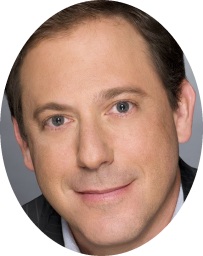}
 	\end{subfigure}
 	\quad
  \caption{ \blue{Five samples from original face data.} }
\label{fig:OriginalFaces}
\end{figure}
\blue{
Let us conclude the section with a comparison of SPF and ATLAS on real-life data. To this end, we choose ten faces from the "10k US Adult Faces Data Base" \cite{bainbridge2013intrinsic}, cf.\ Figure \ref{fig:OriginalFaces}, reduce their resolution to 64x44 pixels and their color range to gray-scale, and apply a $5$-layer wavelet transform with Haar-wavelets to obtain effectively sparse representations. The resulting ten $2836$-dimensional coefficient vectors -- not sparse but with an effective sparsity level of about $s = 66$ -- are then used to build a ground truth matrix $\hat{X} \in \R^{10\times 2836}$ which is re-scaled to unit Frobenius norm. Though not rank-deficient, $\hat{X}$ has effective rank $\boldsymbol{r}(\hat{X}) \approx 1.5$ where
\begin{align*}
    \boldsymbol{r}(\hat{X}) = \frac{\| \hat{X} \|_*}{\| \hat{X} \|} \in [ 1,\rank(\hat{X}) ]
\end{align*}
is a relaxed measure of low-rankness of a matrix similar to effective sparsity for vectors. In particular, $\hat{X}$ is well approximated by low-rank matrices. As in the above experiments $\mathcal{A}$ is a Gaussian operator and the noise level is set to $\| \eta \|_2 = 0.2 \| \hat{X} \|_F$. We choose the number of measurements as ten times the information theoretic limit $R(n_1+s)$ where $R$ is the rank parameter used for SPF and ATLAS in the experiments below. We initialize both algorithms with the leading singular vectors of $\mathcal{A}^*(y)$.\\
\begin{table}[!htb]
		\centering
		\blue{
		\begin{tabular}{ || m{3cm} || m{1.6cm}|m{1.6cm}|m{1.6cm} || } 
			\hline
			Setting & $R = 2$ & $R = 3$ & $R = 4$ \\
			\hline \hline
			SPF (B.A.) & \textbf{0.2239} & 0.2241 & \textbf{0.2098}    \\
			ATLAS (B.A.) & 0.2360 & \textbf{0.2195} & 0.2140 \\
			\hline
			SPF (D.P.) & 0.3043 & 0.2753 & 0.2532 \\
			ATLAS (D.P.) & \textbf{0.2382} & \textbf{0.2247} & \textbf{0.2210} \\
			\hline
			$\| \hat{X} - \hat{X}_R \|_F$ & 0.1310 & 0.1095 & 0.0932 \\
			\hline
		\end{tabular}
		\caption{\blue{Comparison of SPF and ATLAS for different choices of $R$ (the error produced by the best rank-$R$ term approximation $\hat{X}_R$ of $\hat{X}$ is given as a benchmark). Depicted is the relative approximation error measured in Frobenius norm. The parameters $s$ resp.\ $\alpha,\beta$ have been tuned on a discrete search grid either using \textbf{B}est \textbf{A}pproximation, i.e., minimizing $\| X_{\text{rec}} - \hat{X} \|_F$, or \textbf{D}iscrepancy \textbf{P}rinciple, i.e., minimizing $| \| \mathcal{A}(X_{\text{rec}}) - y \|_2 - \| \eta \|_2|$.}} \label{fig:Table1}
		}
\end{table}
\begin{table}[!htb]
		\centering
		\blue{
		\begin{tabular}{ || m{3cm} || m{1.6cm}|m{1.6cm}|m{1.6cm} || } 
			\hline
			\# perm.\ rows & $0$ & $1$ & $2$ \\
			\hline \hline
			SPF (B.A.) & \textbf{0.2239} & 0.2508 & 0.2427    \\
			ATLAS (B.A.) & 0.2360 & \textbf{0.2412} & \textbf{0.2347} \\
			\hline
			SPF (D.P.) & 0.3043 & 0.3140 & 0.3185 \\
			ATLAS (D.P.) & \textbf{0.2382} & \textbf{0.2369} & \textbf{0.2369} \\
			\hline
			$\| \hat{X} - \hat{X}_R \|_F$ & 0.1310 & 0.1367 & 0.1283 \\
			\hline
		\end{tabular}
		\caption{\blue{Comparison of SPF and ATLAS for $\hat{X}$ having different numbers of randomly permuted rows and $R = \lceil \boldsymbol{r}(\hat{X}) \rceil$ (the error produced by the best rank-$R$ term approximation $\hat{X}_R$ of $\hat{X}$ is given as a benchmark). Depicted is the relative approximation error measured in Frobenius norm. The parameters $s$ resp.\ $\alpha,\beta$ have been tuned on a discrete search grid either using \textbf{B}est \textbf{A}pproximation, i.e., minimizing $\| X_{\text{rec}} - \hat{X} \|_F$, or \textbf{D}iscrepancy \textbf{P}rinciple, i.e., minimizing $| \| \mathcal{A}(X_{\text{rec}}) - y \|_2 - \| \eta \|_2|$.}} \label{fig:Table2}
		}
\end{table}
In Table \ref{fig:Table1} we compare the full matrix reconstruction performance of SPF and ATLAS for different choices of the rank hyper-parameter $R$. As a benchmark, the error produced by best rank-$R$ term approximation is reported as well; this is the best achievable error under complete knowledge of $\hat{X}$ and without added noise. If the parameters are tuned under knowledge of the true solution $\hat{X}$ both algorithms perform similarly and allow reconstruction up to noise level while the reconstruction quality improves with increasing $R$. However, SPF performance worsens in the case $\hat{X}$ is not fed as information to the best approximation principle to tune hyper-parameters. When the hyper-parameters are tuned under exclusive knowledge of $\| \eta \|_2$, the reconstruction error produced by SPF is significantly larger than the one of ATLAS. \\
Since the wavelet transform creates a joint row support structure (similar positions of dominant entries) and one of the benefits of ATLAS is not to rely on joint supports, we repeat the experiment but randomly permute the entries of one resp.\ two rows of $\hat{X}$ to create a ground truth of higher effective rank and less joint support structure (in this case we set $R = \lceil \boldsymbol{r}(\hat{X}) \rceil$). Table \ref{fig:Table2} shows that SPF's reconstruction quality suffers more from this loss of structure. In particular, when using the discrepancy principle the performance gap becomes wider. \\
The most important difference between SPF and ATLAS can be observed when comparing the reconstructed images obtained by reverting the wavelet transform on the rows of $\hat{X}$. In this experiment we set the noise level to zero, increase the oversampling factor from ten to twenty, and set $R = 5$. Figures \ref{fig:Face1} and \ref{fig:Face5} reveal that, although SPF achieves a similar Frobenius error in reconstructing $\hat{X}$, it oversimplifies the images encoded in the rows of $\hat{X}$ and produces large pixel areas of uniform gray-level. Moreover, Figure \ref{fig:EigenFaces} (especially comparing the second eigenfaces) proves that the two algorithms search for qualitatively different decompositions of the ground truth. SPF stays closer to the original SVD while ATLAS has more freedom in decomposing $\hat{X}$. This can be seen as well when comparing the Gramians of the matrices $V_\text{SPF}, V_\text{ATLAS} \in \R^{2836\times 5}$ containing the right components reconstructed by SPF and ATLAS. They show that ATLAS is not restricted to orthogonal decompositions (the Gramian of SPF is not perfectly diagonal, since the last orthogonormalization is performed before the last application of Hard Thresholding Pursuit):
\begin{align*}
    V_\text{SPF}^T V_\text{SPF} &=
    \begin{pmatrix}
    0.9720  &  0.0003  &  0.0002 &   0.0000  &  0.0006 \\
    0.0003  &  0.0093  &  0.0000 &   0.0000  &  0.0001 \\
    0.0002  &  0.0000  &  0.0054 &   0.0000  &  0.0000 \\
    0.0000  &  0.0000  &  0.0000 &   0.0031  &  0.0001 \\
    0.0006  &  0.0001  &  0.0000 &   0.0001  &  0.0027 
    \end{pmatrix}
    \\ \vspace{3cm}
    V_\text{ATLAS}^T V_\text{ATLAS} &=
    \begin{pmatrix}
    1.1499  & -0.0038  &  0.0016  & -0.0135  &  0.0014 \\
   -0.0038  &  0.0311  &  0.0038  &  0.0044  &  0.0018 \\
    0.0016  &  0.0038  &  0.0232  &  0.0002  &  0.0008 \\
   -0.0135  &  0.0044  &  0.0002  &  0.0131  &  0.0010 \\
    0.0014  &  0.0018  &  0.0008  &  0.0010  &  0.0132 
    \end{pmatrix}
\end{align*}
\begin{figure}[!htb]
 	\centering
 	\captionsetup{width=0.9\linewidth}
 	\begin{subfigure}[b]{0.2\textwidth}
 		\includegraphics[width=\textwidth]{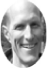}
 		\subcaption{Ground truth}
 	\end{subfigure}
 	\qquad
 	\begin{subfigure}[b]{0.2\textwidth}
 		\includegraphics[width=\textwidth]{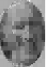}
 		\subcaption{SPF}
 	\end{subfigure}
 	\qquad
 	\begin{subfigure}[b]{0.2\textwidth}
 		\includegraphics[width=\textwidth]{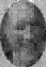}
 		\subcaption{ATLAS}
 	\end{subfigure}
 	\qquad
 	\begin{subfigure}[b]{0.2\textwidth}
 		\includegraphics[width=\textwidth]{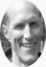}
 		\subcaption{Rank-$R$ SVD}
 	\end{subfigure}
  \caption{ \blue{Comparison of the first face and its reconstructions.} }
\label{fig:Face1}
\end{figure}
\begin{figure}[!htb]
 	\centering
 	\captionsetup{width=0.9\linewidth}
 	\begin{subfigure}[b]{0.2\textwidth}
 		\includegraphics[width=\textwidth]{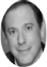}
 		\subcaption{Ground truth}
 	\end{subfigure}
 	\qquad
 	\begin{subfigure}[b]{0.2\textwidth}
 		\includegraphics[width=\textwidth]{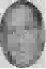}
 		\subcaption{SPF}
 	\end{subfigure}
 	\qquad
 	\begin{subfigure}[b]{0.2\textwidth}
 		\includegraphics[width=\textwidth]{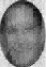}
 		\subcaption{ATLAS}
 	\end{subfigure}
 	\qquad
 	\begin{subfigure}[b]{0.2\textwidth}
 		\includegraphics[width=\textwidth]{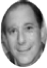}
 		\subcaption{Rank-$R$ SVD}
 	\end{subfigure}
  \caption{ \blue{Comparison of the fifth face and its reconstructions.} }
\label{fig:Face5}
\end{figure}
\begin{figure}[h!]
 	\centering
 	\captionsetup{width=0.9\linewidth}
 	\begin{subfigure}[b]{0.25\textwidth}
 	    \centering
 		\includegraphics[width=0.5\textwidth, height=2.525cm]{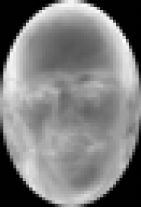}
 		\subcaption{Ground truth}
 	\end{subfigure}
 	\quad
 	\begin{subfigure}[b]{0.25\textwidth}
 	    \centering
 		\includegraphics[width=0.5\textwidth, height=2.525cm]{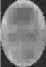}
 		\subcaption{SPF}
 	\end{subfigure}
 	\quad
 	\begin{subfigure}[b]{0.25\textwidth}
 	    \centering
 		\includegraphics[width=0.5\textwidth, height=2.525cm]{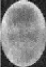}
 		\subcaption{ATLAS}
 	\end{subfigure}
 	\\
 	\begin{subfigure}[b]{0.25\textwidth}
 	    \centering
 		\includegraphics[width=0.5\textwidth, height=2.525cm]{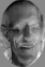}
 		\subcaption{Ground truth}
 	\end{subfigure}
 	\quad
 	\begin{subfigure}[b]{0.25\textwidth}
 	    \centering
 		\includegraphics[width=0.5\textwidth, height=2.525cm]{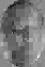}
 		\subcaption{SPF}
 	\end{subfigure}
 	\quad
 	\begin{subfigure}[b]{0.25\textwidth}
 	    \centering
 		\includegraphics[width=0.5\textwidth, height=2.525cm]{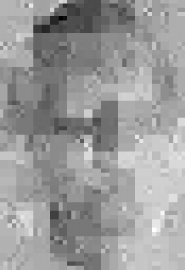}
 		\subcaption{ATLAS}
 	\end{subfigure}
 	\\
 	\begin{subfigure}[b]{0.25\textwidth}
 	    \centering
 		\includegraphics[width=0.5\textwidth, height=2.525cm]{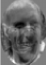}
 		\subcaption{Ground truth}
 	\end{subfigure}
 	\quad
 	\begin{subfigure}[b]{0.25\textwidth}
 	    \centering
 		\includegraphics[width=0.5\textwidth, height=2.525cm]{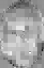}
 		\subcaption{SPF}
 	\end{subfigure}
 	\quad
 	\begin{subfigure}[b]{0.25\textwidth}
 	    \centering
 		\includegraphics[width=0.5\textwidth, height=2.525cm]{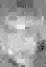}
 		\subcaption{ATLAS}
 	\end{subfigure}
  \caption{ \blue{Comparison of the three leading eigenfaces (back transformed leading right singular vectors) of $\hat{X}$ and the corresponding leading rank-$1$ components reconstructed by SPF and ATLAS.} }
\label{fig:EigenFaces}
\end{figure}
%
%
}

\subsection{{Initialization}} \label{Numerics4}

\begin{figure}[h!]
	\centering
	\captionsetup{width=.8\linewidth}
	\begin{subfigure}[b]{\textwidth}
		\includegraphics[width=\textwidth]{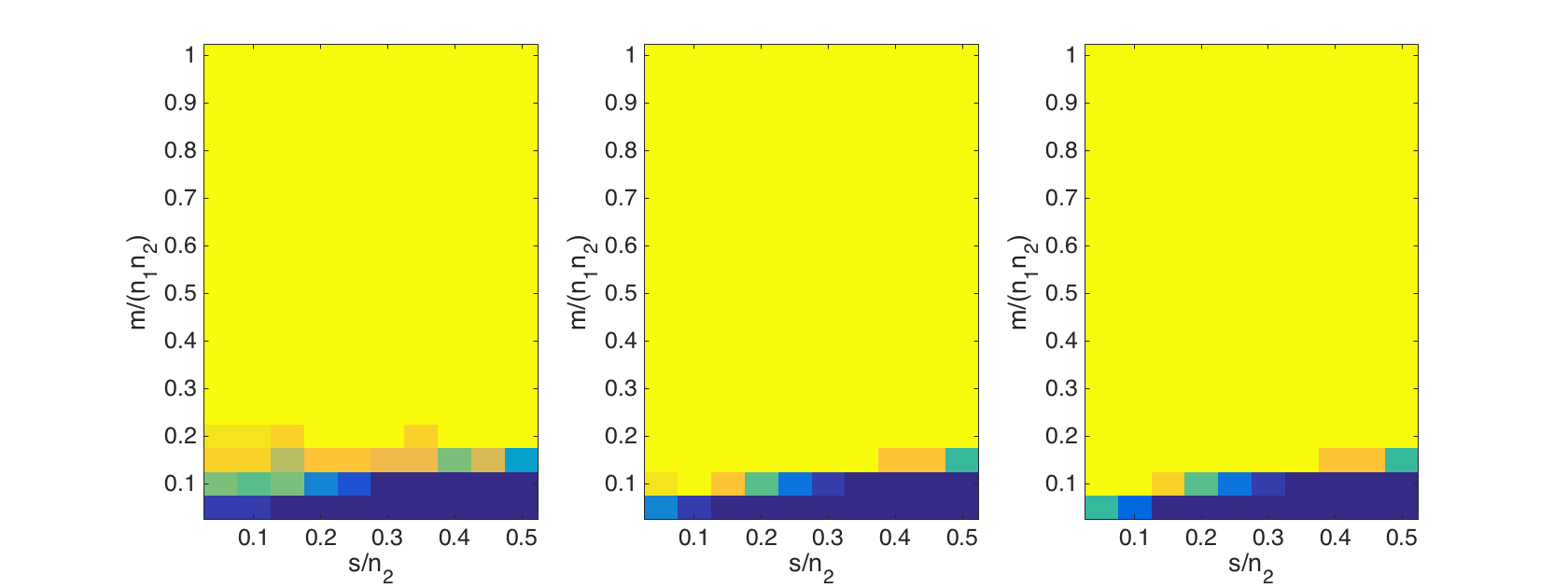}
		\caption{$R = 1$}
	\end{subfigure}
	\quad 
	\begin{subfigure}[b]{\textwidth}
		\includegraphics[width=\textwidth]{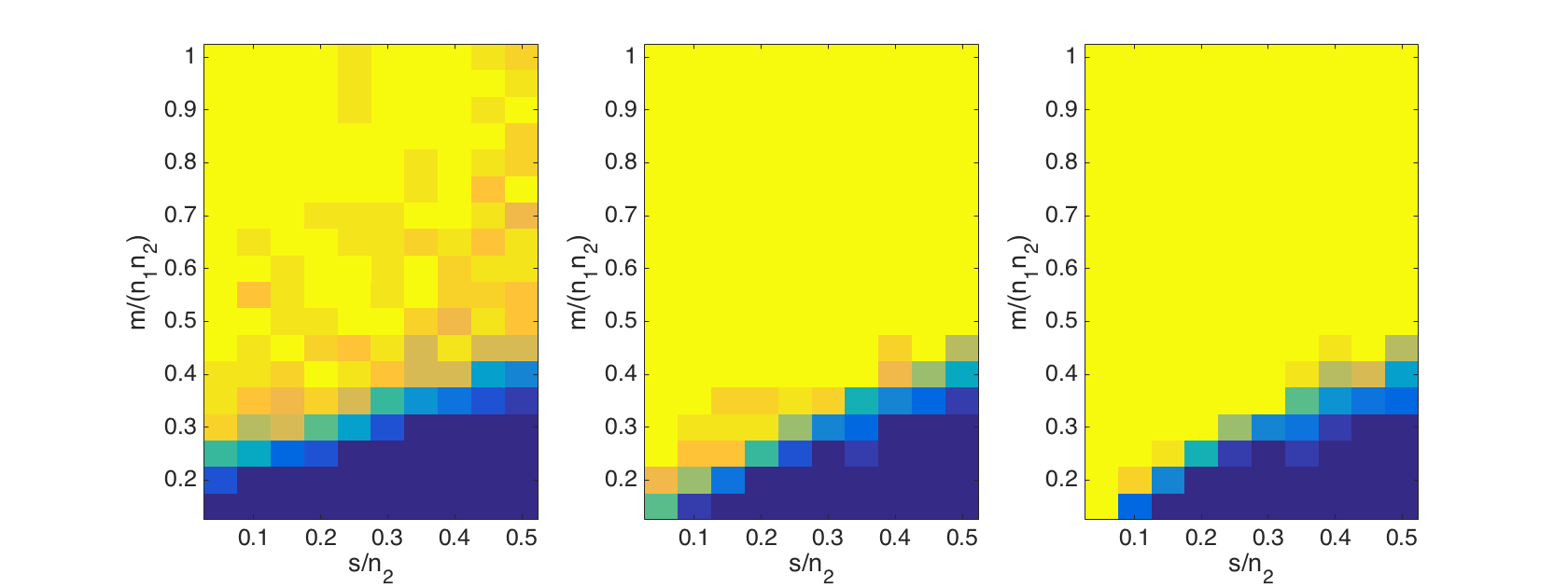}
		\caption{$R = 3$}
	\end{subfigure}
	\caption{Comparison of different initializations for ATLAS for \eqref{A} with noise $\eta \neq 0$ on the measurments (see Section \ref{Numerics4}), namely,  initialization with a strongly perturbed approximation $X_0 \approx \hat X$ (left),  initialization by the leading singular vectors of $\A^*(y)$ (middle), and initialization with a mildly perturbed approximation $X_0 \approx \hat X$  (right). Empirical recovery probability is depicted by color from zero (blue) to one (yellow).}\label{Fig:InitPlot}
\end{figure}

We \red{also perform} a simple test on {the influence of the initialization}. The plots in Figure \ref{Fig:InitPlot} compared for $s/n_2 \in [0,0.5]$ and {$m/(n_1n_2) \in [0,1]$} the number of successful recoveries of $20$ randomly drawn $\hat{X} \in \R^{8 \times 128}$, $\Vert \hat{X} \Vert_F = 10$, with $\rank(\hat{X}) \in \{1,3\}$ and $s$-sparse right singular vectors from $m$ measurements. The noise level was set to $0.3\Vert \hat{X} \Vert_F$ and recovery was counted successful if $\Vert \hat{X} - X_\text{appr} \Vert_F/\Vert \hat{X} \Vert_F \le 0.4$. We compare the initializations by the leading singular vectors of $\A^*(y)$ and by the leading singular vectors of $X + Z$ where $Z$ was drawn at random, and scaled to $\Vert Z \Vert_F = 100$ (strong perturbation) resp. $\Vert Z \Vert_F = 0.2$ (mild perturbation). 
\\
{For $\rank(\hat{X})=1$ we note remarkably that the convergence radius of ATLAS is seemingly very large (yet not global), as the phase transition diagrams in Figure \ref{Fig:InitPlot} do not show significant variations from choosing as initialization the leading singular vectors of $\A^*(y)$ and those of small random perturbation. Instead for $\rank(\hat{X})=3$, initialization plays a more important} role in performance and the \purple{initialization by leading singular vectors of $\A^*(y)$} does not yield optimal performance.


\section{Discussion and Open Questions} \label{Conclusion}

Motivated by challenging examples from recommendation systems and blind demixing in signal processing, in this paper \red{we propose a multi-penalty approach to recover low-rank matrices with sparsity structure from incomplete and inaccurate linear measurements.}
The considered problem stands at the intersection of the compressed sensing and sparse PCA framework, though significantly extending their settings.
Our analysis results in general bounds on the performance of the proposed algorithm, ATLAS, and in a necessary number of subgaussian measurements to approximate effectively sparse and low-rank matrices. These theoretical results are confirmed in numerical experiments.
ATLAS is especially of use and effective in the most realistic setting of ineliminable noise and, hence, it complements the state-of-the-art algorithm SPF of Lee et.\ al.\ in \cite{lee2013near}, which works well for low level of noise or exact measurements only. \red{While the theoretical guarantees for SPF are sharper, ATLAS tackles the recovery of a significantly larger class of matrices, i.e., matrices with non-orthogonal rank-$1$ decompositions and effectively sparse components, which are of interest when turning to more general tasks as, for instance, in machine learning.} Nevertheless, the analysis of ATLAS is remarkably simple and it is easily prone  to several extensions/generalizations. We mention a few of them as follows.
\\
 First, replacing the $\ell_2$- and $\ell_1$-norms by $\ell_p$- and $\ell_q$-(quasi)-norms, for $p \ge 2$ and $0 < q < 2$, yields to the functional
\begin{equation} \label{Jabpq}
J_{\alpha,\beta}^{R,p,q}(u^1, \dots, u^R, v^1, \dots, v^R) := \left\Vert y - \mathcal A\left(\sum_{r=1}^R u^r (v^r)^T\right) 
\right\Vert_2^2 + \alpha \sum_{r=1}^R \Vert u^r \Vert_p^p + \beta \sum_{r=1}^R \Vert v^r  
\Vert_q^q.
\end{equation}
 and, in turn, the algorithm $({\text{A-T-LAS}_{p,q}})$
\begin{equation}
	 \label{ATLASpq}
\begin{cases}
u_{k+1}^1 &=  \argmin_u \left\Vert \left( y - \A \left( \sum_{r=2}^R u_k^r {v_k^r}^T 
\right) \right) - \A(u {v_k^1}^T) \right\Vert^2_2 + \alpha \Vert u \Vert^p_p
+ \frac{1}{2\lambda_k^1} \Vert u - u_k^1 \Vert^2_2, \\
v_{k+1}^1 &=  \argmin_v \left\Vert \left( y - \A \left( \sum_{r=2}^R u_k^r {v_k^r}^T 
\right)	\right) - \A(u_{k+1}^1 v^T) \right\Vert^2 + \beta \Vert v \Vert_q^q
+ \frac{1}{2\mu_k^1} \Vert v - v_k^1 \Vert^2_2, \\
&\vdots \\
u_{k+1}^R &=  \argmin_u \left\Vert \left( y - \A \left( \sum_{r=1}^{R-1} u_{k+1}^r 
{v_{k+1}^r}^T \right) \right) - \A(u {v_k^R}^T) \right\Vert^2_2 + \alpha \Vert u \Vert^p_p
+ \frac{1}{2\lambda_k^R} \Vert u - u_k^R \Vert^2_2, \\
v_{k+1}^R &=  \argmin_v \left\Vert \left( y - \A \left( \sum_{r=1}^{R-1} u_{k+1}^r 
{v_{k+1}^r}^T \right) \right) - \A(u_{k+1}^R v^T) \right\Vert^2 + \beta \Vert v \Vert_q^q
+ \frac{1}{2\mu_k^R} \Vert v - v_k^R \Vert^2_2, \\
\end{cases}
\end{equation}
As for $q < 1$ even the single component minimizations become non-convex, this setting needs special care. One would need non-standard iterative thresholding methods, which have been developed and studied, e.g., in \cite{naumova2014minimization}. As $q$-quasi-norms, for $q < 1$, have proved particularly effective in enforcing sparsity, this additional technical difficulties are worth to overcome. 
\\
Second, in recommendation systems, one usually imposes additionally non negativity constraints on the obtained matrices. We could easily implement them in ATLAS by asymmetric $\ell_1$-regularization. Define for $z \in \R^n$ and $\theta > 0$
\begin{align*}
	\| z \|_{1,\theta}^+ \coloneqq \sum_{i=1}^n |z_i|_{\theta}^+, \quad |x|_\theta^+ \coloneqq \begin{cases}
	x & x \ge 0 \\
	\theta |x| & \text{else.}
	\end{cases}
\end{align*}
For $\theta$ becoming large, the regularization by $\| \cdot \|_{1,\theta}^+$ promotes sparsity and non-negativity.  Replacing the $\ell_1$-norm in ATLAS by $\| \cdot \|_{1,\theta}^+$ would result in the remarkably simple modification of ISTA (Algorithm \ref{A2}), where the soft-thresholding operator $\St_\beta$ is substituted in line 2 with
\begin{align*}
	\St_{\beta,\theta} = \begin{pmatrix}
	S_{\beta, \theta} (z_1) \\ \vdots \\ S_{\beta, \theta} (z_{n_2}) 
	\end{pmatrix}, \tab
	\text{ where } S_{\beta, \theta} (z_i) = \begin{cases}
	z_i - \frac{\beta}{2} & z_i > \frac{\beta}{2} \\
	0 & -\theta \frac{\beta}{2} \le z_i \le \frac{\beta}{2} \\
	z_i + \theta\frac{\beta}{2} & z_i < -\theta \frac{\beta}{2}
	\end{cases}.
\end{align*}
Note that in the limit case $\theta \rightarrow \infty$ the operator $\St_{\beta,\theta}$ is a shifted ReLU function. Choosing $\theta$ sufficiently large or considering the limit $\theta \rightarrow \infty$ would lead to non-negative sparse PCA \cite{zass2007nonnegative} from incomplete and inaccurate measurements with further applications in economics \cite{jagannathan2003risk}, biology \cite{badea2005sparse}, and computer vision \cite{lee1999learning}.\\
Third, \red{as a byproduct of our generalizations}, we introduce, in our view, the right class of matrices $K_{s_1,s_2}^{R,\Gamma}$ and corresponding RIP, Definition \ref{apprRIPDef}, which might allow to generalize SPF to matrices having non-orthogonal, effectively sparse decompositions. In fact, the assumption of SPF of model matrices with sparse SVD is quite restrictive. For $q \to 0,$ the ATLAS algorithm  can be  seen as a  generalization of SPF and be realized by iterative hard-thresholding.
\\	
The current results demand a careful choice of parameters at noise level. This drawback of multi-penalty regularization is well-known and could be attacked by implementing LASSO-path. LASSO-path has been recently extended to the multi-penalty setting in case of superposition of the signals \cite{grklna17}, where the authors provided an efficient procedure for the construction of regions containing structurally similar solutions. In addition, $J_\abg^R$ depends by construction heavily on pre-knowledge of the rank $R$. One might ask how to get good estimates for $R$ in case the rank is unknown.\\
As mentioned above,  initialization  is crucial for  good performances of the algorithm. It is currently unclear how a good initialization can be obtained to guarantee convergence of the whole procedure to global minimizers. This question is closely connected to the fundamental problem in non-convex optimization how to initialize gradient-descent methods. In fact, alternating minimization is somewhat related to gradient-descent. While in gradient-descent one determines an optimal descent direction and then approximates the optimal step size, alternating minimization strongly restricts the directions in space in order to calculate optimal step sizes. Lee et.\ al.\ proposed an initialization, which worked in their setting if one assumes a strong decay of the singular values. 
Possibly one could prove this initialization to be sufficiently good in our setting as well, also in the light of recently \purple{improved analysis} \cite{gekrst17}.


\section{Appendix}

\subsection{Proofs of \Cref{Proofs}}

\paragraph{} \purple{Proofs of} two technical \purple{results} used in \ref{Proofs} are provided here. The first \purple{result estimates} possible coverings of {$S_{s_1,s_2}^{R,\Gamma}$ defined in \eqref{S}} \purple{as stated in \Cref{CoveringNumber}}, while the second \purple{result} contains two integral estimates \purple{used in the proof of \Cref{GaussianRIP}}.

\begin{Proof}[of \Cref{CoveringNumber}]
	Recall, each $Z \in S_{s_1,s_2}^{R,\Gamma}$ can be represented as $Z = U\Sigma V^T$ with $U = (u^1,...,u^R)$, $V = (v^1,...,v^R)$ where all unit norm columns $u^r \in \R^{n_1}$ are $s_1$-sparse, all unit norm columns $v^r \in \R^{n_2}$ are $s_2$-sparse, and $\Vert \Sigma \Vert_F \le \Gamma$. Let us first consider the larger set $S = \{ Z = U\Sigma V^T \colon U \in Q_{n_1,s_1}^R, \Sigma \in D_\Gamma, \text{ and } V \in Q_{n_2,s_2}^R \}$ where $D_\Gamma$ is the set of $R\times R$ diagonal matrices with Frobenius norm less or equal $\Gamma$ and $Q_{n,s}^R = \{ W \in \R^{n\times R} \colon \Vert W \Vert_F \le \sqrt{R} \text{ and all columns } w^r \text{ are } s\text{-sparse} \}$. Then, we know that $S_{s_1,s_2}^{R,\Gamma} \subset S$. We construct an $(\eps/2)$-net $\tilde{S}$ of $S$ by covering the sets of permissible $U$, $\Sigma$, and $V$ and conclude the proof by applying the well-known relation $N(K,\Vert \cdot \Vert, \eps) \le N(K',\Vert \cdot \Vert, \eps/2)$ which holds whenever $K \subset K'$.
	\paragraph{} First note that if $B$ is a unit ball in $D$ dimensions (with respect to some norm $\Vert \cdot \Vert_B$) there exists an $\eps$-net $\tilde{B}$ (i.e., for all $b \in B$ there is some $\tilde{b} \in \tilde{B}$ with $\Vert b - \tilde{b} \Vert_B \le \eps$) with $\tilde{B} \subset B$ and $|\tilde{B}| \le (3/\eps)^D$. See for example \cite[Begin Section 3]{candes2011tight}. Moreover, note that $N(K,\Vert \cdot \Vert,\eps) = N(cK,\Vert \cdot \Vert,c\eps)$ for any set $K$ and $c > 0$. Hence, for any scaled unit ball $cB$ there exists an $\eps$-net $\tilde{B} \subset cB$ and $|\tilde{B}| \le (3c/\eps)^D$.
	\paragraph{} Let $\tilde{D}_\Gamma$ be an $(\eps/(6R))$-net of $D_\Gamma$ which is of size $|\tilde{D}_\Gamma| \le (18\Gamma R/\eps)^R$. For $W \in \R^{n\times R}$ denote by $\supp(W) = \{ \supp(w^1),...,\supp(w^R) \}$ and by $\supp(W) \Subset \supp(W')$ that $\supp(w^r) \subset \supp((w')^r)$, for all $r \in [R]$. Define the set of all possible supports of maximal size
	\begin{align*}
		T_{n,s}^R = \{ \supp(W) \colon W \in \R^{n \times R} \text{ and all columns } w^r \text{ have exactly } s \text{ non-zero entries} \}.
	\end{align*} 
	For any fixed $\theta \in T_{n,s}^R$ the set $\{ W \in Q_{n,s}^R \colon \supp(W) \Subset \theta \}$ is an $\R^{s\times R}$ Frobenius ball of radius $\sqrt{R}$ embedded into $\R^{n\times R}$ and $Q_{n,s}^R = \bigcup_{\theta \in T_{n,s}^R} \{ W \in Q_{n,s}^R \colon \supp(W) \Subset \theta \}$. Hence, there is an $(\eps/(6\Gamma \sqrt{R}))$-net $\tilde{Q}_{n,s}^R$ of $Q_{n,s}^R$ with
	\begin{align*}
		| \tilde{Q}_{n,s}^R | \le | T_{n,s}^R | \left( \frac{18\Gamma R}{\eps} \right)^{Rs} \le \binom{n}{s}^R \left( \frac{18\Gamma R}{\eps} \right)^{Rs} \le  \left( \frac{en}{s} \right)^{Rs}  \left( \frac{18\Gamma R}{\eps} \right)^{Rs}
	\end{align*}
	We define now $\tilde{S} = \{ \tilde{Z} = \tilde{U} \tilde{\Sigma} \tilde{V}^T \colon \tilde{U} \in \tilde{Q}_{n_1,s_1}^R, \tilde{\Sigma} \in \tilde{D}_\Gamma, \text{ and } \tilde{V} \in \tilde{Q}_{n_2,s_2}^R \}$. It is clear that
	\begin{align*}
		|\tilde{S}| \le | \tilde{Q}_{n_1,s_1}^R | \cdot | \tilde{D}_\Gamma | \cdot | \tilde{Q}_{n_2,s_2}^R | \le \left( \frac{18\Gamma R}{\eps} \right)^{R(s_1+s_2+1)} \left( \frac{e n_1}{s_1} \right)^{Rs_1} \left( \frac{e n_2}{s_2} \right)^{Rs_2}.
	\end{align*}
	Let us conclude by showing $\tilde{S}$ is indeed an $(\eps/2)$-net for $S$. Given any $Z = U\Sigma V^T \in S$, there exists ${\tilde{Z} = \tilde{U} \tilde{\Sigma} \tilde{V}^T} \in \tilde{S}$ with $\Vert U - {\tilde{U}} \Vert_F \le \eps/(6\Gamma \sqrt{R})$, $\Vert \Sigma - \tilde{\Sigma} \Vert_F \le \eps/(6R)$, and $\Vert V - {\tilde{V}} \Vert_F \le \eps/(6\Gamma \sqrt{R})$. We can estimate
	\begin{align*}
	\Vert Z - {\tilde{Z}} \Vert_F &\le \Vert (U-{\tilde{U}})\Sigma V^T \Vert_F + \Vert \tilde{U} (\Sigma - \tilde{\Sigma}) V^T \Vert_F + \Vert {\tilde{U}} \tilde{\Sigma} (V - {\tilde{V}})^T \Vert_F \\
	&\le \frac{\eps}{6\Gamma \sqrt{R}} \Gamma \sqrt{R} + \sqrt{R} \frac{\eps}{6R} \sqrt{R} + \sqrt{R} \Gamma \frac{\eps}{6\Gamma \sqrt{R}} \\
	&\le \frac{\eps}{2}
	\end{align*}
	where we used triangle inequality in the first line and $\Vert AB \Vert_F \le \Vert A \Vert_F \Vert B \Vert_F$ in the second.
\end{Proof}

\begin{lemma} \label{GammaBounds}
	If $\Gamma \ge 1$, we have for the sets $S_{s_1,s_2}^{R,\Gamma}$ and $K_{s_1,s_2}^{R,\Gamma}$ defined in \eqref{S} and \eqref{K} that
	\begin{align*}
		\int_{0}^{\frac{\Gamma \sqrt{R}}{\sqrt{m}}} \sqrt{ \log N \left( S_{s_1,s_2}^{R,\Gamma} , \Vert \cdot \Vert_F ,\sqrt{m} \varepsilon \right) } d\varepsilon 
		&\le \sqrt{\frac{C_S \Gamma^2 R^2 (s_1+s_2) \log \left( \max \left\{ n_1,n_2 \right\} \right) }{m}} \\
		\int_{0}^{\frac{\Gamma \sqrt{R}}{\sqrt{m}}} \sqrt{\log N(K_{s_1,s_2}^{R,\Gamma},\Vert \cdot \Vert_F,\sqrt{m}\eps)} \;d\eps &\le \sqrt{\frac{C_K \Gamma^2 R^2 (s_1 + s_2) \log^3(\max\{n_1,n_2\})}{m}}
	\end{align*}
	where $C_S,C_K > 0$ are constants.
\end{lemma}
\begin{Proof}
	For the first estimate apply \Cref{CoveringNumber} to obtain
	\begin{align*}
		\int_{0}^{\frac{\Gamma \sqrt{R}}{\sqrt{m}}} &\sqrt{ \log N \left( S_{s_1,s_2}^{R,\Gamma} , \Vert \cdot \Vert_F ,\sqrt{m} \varepsilon \right) } d\varepsilon
		\le \sqrt{ \int_{0}^{\frac{\Gamma \sqrt{R}}{\sqrt{m}}} 1\; d\eps \int_{0}^{\frac{\Gamma \sqrt{R}}{\sqrt{m}}}\log N \left( S_{s_1,s_2}^{R,\Gamma} , \Vert \cdot \Vert_F ,\sqrt{m} \varepsilon \right) d\varepsilon } \\
		&\le \sqrt{ \frac{\Gamma^2 R^2 (s_1+s_2+1) {\left( 1 +\log\left( 18\sqrt{R} \right) \right)} + \Gamma^2 R^2 s_1 \log \left( \frac{en_1}{s_1} \right) + \Gamma^2 R^2 s_2 \log \left( \frac{en_2}{s_2} \right)}{m} } \\
		&\le \sqrt{\frac{C_S \Gamma^2 R^2 (s_1+s_2) \log \left( \max \left\{ n_1,n_2 \right\}  \right)}{m}},
	\end{align*}
	where we used Cauchy-Schwarz inequality in the first step and the fact that $\sqrt{R} \le \max\{ n_1,n_2 \}$ in the last inequality. $C_S > 0$ is an appropriate constant.
	\paragraph{} To obtain the second estimate let us first assume $s_1/n_1 \le s_2/n_2$. We apply \Cref{CoveringNumber2} and find
	\begin{align*}
	&\int_{0}^{\frac{\Gamma \sqrt{R}}{\sqrt{m}}} \sqrt{\log N(K_{s_1,s_2}^{R,\Gamma},\Vert \cdot \Vert_F,\sqrt{m}\eps)} \;d\eps \\
	&\le \int_{0}^{12\Gamma \sqrt{\frac{Rs_1}{mn_1}}} \sqrt{R(n_1+n_2+1) \log\left( \frac{36\Gamma R}{\sqrt{m}\eps} \right)} \de + \int_{12\Gamma \sqrt{\frac{Rs_1}{mn_1}}}^{12\Gamma \sqrt{\frac{Rs_2}{mn_2}}} \sqrt{\frac{144\Gamma^2R^2s_1}{m\eps^2} \log\left( \frac{9\sqrt{m}\eps n_1}{6\Gamma \sqrt{R} s_1} \right)} \de \\
	&+ \int_{12\Gamma \sqrt{\frac{Rs_1}{mn_1}}}^{12\Gamma \sqrt{\frac{Rs_2}{mn_2}}} \sqrt{R(n_2+1) \log\left( \frac{36\Gamma R}{\sqrt{m}\eps} \right)} \de + \int_{12\Gamma \sqrt{\frac{Rs_2}{mn_2}}}^{\frac{\Gamma \sqrt{R}}{\sqrt{m}}} \sqrt{\frac{144\Gamma^2 R^2(s_1 + s_2)}{m\eps^2} \log\left( \frac{9\sqrt{m}\eps n_1}{6\Gamma \sqrt{R} s_1} \right)} \de \\
	&+ \int_{12\Gamma \sqrt{\frac{Rs_2}{mn_2}}}^{\frac{\Gamma \sqrt{R}}{\sqrt{m}}} \sqrt{R \log\left( \frac{18\Gamma R}{\sqrt{m} \eps} \right)} \\
	&= I_1 + I_2 + I_3 + I_4 + I_5.
	\end{align*}
	We now estimate the five integrals. We use the short notation $a_i = 12\Gamma \sqrt{\frac{Rs_i}{mn_i}}$ for $i = 1,2$ and $b = \frac{\Gamma \sqrt{R}}{\sqrt{m}}$. The first integral can be bounded by
	\begin{align*}
	I_1 &\le \left( \int_0^{a_1} 1 \de \int_0^{a_1} R(n_1+n_2+1) \log\left( \frac{36\Gamma R}{\sqrt{m}\eps} \right) \de \right)^{\frac{1}{2}} \\
	&\le \left( a_1 R (n_1+n_2+1) \left[ \eps \left(1+\log \left( \frac{36\Gamma R}{\sqrt{m}\eps} \right)\right) \right]_{\eps = 0}^{a_1} \right)^{\frac{1}{2}} \\
	&= \left( \frac{144 \Gamma^2 R^2 s_1 (n_1 + n_2+1)}{mn_1} \left(1+\log \left( 3 \sqrt{\frac{Rn_1}{s_1}} \right)\right) \right)^\frac{1}{2} \le \begin{cases}
	\left( \frac{432 \Gamma^2 R^2 s_1}{m} \left(1+\log \left( 3 \sqrt{\frac{Rn_1}{s_1}} \right)\right) \right)^\frac{1}{2} & n_1 \ge n_2 \\
	\left( \frac{432 \Gamma^2 R^2 s_2}{m} \left(1+\log \left( 3 \sqrt{\frac{Rn_1}{s_1}} \right)\right) \right)^\frac{1}{2} & \text{else}
	\end{cases}
	\intertext{where we used in the last step the assumption $s_1/n_1 \le s_2/n_2$. As can be seen later, the case distinction is irrelevant in the final estimate. Let us now turn to the second integral.}
	I_2 &= \sqrt{\frac{144\Gamma^2 R^2 s_1}{m}} \int_{a_1}^{a_2} \frac{1}{\eps} \sqrt{\log\left( \frac{9\sqrt{m}\eps n_1}{6\Gamma \sqrt{R} s_1} \right)} \de = \sqrt{\frac{144\Gamma^2 R^2 s_1}{m}} \left[ \frac{2}{3} \log^\frac{3}{2} \left( \frac{9\sqrt{m}\eps n_1}{6\Gamma \sqrt{R} s_1} \right) \right]_{\eps = a_1}^{a_2} \\
	&{= \left( \frac{64\Gamma^2 R^2 s_1}{m} \right)^{\frac{1}{2}} \left( \log^\frac{3}{2} \left( \frac{18\sqrt{s_2}n_1}{\sqrt{n_2}s_1} \right) - \log^\frac{3}{2} \left( \frac{18\sqrt{n_1}}{\sqrt{s_1}} \right) \right) \le \left( \frac{64\Gamma^2 R^2 s_1}{m} \log^3 (18 n_1) \right)^{\frac{1}{2}}}
	\intertext{The third integral is similar to the first. Again the case distinction does not play a major role in the end.}
	I_3 &\le \left( (a_2-a_1) R (n_2+1) \left[ \eps \left(1+\log \left( \frac{36\Gamma R}{\sqrt{m}\eps} \right)\right) \right]_{\eps = a_1}^{a_2} \right)^{\frac{1}{2}} \\
	&= \left( (a_2-a_1) R (n_2+1) \left[ a_2 \left(1+\log \left( 3 \sqrt{\frac{Rn_2}{s_2}} \right) \right) - a_1 \left(1+\log \left( 3 \sqrt{\frac{Rn_1}{s_1}} \right)\right) \right] \right)^{\frac{1}{2}} \\
	&\le \left( (a_2-a_1)^2 R (n_2+1) \left(1+\log \left( 3 \sqrt{\frac{Rn_1}{s_1}} \right) \right) \right)^{\frac{1}{2}} \le \left( (a_2^2+a_1^2) R (n_2+1) \left(1+\log \left( 3 \sqrt{\frac{Rn_1}{s_1}} \right) \right) \right)^{\frac{1}{2}} \\
	&= \left( \frac{144 \Gamma^2 R^2}{m} \left( \frac{s_2 (n_2+1)}{n_2} + \frac{s_1 (n_2+1)}{n_1} \right) \left(1+\log \left( 3 \sqrt{\frac{Rn_1}{s_1}} \right) \right) \right)^{\frac{1}{2}} \\
	&\le \begin{cases}
	\left( \frac{432 \Gamma^2 R^2 (s_1+s_2)}{m} \left(1+\log \left( 3 \sqrt{\frac{Rn_1}{s_1}} \right) \right) \right)^\frac{1}{2} & n_1 \ge n_2, \\
	\left( \frac{432 \Gamma^2 R^2 s_2}{m} \left(1+\log \left( 3 \sqrt{\frac{Rn_1}{s_1}} \right) \right) \right)^\frac{1}{2} & \text{else.}
	\end{cases}
	\intertext{In the third and the last line we again used $s_1/n_1 \le s_2/n_2$. The fourth integral is similar to the second.}
	I_4 &= \sqrt{\frac{144 \Gamma^2 R^2 (s_1 + s_2)}{m}} \int_{a_2}^b \frac{1}{\eps} \sqrt{\log\left( \frac{9\sqrt{m}\eps n_1}{6\Gamma \sqrt{R} s_1} \right)} \de = \sqrt{\frac{144 \Gamma^2 R^2 (s_1 + s_2)}{m}} \left[ \frac{2}{3} \log^\frac{3}{2} \left( \frac{9\sqrt{m}\eps n_1}{6\Gamma \sqrt{R} s_1} \right) \right]_{\eps = a_2}^b \\
	&{= \left( \frac{64 \Gamma^2 R^2 (s_1 + s_2)}{m} \right)^{\frac{1}{2}} \left( \log^\frac{3}{2} \left( \frac{3n_1}{2s_1} \right) - \log^\frac{3}{2} \left( \frac{18\sqrt{s_2}n_1}{\sqrt{n_2}s_1} \right) \right) \le \left( \frac{64 \Gamma^2 R^2 (s_1 + s_2)}{m} \log^3 (18 n_1) \right)^{\frac{1}{2}}}
	\intertext{The last integral is similar to the third.}
	I_5 &\le \left( (b-a_2) R \left[ \eps \left(1 + \log\left( \frac{18\Gamma R}{\sqrt{m}\eps} \right)\right) \right]_{\eps = a_2}^b \right)^\frac{1}{2} \le \left( (b-a_2)^2 R \left( 1 + \log\left( 18 \sqrt{\frac{Rn_2}{s_2}} \right) \right) \right)^\frac{1}{2} \\
	&\le \left( (b^2+a_2^2) R \left( 1 + \log\left( 18 \sqrt{\frac{Rn_2}{s_2}} \right) \right) \right)^\frac{1}{2} = \left( \left( \frac{\Gamma^2 R^2}{m} + \frac{144 \Gamma^2 R^2 s_2}{m n_2} \right) \left( 1 + \log\left( 18 \sqrt{\frac{Rn_2}{s_2}} \right) \right) \right)^\frac{1}{2} \\
	&\le \left( \frac{145 \Gamma^2 R^2}{m} \left( 1 + \log\left( 18 \sqrt{\frac{Rn_2}{s_2}} \right) \right) \right)^\frac{1}{2}
	\end{align*}
	Let us now put all estimates together. If $s_1/n_1 \ge s_2/n_2$, the involved entities would just switch their roles. Hence, we obtain
	\begin{align*}
	\int_{0}^{\frac{\Gamma \sqrt{R}}{\sqrt{m}}} \sqrt{\log N(K,\Vert \cdot \Vert_F,\sqrt{m}\eps)} \;d\eps \le \sqrt{\frac{C_K \Gamma^2 R^2 (s_1 + s_2) \log^3(\max\{n_1,n_2\})}{m}}
	\end{align*}
	for some constant $C_K > 0$.
\end{Proof}

\subsection{Proof of \Cref{LocalConvergence}}

\purple{In this subsection we show the convergence of ATLAS to global minimizers as presented in \Cref{LocalConvergence}. To do so, we make use of the results from  \cite{attouch2010proximal}. In particular, we first present two technical  lemmas (Lemma \ref{Lemma 5} \& Lemma \ref{proposition6}), which are essentially generalizations of work  \cite{attouch2010proximal}. These lemmas would be useful to prove the central theorem (here \Cref{theorem8}) of Attouch et.\ al.\ in our general setting.  Finally,  the theorem on local convergence, \Cref{LocalConvergence}, can essentially be derived from \Cref{theorem8} and combines the two statements \cite[Theorem 9]{attouch2010proximal} and \cite[Theorem 10]{attouch2010proximal}. 
We refer the interested reader to \cite{attouch2010proximal} for further details. We also provide a reference to the original work in brackets.
}

\begin{lemma}[{\cite[Lemma 5]{attouch2010proximal}}] \label{Lemma 5}
	Under assumptions $(H)$ and $(H1)$ the sequences $u_k^1,\dots ,v_k^R$ are well-posed
	in the sense that all minimizations in (\ref{ProxAlgo0}) have unique and finite solutions.
	Moreover, \smallskip \\
	$(i)$
	\begin{align*}
	L(u_k^1,\dots ,v_k^R) + \sum_{r=1}^R \frac{1}{2\lambda_{k-1}^r} \Vert u_k^r -
	u_{k-1}^r \Vert_2^2 + \sum_{r=1}^R \frac{1}{2\mu_{k-1}^r} \Vert v_k^r - v_{k-1}^r
	\Vert_2^2 \leq L(u_{k-1}^1,\dots ,v_{k-1}^R),
	\end{align*}
	for all $k \geq 1$, hence $L(u_k^1,\dots ,v_k^R)$ is non-increasing. \smallskip \\
	$(ii)$
	\begin{align*}
	\sum_{k=1}^\infty \left( \Vert u_k^1 - u_{k-1}^1 \Vert_2^2 + \cdots +
	\Vert v_k^R - v_{k-1}^R \Vert_2^2 \right) < \infty,
	\end{align*}
	hence $\lim_{k \rightarrow \infty} \left( \Vert u_k^1 - u_{k-1}^1 \Vert_2 + \cdots +
	\Vert v_k^R - v_{k-1}^R \Vert_2 \right) = 0$. \smallskip \\
	$(iii)$ For $k \geq 1$, define
	\begin{align*}
	(\tilde{u}_k^1,\dots,\tilde{v}_k^R) :=
	\begin{pmatrix}
	\nabla_{u^1} Q(u_k^1,\dots,v_k^R) - \nabla_{u^1} Q(u_k^1,u_{k-1}^2,\dots,
	v_{k-1}^R) \\
	\vdots \\
	0
	\end{pmatrix} -
	\begin{pmatrix}
	\frac{1}{\lambda_{k-1}^1} (u_k^1 - u_{k-1}^1) \\
	\vdots \phantom{\frac{1}{2}} \\
	\frac{1}{\mu_{k-1}^R} (v_k^R - v_{k-1}^R)
	\end{pmatrix}.
	\end{align*}
	Then $(\tilde{u}_k^1,\dots,\tilde{v}_k^R) \in \partial L(u_k^1,\dots,v_k^R)$ and for all
	bounded subsequences $(u_{k'}^1,\dots,v_{k'}^R)$ we have $(\tilde{u}_{k'}^1,\dots,
	\tilde{v}_{k'}^R) \rightarrow 0$, hence $\dist (0,\partial L(u_{k'}^1,\dots,v_{k'}^R))
	\rightarrow 0$, for $k' \rightarrow \infty$.
\end{lemma}

\begin{Proof}
	From $\inf L > -\infty$ and $(H)$ it follows that the functions to be minimized in
	(\ref{ProxAlgo0}) are bounded below, coercive and lower semicontinuous and, therefore,
	the sequence $(u_k^1,\dots,v_k^R)$ is well-posed. 
	
	\paragraph{$(i)$} \purple{Using the} minimizing properties of $u_k^1,\dots,v_k^R$ from
	(\ref{ProxAlgo0}), \purple{we obtain} 
	\begin{align*}
	L(u_{k-1}^1,\dots,v_{k-1}^R) &\geq L(u_{k}^1,u_{k-1}^2,\dots,v_{k-1}^1, \dots, v_{k-1}^R) + \frac{1}{2\lambda_{k-1}^1}
	\Vert u_{k}^1 - u_{k-1}^1 \Vert_2^2 \\
	&\geq \left( L(u_{k}^1,u_{k-1}^2,\dots,u_{k-1}^R,v_{k}^1,v_{k-1}^2,\dots,v_{k-1}^R) +
	\frac{1}{2\mu_{k-1}^1} \Vert v_{k}^1 - v_{k-1}^1 \Vert_2^2 \right) + \frac{1}{2\lambda_{k-1}^1}
	\Vert u_{k}^1 - u_{k-1}^1 \Vert_2^2 \\
	&\phantom{\tab \tab \tab \tab \tab \tab} \vdots \\
	&\geq L(u_{k}^1,\dots,v_{k}^R) + \sum_{r=1}^R \frac{1}{2\lambda_{k-1}^r} \Vert 
	u_{k}^r - u_{k-1}^r \Vert_2^2 + \sum_{r=1}^R \frac{1}{2\mu_{k-1}^r} \Vert v_{k}^r - v_{k-1}^r
	\Vert_2^2.
	\end{align*}
	
	\paragraph{$(ii)$} From $(i)$ and $(H1)$ one has, for every $K \in \mathbb{N}$,
	\begin{align*}
	\frac{1}{2r_+} \sum_{k=1}^K \left( \Vert u_k^1 - u_{k-1}^1 \Vert_2^2 + \cdots +
	\Vert v_k^R - v_{k-1}^R \Vert_2^2 \right) &\leq
	\sum_{k=1}^K \left( L(u_{k-1}^1,\dots,v_{k-1}^R) - L(u_k^1,\dots,v_k^R) \right) \\ 
	&= L(u_0^1,\dots,v_0^R) - L(u_K^1,\dots,v_K^R) \\
	&< L(u_0^1,\dots,v_0^R) - \inf L < \infty.
	\end{align*}
	By letting $K \rightarrow \infty$ we get the claim.
	
	\paragraph{$(iii)$} By definition of $u_k^1$, $0$ must lie in the subdifferential of $\xi
	\mapsto L(\xi,u_{k-1}^2,\dots,v_{k-1}^R) + \frac{1}{2\lambda_{k-1}^1} \Vert \xi - u_{k-1}^1
	\Vert_2^2$ at $u_k^1$. As a similar fact holds true for the other sequences, one gets, 
	for all $1 \leq r \leq R$
	\begin{align*}
	0 &\in \frac{1}{\lambda_{k-1}^r} (u_k^r - u_{k-1}^r) + \partial_{u^r} L(u_k^1,\dots,
	u_k^{r},u_{k-1}^{r+1},\dots,u_{k-1}^R,v_k^1,\dots,v_k^{r-1},v_{k-1}^r,\dots,v_{k-1}^R), \\
	0 &\in \frac{1}{\mu_{k-1}^r} (v_k^r - v_{k-1}^r) + \partial_{v^r} L(u_k^1,\dots,
	u_k^{r},u_{k-1}^{r+1},\dots,u_{k-1}^R,v_k^1,\dots,v_k^{r},v_{k-1}^{r+1},\dots,
	v_{k-1}^R). 
	\end{align*}
	The structure of $L$ implies $\partial_{u^r} L(u_k^1,\dots,u_k^{r},u_{k-1}^{r+1}
	\dots,u_{k-1}^R,v_k^1,\dots,v_k^{r-1},v_{k-1}^r,\dots,v_{k-1}^R) = \partial f_r(u_k^r) +
	\nabla_{u^r} Q(u_k^1,\dots,u_k^{r},u_{k-1}^{r+1} \dots,u_{k-1}^R,v_k^1,
	\dots,v_k^{r-1},v_{k-1}^r,\dots,v_{k-1}^R)$ and a similar equation for the $v$-components.
	Hence, one may rewrite the inclusions above:
	\begin{align*}
	-\frac{1}{\lambda_{k-1}^1} (u_k^1 - u_{k-1}^1) - ( \nabla_{u^1} Q(u_k^1,u_{k-1}^2,\dots,
	v_{k-1}^R) - \nabla_{u^1} Q(u_k^1,\dots,v_k^R) ) &\in \partial f_1(u_k^1) + \nabla_{u^1}
	Q(u_k^1,\dots,v_k^R) , \\
	\vdots \\
	-\frac{1}{\mu_{k-1}^R} (v_k^R - v_{k-1}^R) &\in \partial g_R(v_k^R) + \nabla_{v^R}
	Q(u_k^1,\dots,v_k^R).
	\end{align*}
	This, together with Proposition 3 in the paper, yields the claim.
\end{Proof}
\begin{lemma}[{\cite[Proposition 6]{attouch2010proximal}}] \label{proposition6}
	Assume $(H)$ and $(H1)$ hold. Let $(u_k^1,\dots,v_k^R)$ be a sequence defined by
	(\ref{ProxAlgo0}) and $\omega (u_0^1,\dots,v_0^R)$ be \purple{a} (possibly empty) set of limit
	points. Then,
	\begin{enumerate}
		\item[(i)] if $(u_k^1,\dots,v_k^R)$ is bounded, \purple{then} $\omega (u_0^1,\dots,v_0^R)$ is nonempty,
		compact and connected\\ 
		and $\dist ((u_k^1,\dots,v_k^R),\omega (u_0^1,\dots,v_0^R))
		\rightarrow 0$ \purple{as $k \rightarrow \infty$},
		\item[(ii)] $\omega (u_0^1,\dots,v_0^R) \subset \crit L$, \purple{where $\crit L$ denotes a set of critical points of $L$,}
		\item[(iii)] $L$ is finite and constant on $\omega (u_0^1,\dots,v_0^R)$, equal to
		$\purple{\inf_{k \in \mathbb N}} L(u_k^1,\dots,v_k^R) = \purple{\lim_{k \rightarrow \infty}} L(u_k^1,\dots,v_k^R)$.
	\end{enumerate}
\end{lemma}

\begin{Proof}
	$(i)$ If $(u_k^1,\dots,v_k^R)$ is bounded, there exists a convergent subsequence, which
	implies $\omega (u_0^1,\dots,v_0^R)$ is nonempty. It also follows $\omega
	(u_0^1,\dots,v_0^R)$ is bounded. \\ 
	Let now $(\hat{u}^1,\dots,\hat{v}^R) \notin \omega (u_0^1,\dots,v_0^R)$ be given. There 
	must exist some $\varepsilon > 0$ with $(u_k^1,\dots,v_k^R) \notin B((\hat{u}^1,\dots,
	\hat{v}^R),\varepsilon)$, for all $k \in \mathbb{N}$. But then $\omega (u_0^1,\dots,v_0^R) 
	\cap B((\hat{u}^1,\dots,\hat{v}^R), \varepsilon) = \emptyset$. This proves $\omega (u_0^1,
	\dots,v_0^R)$ is closed and, hence, compact. \\
	{Let us assume $\omega (u_0^1,\dots,v_0^R)$ is not connected and let $\omega_c (u_0^1,\dots,v_0^R) \subset \omega (u_0^1,\dots,v_0^R)$ be a connected component. Then, $\omega (u_0^1,\dots,v_0^R) \setminus \omega_c (u_0^1,\dots,v_0^R) \neq \emptyset$ and there exists some
	$\varepsilon > 0$ such that 
	$$\omega_c^\eps (u_0^1,\dots,v_0^R) \cap \omega (u_0^1,\dots,v_0^R) \setminus \omega_c (u_0^1,\dots,v_0^R) = \emptyset,$$
	 where $\omega_c^\eps (u_0^1,\dots,v_0^R)$ is an $\eps$-neighborhood of $\omega_c (u_0^1,\dots,v_0^R)$. We know from \Cref{Lemma 5} $(ii)$ that
	$$\lim_{k \rightarrow \infty} \left( \Vert u_k^1 - u_{k-1}^1 \Vert_2 + \cdots + \Vert v_k^R -
	v_{k-1}^R \Vert_2 \right) = 0.$$
	Combined with $\omega_c (u_0^1,\dots,v_0^R)$ and $\omega (u_0^1,\dots,v_0^R) \setminus \omega_c (u_0^1,\dots,v_0^R)$ being
	sets of limit points of $(u_k^1,\dots,v_k^R),$ it implies the existence of a subsequence
	$(u_{k'}^1,\dots,v_{k'}^R) \subset \omega_c^\eps (u_0^1,\dots,v_0^R) 
	\setminus \omega_c^{\frac{\eps}{2}} (u_0^1,\dots,v_0^R)$. As this subsequence
	is bounded, it must have a limit point and $\omega (u_0^1,\dots,v_0^R) \cap
	\omega_c^\eps (u_0^1,\dots,v_0^R) 
	\setminus \omega_c^{\frac{\eps}{2}} (u_0^1,\dots,v_0^R) \neq \emptyset$. Contradiction.} \\
	The last part of $(i)$ can be proven in a similar way. If $\dist ((u_k^1,\dots,v_k^R),\omega
	(u_0^1,\dots,v_0^R)) \nrightarrow 0$, there must exist a subsequence that keeps distance
	to $\omega (u_0^1,\dots,v_0^R)$. But this subsequence again must have a limit point which
	obviously lies in $\omega (u_0^1,\dots,v_0^R)$. Contradiction.
	
	\paragraph{$(ii)$} We have, for all $k \geq 1$, $\xi^r \in \mathbb{R}^{n_1}$, $\eta^r \in
	\mathbb{R}^{n_2}$
	\begin{align*}
	L(u_k^1,u_{k-1}^2,\dots,v_{k-1}^R) + \frac{1}{2\lambda_{k-1}^1} \Vert u_k^1 - u_{k-1}^1
	\Vert_2^2 &\leq L(\xi^1,u_{k-1}^2,\dots,v_{k-1}^R) + \frac{1}{2\lambda_{k-1}^1} \Vert
	\xi^1 - u_{k-1}^1 \Vert_2^2 \\
	\vdots \\
	L(u_k^1,\dots,v_k^R) + \frac{1}{2\mu_{k-1}^R} \Vert v_k^R - v_{k-1}^R \Vert_2^2 &\leq 
	L(u_k^1,\dots,v_k^{R-1},\eta^R) + \frac{1}{2\mu_{k-1}^R} \Vert \eta^R - v_{k-1}^R
	\Vert_2^2
	\end{align*}
	Using the bounds on $\lambda_k^r$ and $\mu_k^r$ and the special form of $L$ one
	gets
	\begin{align*}
	f_1(u_k^1) + Q(u_k^1,u_{k-1}^2,\dots,v_{k-1}^R) + \frac{1}{2r_+} \Vert u_k^1 - u_{k-1}^1
	\Vert_2^2 &\leq f_1(\xi^1) + Q(\xi^1,u_{k-1}^2,\dots,v_{k-1}^R) + \frac{1}{2r_-} \Vert \xi^1
	- u_{k-1}^1 \Vert_2^2 \\
	\vdots \\
	g_R(v_k^R) + Q(u_k^1,\dots,v_k^R) + \frac{1}{2r_+} \Vert v_k^R - v_{k-1}^R \Vert_2^2 &\leq
	g_R(\eta^R) + Q(u_k^1,\dots,v_k^{R-1},\eta^R) + \frac{1}{2r_-} \Vert \eta^R - v_{k-1}^R
	\Vert_2^2
	\end{align*}
	Let $(\overline{u}^1,\dots,\overline{v}^R) \in \omega (u_0^1,\dots,v_0^R)$. There exists a
	subsequence $(u_{k'}^1,\dots,v_{k'}^R)$ of $(u_k^1,\dots,v_k^R)$ with $(u_{k'}^1
	\dots,v_{k'}^R) \rightarrow (\overline{u}^1,\dots,\overline{v}^R)$. Together with Lemma
	\ref{Lemma 5}.$(ii)$ this gives
	\begin{align*}
	\liminf_{k' \rightarrow \infty} f_r(u_{k'}^r) + Q(\overline{u}^1,\dots,\overline{v}^R) \leq
	f_r(\xi^r) +
	Q(\overline{u}^1,\dots,\xi^r,\dots,\overline{v}^R) + \frac{1}{2r_-} \Vert \xi^r - 
	\overline{u}^r \Vert_2^2,
	\end{align*}
	for all $1 \leq r \leq R$. We can now set $\xi^r = \overline{u}^r$ to obtain
	\begin{align*}
	\liminf_{k' \rightarrow \infty} f_r(u_{k'}^r) \leq f_r(\overline{u}^r).
	\end{align*}
	This and $f_r$ being lower semicontinuous yields
	\begin{align*}
	\lim_{k' \rightarrow \infty} f_r(u_{k'}^r) = f_r(\overline{u}^r).
	\end{align*}
	Repeating this for $g_r$, $1 \leq r \leq R$, and recalling the continuity of $Q$ we obtain
	$L(u_{k'}^1,\dots,v_{k'}^R) \rightarrow L(\overline{u}^1,\dots,\overline{v}^R)$. Combined with
	Lemma \ref{Lemma 5}.$(iii)$ and the closedness properties of $\partial L$(see Remark 1(b) {in
	\cite{attouch2010proximal}}) proves $0 \in \partial L(\overline{u}^1,\dots,\overline{v}^R)$.
	
	\paragraph{$(iii)$} \purple{As we just seen, for any point $(\overline{u}^1,\dots,\overline{v}^R) \in \omega (u_0^1,\dots,v_0^R),$ there exists a subsequence $(u_{k'}^1,\dots,v_{k'}^R)$ of $(u_k^1,\dots,v_k^R)$  with $L(u_{k'}^1
	\dots,v_{k'}^R) \rightarrow L(\overline{u}^1,\dots,\overline{v}^R).$ Then} $L(\overline{u}^1,\dots,\overline{v}^R) = \inf L(u_k^1,\dots,v_k^R)$ as
	$L(u_k^1,\dots,v_k^R)$ is non-increasing. This holds for every limit point. Hence, $L$ is finite and constant on the set of limit points.
\end{Proof}

\paragraph{} As in \cite{attouch2010proximal} we use the notation
\begin{align*}
z_k := (u_k^1,\dots,v_k^R), \tab \tab &l_k := L(z_k), \\
\overline{z} := (\overline{u}^1,\dots, \overline{v}^R), \tab \tab &\overline{l} := L(\overline{z}).
\end{align*}
\purple{The next theorem essentially says that a sequence $z_k$ that starts in the neighborhood of a point $\overline{z}$ as described in \eqref{rho} and that does not improve $L(\overline{z})$ as given in \eqref{eta} converges to a critical point near $\overline{z}.$}
\begin{theorem}[{\cite[Theorem 8]{attouch2010proximal}}] \label{theorem8}
	Let $L$ satisfy $(H)$, $(H1)$ and have the KL-property at some $\overline{z}$. Denote by
	$U$, $\eta$ and $\varphi : \left[ 0,\eta \right) \rightarrow \mathbb{R}$ the objects
	connected to the KL-property of $L$ at $\overline{z}$. Let $\rho > 0$ be chosen such that
	$B(\overline{z}, \rho) \subset U$. 
	Let $z_k$ be generated by (\ref{ProxAlgo0}) with $z_0$ as initial point. Let us assume that
	\begin{align} \label{eta}
	\overline{l} < l_k < \overline{l} + \eta,
	\end{align}
	for all $k \geq 0$, and
	\begin{align} \label{rho}
	M \varphi(l_0-\overline{l}) + 2 \sqrt{2r_+} \sqrt{l_0 - \overline{l}} + \Vert z_0 - 
	\overline{z} \Vert_2 < \rho
	\end{align}
	where $M = 2r_+ (C{\sqrt{2R}} + \frac{1}{r_-})$ and $C$ is a Lipschitz-constant for $\nabla Q$ on
	$B(\overline{z}, \sqrt{2R} \rho)$.Then, the sequence $z_k$ converges to a critical point of $L$
	and the following holds, for all $k \geq 0$:
	\begin{enumerate}
		\item[$(i)$] $z_k \in B(\overline{z}, \rho)$ 
		\item[$(ii)$] $\sum_{i = k+1}^\infty \Vert z_{i+1} - z_i \Vert_2 \leq M \varphi(l_k - 
		\overline{l}) + \sqrt{2r_+} \sqrt{l_k-\overline{l}}$.
	\end{enumerate}
\end{theorem}

\begin{Proof}
	We may without loss of generality assume $L(\overline{z}) = 0$ (just replace $L$ by $L - L(
	\overline{z})$ ). With Lemma \ref{Lemma 5}.$(i)$ we have
	\begin{align} \label{19}
	l_i - l_{i+1} \geq \frac{1}{2r_+} \Vert z_{i+1} - z_i \Vert_2^2,
	\end{align}
	for all $i \geq 0$. Moreover, $\varphi'(l_i)$ makes sense in view of (\ref{eta}) and
	$\varphi'(l_i) > 0$. Hence,
	\begin{align*}
	\varphi'(l_i) (l_i - l_{i+1}) \geq \frac{\varphi'(l_i)}{2r_+} \Vert z_{i+1} - z_i \Vert_2^2.
	\end{align*}
	Owing to $\varphi$ being concave, we obtain
	\begin{align} \label{20}
	\varphi(l_i) - \varphi(l_{i+1}) \geq \frac{\varphi'(l_i)}{2r_+} \Vert z_{i+1} - z_i \Vert_2^2,
	\end{align}
	for all $i \geq 0$. Let us first check $(i)$ for $k = 0$ and $k = 1$. We know from (\ref{rho})
	that $z_0$ lies in $B(\overline{z},\rho)$. Furthermore, (\ref{19}) yields
	\begin{align*}
	\frac{1}{2r_+} \Vert z_1 - z_0 \Vert_2^2 \leq l_0 - l_1 \leq l_0
	\end{align*}
	which gives
	\begin{align*}
	\Vert z_1 - \overline{z} \Vert_2 \leq \Vert z_1 - z_0 \Vert_2 + \Vert z_0 - \overline{z} \Vert_2
	\leq \sqrt{2r_+} \sqrt{l_0} + \Vert z_0 - \overline{z} \Vert_2 < \rho.
	\end{align*}
	
	\paragraph{} Let us now prove by induction that $z_k \in B(\overline{z},\rho)$, for all 
	$k \geq 0$. We assume this holds true up to some $k \geq 0$. Hence, for $0 \leq i \leq k$,
	using $z_i \in B(\overline{z},\rho)$ and $0 < l_i < \eta$ we can write the KL-inequality
	\begin{align*}
	\varphi'(l_i) \dist(0,\partial L(z_i)) \geq 1.
	\end{align*}
	Lemma \ref{Lemma 5}.$(iii)$ says
	\begin{align*}
	z_i^\ast := 
	\begin{pmatrix}
	\nabla_{u^1} Q(u_i^1,\dots,v_i^R) - \nabla_{u^1} Q(u_i^1,u_{i-1}^2,\dots,
	v_{i-1}^R) \\
	\vdots \\
	0
	\end{pmatrix} -
	\begin{pmatrix}
	\frac{1}{\lambda_{i-1}^1} (u_i^1 - u_{i-1}^1) \\
	\vdots \phantom{\frac{1}{2}} \\
	\frac{1}{\mu_{i-1}^R} (v_i^R - v_{i-1}^R)
	\end{pmatrix}.
	\end{align*}	
	is an element of $\partial L(z_i)$. So, we have
	\begin{align} \label{KL}
	\varphi'(l_i) \Vert z_i^\ast \Vert_2 \geq 1,
	\end{align}
	for all $1 \leq i \leq k$. \purple{Let us now examine}  $\Vert z_i^\ast \Vert_2$, for $1 \leq i \leq k$. On the one
	hand,
	\begin{align*}
	\left\Vert \left( \frac{1}{\lambda_{i-1}^1} (u_i^1 - u_{i-1}^1), \dots ,\frac{1}{\mu_{i-1}^R}
	(v_i^R - v_{i-1}^R) \right) \right\Vert_2 \leq \frac{1}{r_-} \Vert z_i - z_{i-1} \Vert_2.
	\end{align*}
	On the other hand, for arbitrary $s_t \in \{ i-1,i \}$, $t \in \{ 1,\dots,2R \}$,
	\begin{align*}
	\Vert (u_{s_1}^1,\dots,v_{s_{2R}}^R) - (\overline{u}^1,\dots,\overline{v}^R) \Vert_2^2 &= 
	\Vert u_{s_1}^1 - \overline{u}^1 \Vert_2^2 + \cdots + \Vert v_{s_{2R}}^R - \overline{v}^R
	\Vert_2^2 \\
	&\leq \Vert z_{s_1} - \overline{z} \Vert_2^2 + \cdots + \Vert z_{s_{2R}} - \overline{z} \Vert_2^2
	\leq 2R\rho^2.
	\end{align*}
	Hence, $(u_{s_1}^1,\dots,v_{s_{2R}}^R)$ and $z_i$ lie in $B(\overline{z}, \sqrt{2R} \rho)$. We 
	can use Lipschitz-continuity of $\nabla Q$ to obtain
	\begin{align*}
	\Vert \nabla_\theta Q(u_{s_1}^1,\dots,v_{s_{2R}}^R) - \nabla_\theta Q(u_i^1,\dots,v_i^R)
	\Vert_2 \leq C \Vert z_i - z_{i-1} \Vert_2,
	\end{align*}
	for any $\theta \in \{ u^1,\dots,v^R \}$, which implies
	\begin{align*}
	\left\Vert
	\begin{pmatrix}
	\nabla_{u^1} Q(u_i^1,\dots,v_i^R) - \nabla_{u^1} Q(u_i^1,u_{i-1}^2,\dots,
	v_{i-1}^R) \\
	\vdots \\
	0
	\end{pmatrix}
	\right\Vert_2
	\leq C{\sqrt{2R}} \Vert z_i - z_{i-1} \Vert_2.
	\end{align*}
	We get
	\begin{align*}
	\Vert z_i^\ast \Vert_2 \leq (C{\sqrt{2R}} + \frac{1}{r_-}) \Vert z_i - z_{i-1} \Vert_2,
	\end{align*}
	for all $1 \leq i \leq k$. Now (\ref{KL}) yields
	\begin{align*}
	\varphi'(l_i) \geq \frac{1}{C{\sqrt{2R}} + \frac{1}{r_-}} \Vert z_i - z_{i-1} \Vert_2^{-1}, \tab 
	1 \leq i \leq k,
	\end{align*}
	and combined with (\ref{20})
	\begin{align*}
	\varphi(l_i) - \varphi(l_{i+1}) \geq \frac{1}{M} \frac{\Vert z_{i+1} - z_i \Vert_2^2}{\Vert z_i
		- z_{i-1} \Vert_2}, \tab 1 \leq i \leq k.
	\end{align*}
	This is equivalent to
	\begin{align*}
	\Vert z_i - z_{i-1} \Vert_2^\frac{1}{2} (M(\varphi(l_i) - \varphi(l_{i+1})))^\frac{1}{2}
	\geq \Vert z_{i+1} - z_i \Vert_2
	\end{align*}
	and, using $ab \leq (a^2 + b^2)/2$, gives
	\begin{align} \label{phi}
	\Vert z_i - z_{i-1} \Vert_2 + M(\varphi(l_i) - \varphi(l_{i+1}))
	\geq 2\Vert z_{i+1} - z_i \Vert_2, \tab 1 \leq i \leq k.
	\end{align}
	Summation over $i$ leads to
	\begin{align*}
	\Vert z_1 - z_0 \Vert_2 + M(\varphi(l_1) - \varphi(l_{k+1})) \geq 
	\sum_{i=1}^k \Vert z_{i+1} - z_i \Vert_2 + \Vert z_{k+1} - z_k \Vert_2.
	\end{align*}
	Therefore, by using the monotonicity properties of $\varphi$ and $l_k$
	\begin{align*}
	\Vert z_1 - z_0 \Vert_2 + M\varphi(l_0) \geq \sum_{i = 1}^k \Vert z_{i+1} - z_i \Vert_2.
	\end{align*}
	Finally,
	\begin{align*}
	\Vert z_{k+1} - \overline{z} \Vert_2 \leq \sum_{i = 1}^k \Vert z_{i+1} - z_i \Vert_2
	+ \Vert z_1 - \overline{z} \Vert_2 \leq M\varphi(l_0) + 2 \sqrt{2r_+} \sqrt{l_0} +
	\Vert z_0 - \overline{z} \Vert_2 < \rho
	\end{align*}
	which closes the induction and proves $(i)$. Moreover, (\ref{phi}) holds for all $i \geq 1$. We
	can sum from $k$ to $K$ and get
	\begin{align*}
	\Vert z_k - z_{k-1} \Vert_2 + M(\varphi(l_k) - \varphi(l_{K+1})) \geq \sum_{i = k}^K
	\Vert z_{i+1} - z_i \Vert_2 + \Vert z_{K+1} - z_K \Vert_2.
	\end{align*}
	For $K \rightarrow \infty$, this becomes
	\begin{align*}
	\sum_{i = k}^\infty \Vert z_{i+1} - z_i \Vert_2 \leq \Vert z_k - z_{k-1} \Vert_2 + M\varphi(l_k).
	\end{align*}
	We conclude with (\ref{19}) proving $(ii)$
	\begin{align*}
	\sum_{i = k}^\infty \Vert z_{i+1} - z_i \Vert_2 \leq M\varphi(l_k) + \sqrt{2r_+} \sqrt{l_{k-1}}
	\leq M\varphi(l_{k-1}) + \sqrt{2r_+} \sqrt{l_{k-1}}.
	\end{align*}
	This implies $z_k$ is convergent and, therefore, its limit is a critical point. This was
	guaranteed by Lemma \ref{proposition6}.
	
\end{Proof}


\section*{Acknowledgments}

MF and JM acknowledge the support of the DFG project “Information Theory and Recovery Algorithms for Quantized and Distributed Compressed Sensing”. MF and VN acknowledge the support of the DFG-FWF project “Multipenalty Regularization for High-Dimensional Learning”. VN acknowledges the support of the project No 251149/O70 “Function-driven Data Learning in High Dimension” (FunDaHD) funded by the Research Council of Norway. The authors thank Dominik Stöger for providing the very helpful counterexample in Remark \ref{DominiksRemark}.

\bibliographystyle{IEEEtranS}
\bibliography{mybibM,CSBib}

\begin{thebibliography}{10}
\providecommand{\url}[1]{#1}
\csname url@samestyle\endcsname
\providecommand{\newblock}{\relax}
\providecommand{\bibinfo}[2]{#2}
\providecommand{\BIBentrySTDinterwordspacing}{\spaceskip=0pt\relax}
\providecommand{\BIBentryALTinterwordstretchfactor}{4}
\providecommand{\BIBentryALTinterwordspacing}{\spaceskip=\fontdimen2\font plus
\BIBentryALTinterwordstretchfactor\fontdimen3\font minus
  \fontdimen4\font\relax}
\providecommand{\BIBforeignlanguage}[2]{{%
\expandafter\ifx\csname l@#1\endcsname\relax
\typeout{** WARNING: IEEEtranS.bst: No hyphenation pattern has been}%
\typeout{** loaded for the language `#1'. Using the pattern for}%
\typeout{** the default language instead.}%
\else
\language=\csname l@#1\endcsname
\fi
#2}}
\providecommand{\BIBdecl}{\relax}
\BIBdecl

\bibitem{ahmed2014blind}
A.~Ahmed, B.~Recht, and J.~Romberg, ``Blind deconvolution using convex
  programming,'' \emph{IEEE Transactions on Information Theory}, vol.~60,
  no.~3, pp. 1711--1732, 2014.

\bibitem{attouch2010proximal}
H.~Attouch, J.~Bolte, P.~Redont, and A.~Soubeyran, ``Proximal alternating
  minimization and projection methods for nonconvex problems: An approach based
  on the kurdyka-{\l}ojasiewicz inequality,'' \emph{Mathematics of Operations
  Research}, vol.~35, no.~2, pp. 438--457, 2010.

\bibitem{badea2005sparse}
L.~Badea and D.~Tilivea, ``Sparse factorizations of gene expression data guided
  by binding data,'' in \emph{Biocomputing 2005}.\hskip 1em plus 0.5em minus
  0.4em\relax World Scientific, 2005, pp. 447--458.

\bibitem{bahmani2016near}
S.~Bahmani and J.~Romberg, ``Near-optimal estimation of simultaneously sparse
  and low-rank matrices from nested linear measurements,'' \emph{Information
  and Inference: A Journal of the IMA}, vol.~5, no.~3, pp. 331--351, 2016.

\bibitem{bainbridge2013intrinsic}
W.~A. Bainbridge, P.~Isola, and A.~Oliva, ``The intrinsic memorability of face
  photographs.'' \emph{Journal of Experimental Psychology: General}, vol. 142,
  no.~4, p. 1323, 2013.

\bibitem{Bennett07thenetflix}
J.~Bennett, S.~Lanning, and N.~Netflix, ``The netflix prize,'' in \emph{In KDD
  Cup and Workshop in conjunction with KDD}, 2007.

\bibitem{BLUMENSATH2009265}
T.~Blumensath and M.~E. Davies, ``Iterative hard thresholding for compressed
  sensing,'' \emph{Applied and Computational Harmonic Analysis}, vol.~27,
  no.~3, pp. 265 -- 274, 2009.

\bibitem{candes2011tight}
E.~J. Candes and Y.~Plan, ``Tight oracle inequalities for low-rank matrix
  recovery from a minimal number of noisy random measurements,'' \emph{IEEE
  Transactions on Information Theory}, vol.~57, no.~4, pp. 2342--2359, 2011.

\bibitem{Cand2009}
E.~J. Cand{\`e}s and B.~Recht, ``Exact matrix completion via convex
  optimization,'' \emph{Foundations of Computational Mathematics}, vol.~9,
  no.~6, 2009.

\bibitem{d2005direct}
A.~d'Aspremont, L.~E. Ghaoui, M.~I. Jordan, and G.~R. Lanckriet, ``A direct
  formulation for sparse {PCA} using semidefinite programming,'' in
  \emph{Advances in neural information processing systems}, 2005, pp. 41--48.

\bibitem{daubechies2016}
I.~Daubechies, M.~Defrise, and C.~De~Mol, ``Sparsity-enforcing regularisation
  and {ISTA} revisited,'' \emph{Inverse Problems}, vol.~32, no.~10, p. 104001,
  2016.

\bibitem{dili01}
Z.~Ding and Y.~Li, \emph{Blind equalization and identification}.\hskip 1em plus
  0.5em minus 0.4em\relax CRC press, 2001.

\bibitem{gekrst17}
J.~A. Geppert, F.~Krahmer, and D.~Stöger, ``Refined performance guarantees for
  sparse power factorization,'' in \emph{2017 International Conference on
  Sampling Theory and Applications (SampTA)}, July 2017, pp. 509--513.

\bibitem{godard1980self}
D.~Godard, ``Self-recovering equalization and carrier tracking in
  two-dimensional data communication systems,'' \emph{IEEE transactions on
  communications}, vol.~28, no.~11, pp. 1867--1875, 1980.

\bibitem{grklna17}
M.~Grasmair, T.~Klock, and V.~Naumova, ``Adaptive multi-penalty regularization
  based on a generalized lasso path,'' \emph{to appear in Applied and
  Computational Harmonic Analysis}, 2018.

\bibitem{naumova2016}
M.~Grasmair and V.~Naumova, ``Conditions on optimal support recovery in
  unmixing problems by means of multi-penalty regularization,'' \emph{Inverse
  Problems}, vol.~32, no.~10, p. 104007, 2016.

\bibitem{haykin1994blind}
S.~Haykin, ``The blind deconvolution problem,'' \emph{Blind Deconvolution},
  p.~1, 1994.

\bibitem{jagannathan2003risk}
R.~Jagannathan and T.~Ma, ``Risk reduction in large portfolios: Why imposing
  the wrong constraints helps,'' \emph{The Journal of Finance}, vol.~58, no.~4,
  pp. 1651--1683, 2003.

\bibitem{jain2013low}
P.~Jain, P.~Netrapalli, and S.~Sanghavi, ``Low-rank matrix completion using
  alternating minimization,'' \emph{Proceedings of the forty-fifth annual ACM
  symposium on Theory of computing}, pp. 665--674, 2013.

\bibitem{jolliffe2011principal}
I.~Jolliffe, ``Principal component analysis,'' in \emph{International
  encyclopedia of statistical science}.\hskip 1em plus 0.5em minus 0.4em\relax
  Springer, 2011, pp. 1094--1096.

\bibitem{JungKrahmerStoeger2017}
P.~Jung, F.~Krahmer, and D.~St\"oger, ``{Blind Demixing and Deconvolution at
  Near-Optimal Rate},'' \emph{ArXiv: 1704.04178}, 2017.

\bibitem{krahmer2014suprema}
F.~Krahmer, S.~Mendelson, and H.~Rauhut, ``Suprema of chaos processes and the
  restricted isometry property,'' \emph{Communications on Pure and Applied
  Mathematics}, vol.~67, no.~11, pp. 1877--1904, 2014.

\bibitem{lee1999learning}
D.~D. Lee and H.~S. Seung, ``Learning the parts of objects by non-negative
  matrix factorization,'' \emph{Nature}, vol. 401, no. 6755, p. 788, 1999.

\bibitem{lee2016blind}
K.~Lee, Y.~Li, M.~Junge, and Y.~Bresler, ``Blind recovery of sparse signals
  from subsampled convolution,'' \emph{IEEE Transactions on Information
  Theory}, vol.~63, no.~2, pp. 802--821, 2016.

\bibitem{lee2013near}
K.~Lee, Y.~Wu, and Y.~Bresler, ``Near-optimal compressed sensing of a class of
  sparse low-rank matrices via sparse power factorization,'' \emph{IEEE
  Transactions on Information Theory}, vol.~64, no.~3, pp. 1666--1698, 2018.

\bibitem{li2013global}
G.~Li, ``Global error bounds for piecewise convex polynomials,''
  \emph{Mathematical Programming}, vol. 137, no. 1-2, pp. 37--64, 2013.

\bibitem{li2018rapid}
X.~Li, S.~Ling, T.~Strohmer, and K.~Wei, ``Rapid, robust, and reliable blind
  deconvolution via nonconvex optimization,'' \emph{Applied and computational
  harmonic analysis}, 2018.

\bibitem{ling2017blind}
S.~Ling and T.~Strohmer, ``Blind deconvolution meets blind demixing: Algorithms
  and performance bounds,'' \emph{IEEE Transactions on Information Theory},
  vol.~63, no.~7, pp. 4497--4520, 2017.

\bibitem{ling2017regularized}
------, ``Regularized gradient descent: a non-convex recipe for fast joint
  blind deconvolution and demixing,'' \emph{Information and Inference: A
  Journal of the IMA}, 2017.

\bibitem{ma2017implicit}
C.~Ma, K.~Wang, Y.~Chi, and Y.~Chen, ``Implicit regularization in nonconvex
  statistical estimation: Gradient descent converges linearly for phase
  retrieval, matrix completion and blind deconvolution,'' \emph{arXiv preprint
  arXiv:1711.10467}, 2017.

\bibitem{naumova2014minimization}
V.~Naumova and S.~Peter, ``Minimization of multi-penalty functionals by
  alternating iterative thresholding and optimal parameter choices,''
  \emph{Inverse Problems}, vol.~30, no.~12, p. 125003, 2014.

\bibitem{oymak2015simultaneously}
S.~Oymak, A.~Jalali, M.~Fazel, Y.~C. Eldar, and B.~Hassibi, ``Simultaneously
  structured models with application to sparse and low-rank matrices,''
  \emph{IEEE Transactions on Information Theory}, vol.~61, no.~5, pp.
  2886--2908, 2015.

\bibitem{Plan2013LP}
Y.~Plan and R.~Vershynin, ``One-bit compressed sensing by linear programming,''
  \emph{Communications on Pure and Applied Mathematics}, vol.~66, no.~8, pp.
  1275--1297, 2013.

\bibitem{Plan2014}
------, ``Dimension reduction by random hyperplane tessellations,''
  \emph{Discrete {\&} Computational Geometry}, vol.~51, no.~2, pp. 438--461,
  2014.

\bibitem{recht2010}
B.~Recht, M.~Fazel, and P.~A. Parrilo, ``Guaranteed minimum-rank solutions of
  linear matrix equations via nuclear norm minimization,'' \emph{SIAM Review},
  vol.~52, no.~3, pp. 471--501, 2010.

\bibitem{Recht:2010fk}
------, ``Guaranteed minimum-rank solutions of linear matrix equations via
  nuclear norm minimization,'' \emph{SIAM Rev.}, vol.~52, no.~3, pp. 471--501,
  2010.

\bibitem{stockham1975blind}
T.~G. Stockham, T.~M. Cannon, and R.~B. Ingebretsen, ``Blind deconvolution
  through digital signal processing,'' \emph{Proceedings of the IEEE}, vol.~63,
  no.~4, pp. 678--692, 1975.

\bibitem{vershynin2010introduction}
R.~Vershynin, ``Introduction to the non-asymptotic analysis of random
  matrices,'' in \emph{Compressed Sensing: Theory and Applications}.\hskip 1em
  plus 0.5em minus 0.4em\relax Cambridge Univ. Press, 2012, pp. 210--268.

\bibitem{zass2007nonnegative}
R.~Zass and A.~Shashua, ``Nonnegative sparse {PCA},'' in \emph{Advances in
  neural information processing systems}, 2007, pp. 1561--1568.

\bibitem{zou2006sparse}
H.~Zou, T.~Hastie, and R.~Tibshirani, ``Sparse principal component analysis,''
  \emph{Journal of computational and graphical statistics}, vol.~15, no.~2, pp.
  265--286, 2006.

\end{thebibliography}
\end{document}